  \documentstyle[12pt]{article}
\newtheorem{ournothing}{}[subsection]

\newtheorem{ourlemma}[ournothing]{Lemma}
\newtheorem{ourtheorem}[ournothing]{Theorem}

\newtheorem{ourcorollary}[ournothing]{Corollary}
\newtheorem{ourproposition}[ournothing]{Proposition}
\newtheorem{ourconjecture}[ournothing]{Conjecture}

\newcommand{\numero}[1]{
\addtocounter{section}{1}
\begin{center}{\bf \thesection .\
#1\vspace{-.1in}}\end{center}
\setcounter{subsection}{0}
\setcounter{subsubsection}{0}
\setcounter{ournothing}{0}\indent}

\newcommand{\subnumero}[1]{
\addtocounter{subsection}{1}
\noindent{\bf \thesubsection \ }
{\sc #1}
\setcounter{subsubsection}{0}
\setcounter{ournothing}{0}\indent}

\newcommand{\subsubnumero}[1]{
\addtocounter{subsubsection}{1}
\addtocounter{ournothing}{1}
\noindent
\hspace*{.14cm}{\bf \thesubsubsection \ }
{\em #1}}

\newenvironment{parag}{
\addtocounter{subsubsection}{1}
\begin{ournothing}\rm }{
\end{ournothing}}

\newcommand{\eop}{\hfill $/$\hspace*{-.1cm}$/$\hspace*{-.1cm}$/$\vspace{.1in}}

\newenvironment{lemma}{
\addtocounter{subsubsection}{1}
\begin{ourlemma}}{
\end{ourlemma}}

\newenvironment{corollary}{
\addtocounter{subsubsection}{1}
\begin{ourcorollary}}{
\end{ourcorollary}}

\newenvironment{theorem}{
\addtocounter{subsubsection}{1}
\begin{ourtheorem}}{
\end{ourtheorem}}

\newenvironment{proposition}{
\addtocounter{subsubsection}{1}
\begin{ourproposition}}{
\end{ourproposition}}

\newenvironment{definition}{
\addtocounter{subsubsection}{1}
\begin{ournothing}\hspace*{-.2cm}---{\bf Definition:}\, \rm }{
\end{ournothing}}

\newenvironment{conjecture}{
\addtocounter{subsubsection}{1}
\begin{ourconjecture}}{
\end{ourconjecture}}

\newenvironment{remark}{
\addtocounter{subsubsection}{1}
\begin{ournothing}\hspace*{-.2cm}---Remark: \rm }{
\end{ournothing}}

\newenvironment{Theorem}{\vspace{.05in} \par
\noindent {\bf Theorem} \hspace*{.05in} \em }{
\vspace{.05in}\newline }

\newcommand{\cc}{{\bf C}}
\newcommand{\rr}{{\bf R}}
\newcommand{\qq}{{\bf Q}}

\newcommand{\Cc}{{\cal C}}

\newcommand{\Bb}{{\cal B}}

\newcommand{\Pp}{{\cal P}}

\newcommand{\Xx}{{\cal X}}

\newcommand{\tworightarrows}{\stackrel{\displaystyle \rightarrow}{\rightarrow}}

\begin{document}

\section*{Limits in $n$-categories}

Carlos Simpson\newline
CNRS, UMR 5580, Universit\'e Paul Sabatier, 31062 Toulouse CEDEX, France.

\bigskip

\numero{Introduction}

One of the main notions in category theory is the notion of limit. Similarly,
one of the most commonly used techniques in homotopy theory is the notion of
``homotopy limit'' commonly called ``holim'' for short.  The purpose of the
this paper is to begin to develop the notion of limit for $n$-categories,
which should be a bridge between the categorical notion of limit
and the homotopical notion of holim.

We treat Tamsamani's notion of $n$-category \cite{Tamsamani}, but similar
arguments and results should hold for the Baez-Dolan approach
\cite{BaezDolanLetter}, \cite{BaezDolanIII}, or the Batanin approach
\cite{Batanin}, \cite{Batanin2}.

We define the notions of direct and inverse limits in an arbitrary (fibrant cf
\cite{nCAT}) $n$-category $C$.  Suppose $A$ is an $n$-category, and suppose
$\varphi : A\rightarrow C$ is a morphism, which we think of as a family of
objects of $C$ indexed by $A$. For any object $U\in C$ we can define the
$(n-1)$-category $Hom(\varphi , U)$ of morphisms from $\varphi$ to $U$.
We say that a morphism $\epsilon : \varphi \rightarrow U$ (i.e. an object
of this
$(n-1)$-category) is a {\em direct limit of $\varphi$} (cf \ref{defdirect}
below)
if, for every other object $V\in C$ the (weakly defined) composition with
$\epsilon$ induces an equivalence of $n-1$-categories from $Hom _C(U,V)$ to $Hom
(\varphi , V)$.

An analogous definition holds for saying that a morphism $U\rightarrow
\varphi$ is an {\em inverse limit of $\varphi$} (cf \ref{definverse} below).

The
main theorems concern the case where $C$ is the $n+1$-category $nCAT'$ (fibrant
replacement of that of) of $n$-categories.

\begin{Theorem}
{\rm (\ref{inverse} \ref{direct})}
The $n+1$-category $nCAT'$
admits
arbitrary inverse and direct limits.
\end{Theorem}

This is the analogue of the
classical statement that
the category $Sets$ admits inverse and direct limits---which is the case $n=0$
of our theorem.

The fact that we work in an $n$-category means that we automatically keep track
of ``higher homotopies'' and the like. This brings the ideas much closer to the
relatively simple notion of limits in a category.

I first learned of the notion of ``2-limit'' from the paper of Deligne and
Mumford \cite{DeligneMumford}, where it appears at the beginning with very
little explanation.   Unfortunately at the writing of the present paper I
have not been able to investigate the history of the notion of $n$-limits,
and I apologize in advance for any references left out.

At the end of the paper we propose many applications of the notion of limit.
Most of these are as of yet in an embryonic stage of development and we
don't pretend to give complete proofs.

\begin{center}
{\em Organization}
\end{center}
The paper is organized as follows: we start in \S 2
with some preliminary remarks
recalling  the notion of $n$-category from \cite{Tamsamani} and the closed
model structure from \cite{nCAT}. At the end of \S 2 we define and discuss one
of our main technical tools, the $n+1$-precat $\Upsilon ^k(E_1,\ldots , E_k)$
which can be seen as a $k$-simplex with $n$-precats $E_1,\ldots , E_k$
attached to the principal edges.

In \S 3 we give the basic definitions of
inverse and direct limits, and treat some general properties such as
invariance under equivalence, and variation with parameters.

In \S 4 we start into the main
result of the paper which is the existence of inverse limits in $nCAT'$.
Here, the construction is relatively straightforward: if $\varphi :
A\rightarrow nCAT'$ is a morphism then the inverse limit of $\varphi$ is just
the $n$-category
$$
\lambda = Hom _{\underline{Hom}(A,C)}(\ast , \varphi ).
$$
This is in perfect accord with the usual situation for inverse limits of
families of sets. Our only problem is to prove that this satisfies the
definition of being an inverse limit. Thus the reader could read up to here and
then skip the proof and move on to direct limits.

We treat direct limits in $nCAT'$ by a trick in \S 5: given $\psi : A
\rightarrow
nCAT'$ we construct an $n+1$-category $D$ parametrizing all morphisms $\psi
\rightarrow B$ to objects of $nCAT'$, and then construct the direct limit $U$
as the inverse limit of the functor $D\rightarrow nCAT'$. The main problem here
is that, because of set-theoretic considerations, we must restrict to a
category $D_{\alpha}$ of morphisms to objects $B$ with cardinality bounded by
$\alpha$.  We mimic a possible construction of direct limits in $Sets$
and encounter a few of the same difficulties as with inverse limits.
Again, the reader might want to just look at the proof for $Sets$
and skip the difficulties encountered in extending this to $nCAT'$.

At the end we discuss some proposed applications:
\newline
---First, the notions of homotopy coproduct and fiber product, and
their relation to the usual notions which can be calculated using the
closed model structure.
\newline
---Then we discuss representable functors, give a conjectural criterion
for when a functor should be representable, and apply it to
the problem of finding internal $\underline{Hom}$.
\newline
---The next subsection concerns $n$-stacks, defined using certain inverse
limits.
\newline
---We give a very general discussion of the notion of stack in any setting
where one knows what limits mean.
\newline
---We discuss direct images of families of $n$-categories by functors of the
underlying $n+1$-categories, and apply this to give a notion of ``realization''.
\newline
---Finally, we use limits to propose a notion of {\em relative Malcev
completion of the higher homotopy type}.

In all of the above applications except the first, most of the statements which
we need are left as conjectures. Thus, this discussion of applications is
still only at a highly speculative stage. One recurring theme is that the
argument given in \S 5 should work in a fairly general range of situations.

I would like to thank A. Hirschowitz, for numerous discussions about stacks
which contributed to the development of the ideas in this paper. I would like
to thank J. Tapia and J. Pradines for a helpful discussion concerning the
argument in \S 5.

\bigskip

\numero{Preliminary remarks}

\subnumero{$n$-categories}

\begin{parag}
\label{catnerve}
We begin by recalling the correspondence between categories and their nerves.
Let $\Delta$ denote the simplicial category whose objects are finite ordered
sets $p= \{ 0,\ldots , p\}$ and morphisms are order-preserving maps.
If $C$ is a category then its nerve is the simplicial set (i.e.
a functor $A:\Delta ^o\rightarrow Sets$) defined by setting $A_p$ equal to the
set of composable $p$-uples of arrows in $C$.  This satisfies the property that
the ``Segal maps'' (cf the discussion of Segal's delooping machine \cite{Segal}
in \cite{Adams} for the origin of this terminology)
$$
A_p \rightarrow A_1 \times _{A_0} \ldots \times _{A_0} A_1
$$
are isomorphisms. To be precise this  map is given by the $p$-uple
of face maps $1\rightarrow p$ which take $0$ to $i$ and $1$ to $i+1$
for $i=0,\ldots , p-1$. Conversely, given a simplicial set $A$ such that the
Segal maps are isomorphisms we obtain a category $C$ by taking  $$
Ob (C) := A_0
$$
and
$$
Hom _C(x,y):= A_1(x,y)
$$
with the latter defined as the inverse image of $(x,y)$ under the map (given
by the pair of face maps) $A_1\rightarrow A_0 \times A_0$.  The condition on
the Segal maps implies that (with a similar notation)
$$
A_2(x,y,z)\stackrel{\cong}{\rightarrow} A_1(x,y)\times A_1(y,z)
$$
and the third face map $A_2(x,y,z)\rightarrow A_1(x,z)$ thus provides the
composition of morphisms for $C$. By looking at $A_3(x,y,z,w)$ one sees that the
composition is associative and the degeneracy maps in the simplicial set
provide the identity elements.
\end{parag}

\begin{parag}
\label{ncatsdef0}
The notion of weak $n$-category of Tamsamani \cite{Tamsamani}
is a generalization of the above point of view on categories. We present
the definition in a highly recursive way, using the notion of $n-1$-category in
the definition of $n$-category. See \cite{Tamsamani} for a more direct approach.
This definition is based on Segal's delooping machine \cite{Segal} \cite{Adams}.
\end{parag}

\begin{parag}
\label{notstrict}
Note that Tamsamani uses the terminology {\em $n$-nerve} for what we will call
``$n$-category'' since he needed to distinguish this from the notion of
strict $n$-category.  In the present paper we will never speak of strict
$n$-categories and our terminology ``$n$-category'' means weak $n$-category or
$n$-nerve in the sense of \cite{Tamsamani}.
\end{parag}

\begin{parag}
\label{ncatsdef1}
 An {\em
$n$-category} \cite{Tamsamani} is a functor $A$ from $\Delta ^o$ to the category
of $n-1$-categories denoted
$$
p\mapsto A_{p/}
$$
such that $0$ is mapped to a set
\footnote{
Recursively an $n$-category which is a set is  a constant functor
where the $A_{p/}$ are all the same set---considered as $n-1$-categories.}
$A_0$
and such that the {\em Segal maps}
$$
A_{p/} \rightarrow A_{1/} \times _{A_0} \ldots \times _{A_0}A_{1/}
$$
are equivalences of $n-1$-categories (cf \ref{defequiv1} below).
\end{parag}

\begin{parag}
\label{multisimplicial}
The {\em category of $n$-categories} denoted $n-Cat$ is just the category whose
objects are as above and whose morphisms are the morphisms
strictly preserving the structure.
It is a subcategory of $Hom (\Delta ^o, (n-1)-Cat)$. Working this out
inductively we find in the end that $n-Cat$ is a subcategory of
$Hom ((\Delta ^n)^o, Sets)$, in other words an $n$-category is a
certain kind of multisimplicial set.
The multisimplicial set is denoted
$$
(p_1,\ldots , p_n )\mapsto A_{p_1,\ldots
p_n}
$$
and the $(n-1)$-category $A_{p/}$  itself considered as a multisimplicial
set has the expression
$$
A_{p/} = \left( (q_1,\ldots , q_{n-1})\mapsto A_{p,q_1,\ldots ,
q_{n-1}}\right) .
$$
\end{parag}

\begin{parag}
\label{theta}
The condition that $A_0$ be a set yields by induction the condition that
if $p_i=0$ then the functor $A_{p_1,\ldots , p_n}$ is independent of the
$p_{i+1}, \ldots , p_n$. We call this the {\em constancy condition}. In
\cite{nCAT} we introduce the category $\Theta ^n$ which is the quotient of
$\Delta ^n$ having the property that functors $(\Theta ^n)^o\rightarrow Sets$
correspond to functors on $\Delta ^n$ having the above constancy property.
Now $n-Cat$ is a subcategory of the category of presheaves of
sets on $\Theta ^n$.
\end{parag}

\begin{parag}
Before discussing the notion of equivalence which enters into the above
definition we take note of the relationship with \ref{catnerve}.  If $A$ is an
$n$-category then its {\em set of objects} is the set $A_0$.  The face maps
give a morphism from $n-1$-categories to sets
$$
A_{p/}\rightarrow A_0 \times \ldots \times A_0
$$
and  we denote by $A_{p/}(x_0,\ldots , x_p)$ the $n-1$-category
inverse image of $(x_0,\ldots , x_p)$ under this map. For two objects $x,y\in
A_0$ the $n-1$-category $A_{1/}(x,y)$ is the {\em $n-1$-category of morphisms
from $x$ to $y$}.  This is the essential part of the structure which
corresponds, in the case of categories, to the $Hom$ sets.  One could adopt
the notation
$$
Hom _A(x,y):= A_{1/}(x,y).
$$
The condition that the Segal maps are equivalences of
$n-1$-categories says that the  $A_{p/}(x_0,\ldots , x_p)$ are determined up to
equivalence by the $A_{1/}(x,y)$. The role of the higher $A_{p/}(x_0,\ldots ,
x_p)$ is to  provide the composition (in the case $p=2$) and to keep track of
the higher homotopies of associativity $(p\geq 3$). Contrary to the case of
$1$-categories, here we need to go beyond
\footnote{One might conjecture
that it suffices to stop at $p=n+2$.}
$p=3$.
\end{parag}

\begin{parag}
\label{defequiv1}
In order for the recursive definition of $n$-category given in \ref{ncatsdef1}
to make sense, we need to know what an {\em equivalence} of $n$-categories is.
For this we generalize the usual notion for categories: an equivalence of
categories is a morphism which is (1) fully faithful and (2) essentially
surjective. We would like to define what it means for a functor between
$n$-categories $f:A\rightarrow B$ to be an equivalence. The generalization
of the
fully faithful condition is immediate: we require that for any objects
$x,y\in A_0$ the morphism
$$
f: A_{1/} (x,y) \rightarrow B_{1/}(f(x), f(y))
$$
be an equivalence of $n-1$-categories (and we are supposed to know what that
means by recurrence).
\end{parag}

\begin{parag}
\label{essentialsurjectivity}
The remaining question is how to define the notion
of essential surjectivity.  Tamsamani does this by defining a truncation
operation $T$ from $n$-categories to $n-1$-categories (a generalization of the
truncation of topological spaces used in the Postnikov tower).  Applying this
$n$ times to an $n$-category $A$ we obtain a set $T^nA$ which can also be
denoted $\tau _{\leq 0}A$.  This set is the set of ``objects of $A$ up to
equivalence'' where equivalence of objects is thought of in the
$n$-categorical sense.  We say that $f:A\rightarrow B$ is {\em essentially
surjective} if the induced map
$$
\tau _{\leq 0} (f) : \tau _{\leq 0} A \rightarrow \tau _{\leq 0} B
$$
is a surjection of sets. One has in fact that if $f$ is an equivalence according
to the above definition then $\tau _{\leq 0} f$ is an isomorphism.
\end{parag}

\begin{parag}
\label{anotherapproach}
Another way to approach the definition of
$\tau _{\leq 0}A$ is by induction in the following way. Suppose we know what
$\tau _{\leq 0}$ means for $n-1$-categories. Then for an $n$-category $A$ the
simplicial set $p\mapsto \tau _{\leq  0} (A_{p/})$ satisfies the condition that
the Segal maps are isomorphisms, so it is the nerve of a $1$-category. This
category may be denoted $\tau _{\leq 1} A$.  We then define $\tau _{\leq 0}A$
to be the set of isomorphism  classes of objects in the $1$-category
$\tau _{\leq 1} A$.
\end{parag}

The above definition is highly recursive. One must check that everything is
well defined and available when it is needed. This is done in \cite{Tamsamani}
although the approach there avoids some of the inductive definitions above.

\bigskip

\subnumero{The closed model structure}

An $n$-category is a presheaf of sets on $\Theta ^n$ (\ref{theta}) satisfying
certain conditions as described above. Unfortunately $n-Cat$ considered as a
subcategory of the category of presheaves, is not closed under pushout or fiber
product. This remark is the starting point for \cite{nCAT}.  There, one
considers the full category of presheaves of sets on $\Theta$ (these presheaves
are called {\em $n$-precats}) and \cite{nCAT} provides a closed model structure
(cf \cite{Quillen} \cite{QuillenAnnals} \cite{Jardine})
on the category $nPC$ of $n$-precats, corresponding to
the homotopy theory of $n$-categories. In this section we briefly recall how
this works.

\begin{parag}
\label{theta2}
It is more convenient for the purposes of the closed model structure
to work with presheaves over the category $\Theta ^n$ (cf \ref{theta} above),
defined
be the quotient of the cartesian product $\Delta ^n$ obtained by identifying
all of the objects $(M, 0, M')$ for fixed $M = (m_1,\ldots , m_k)$ and variable
$M'= (m'_1, \ldots , m'_{n-k-1})$.  The object of $\Theta ^n$ corresponding
to the class of $(M,0,M')$ with all $m_i >0$ will be denoted $M$.
Two morphisms from $M$ to $M'$ in $\Delta ^n$ are identified if they both factor
through something of the form $(u_1,\ldots , u_i, 0, u_{i+2}, \ldots , u_n)$
and if their first $i$ components are the same.
\end{parag}

\begin{parag}
\label{precat}
An {\em $n$-precat} is defined to be a presheaf on the category $\Theta ^n$.
This corresponds to an $n$-simplicial set $(\Delta ^n)^o\rightarrow Sets$
which satisfies the constancy condition (cf \ref{theta}). The category $nPC$
of $n$-precats (with morphisms being the morphisms of presheaves) is to be given
a closed model structure.
\end{parag}

\begin{parag}
\label{prelims}
Note for a start that $nPC$ is closed under arbitrary products and coproducts,
what is more (and eventually important for our purposes) it admits an internal
$\underline{Hom}(A,B)$. These statements  come simply from the fact that $nPC$
is a category of presheaves over something.

We denote the coproduct or pushout of $A\rightarrow B$ and $A\rightarrow C$
by $B\cup ^AC$. We denote fiber products by the usual notation.
\end{parag}

\begin{parag}
\label{cofibs}
{\em Cofibrations:}
A morphism $A\rightarrow B$ of $n$-precats is a {\em cofibration} if the
morphisms $A_M \rightarrow B_M$ are injective whenever $M\in \Theta ^n$ is
an object of non-maximal length, i.e. $M= (m_1,\ldots , m_k, 0,\ldots , 0)$
for $k< n$.  The case of sets ($n=0$) shows that we can't require
injectivity at the top level $n$, nor do we need to.

We often use the notation $A\hookrightarrow B$ for a cofibration, not meaning
to imply injectivity at the top level.
\end{parag}

\begin{parag}
\label{we}
{\em Weak equivalences:}
In order to say when a morphism $A\rightarrow B$ of $n$-precats is a ``weak
equivalence'' we have to do some work. In \cite{Tamsamani} was defined the
notion of equivalence between $n$-categories (cf \ref{defequiv1} above), but an
$n$-precat is not yet an $n$-category. We need an operation which specifies the
intended relationship between our $n$-precats and $n$-categories. This is the
operation $A\mapsto Cat(A)$ which to any $n$-precat associates an $n$-category
together with  morphism of precats $A\rightarrow Cat(A)$, basically by throwing
onto $A$ in a minimal way all of the elements which are needed in order to
satisfy the  definition of being an $n$-category. See \cite{nCAT} \S 2 for the
details of this. Now we say that a morphism
$$
A\rightarrow B
$$
of $n$-precats is  a {\em weak equivalence} if the induced morphism of
$n$-categories
$$
Cat(A)\rightarrow Cat(B)
$$
is an equivalence as defined in \cite{Tamsamani}---described in \ref{defequiv1}
and \ref{essentialsurjectivity} above.
\end{parag}

\begin{parag}
\label{trivcofibs}
{\em Trivial cofibrations:}
A morphism $A\rightarrow B$ is said to be a {\em trivial cofibration}
if it is a cofibration and a weak equivalence.
\end{parag}

\begin{parag}
\label{fibs}
{\em Fibrations:}
A morphism $A\rightarrow B$ of $n$-precats is said to be a {\em fibration}
if it satisfies the following lifting property: for every trivial cofibration
$E'\hookrightarrow E$ and every morphism $E\rightarrow B$ provided with a
lifting over $E'$ to a morphism $E'\rightarrow A$, there exists an extension of
this to a lifting $E\rightarrow A$.

An $n$-precat $A$ is said to be {\em fibrant} if the canonical (unique)
morphism $A\rightarrow \ast$ to the constant presheaf with values one point,
is a fibration.

A fibrant $n$-precat is, in particular, an $n$-category. This is because the
elements which need to exist to give an $n$-category may be obtained as
liftings of certain standard trivial cofibrations (those denoted $\Sigma
\rightarrow h$ in \cite{nCAT}).
\end{parag}

\begin{theorem}
\label{cmc}
{\rm (\cite{nCAT} Theorem 3.1)}
The category $nPC$ of $n$-precats with the above classes of cofibrations, weak
equivalences and fibrations, is a closed model category.
\end{theorem}

The basic ``yoga'' of the situation is that when we want to look at coproducts,
one of the morphisms should be a cofibration; when we want to look at fiber
products, one of the morphisms should be a fibration; and when we want look at
the space of morphisms from $A$ to $B$, the first object $A$ should be cofibrant
(in our case all objects are cofibrant) and the second object $B$ should be
fibrant.

\begin{parag}
\label{explaincmc}
We explain more precisely what information is contained in the above theorem,
by explaining the axioms for a closed model category structure
(CM1--CM5 of \cite{QuillenAnnals}).  These are proved as such in \cite{nCAT}.

\noindent
{\em CM1}---This says that $nPC$ is closed under finite (and in our
case, arbitrary) direct and inverse limits (\ref{prelims}).

\noindent
{\em CM2}---Given composable morphisms
$$
X\stackrel{f}{\rightarrow}Y\stackrel{g}{\rightarrow}Z,
$$
if any two of $f$ or $g$ or $g\circ f$
are weak equivalences  then the  third is
also a weak equivalence.

\noindent
{\em CM3}---The classes of cofibrations, fibrations and
weak equivalences are closed under retracts. We don't explicitly use this
condition (however it is the basis for the property \ref{pushoutA} below).

\noindent
{\em CM4}---This says that a pair of a cofibration
$E'\rightarrow E$ and a fibration $A\rightarrow B$ have the lifting property (as
in the definition of fibration \ref{fibs}) if either one of the morphisms
is a weak equivalence. Note that the lifting property when $E'\rightarrow E$ is
a weak equivalence (i.e. trivial cofibration) is just the definition that
$A\rightarrow B$ be fibrant \ref{fibs}. The other half, the lifting property
for an arbitrary cofibration when $A\rightarrow B$ is a weak equivalence,
comes from what Jardine calls ``Joyal's trick'' \cite{Jardine}.

\noindent
{\em CM5}---This says that any morphism $f$ may be factored as a
composition $f= p\circ i$ of a cofibration followed by fibration, and either
one of $p$ or $i$ may be assumed to be a weak equivalence.
\end{parag}

\begin{parag}
\label{pushoutA}
Another axiom in Quillen's original point of view (Axiom M3 on
page 1.1 of \cite{Quillen}) is that if $A\rightarrow
B$ is a trivial cofibration and $A\rightarrow C$ is any morphism then
$C\rightarrow B\cup ^AC$ is again a trivial cofibration; and similarly the dual
condition for fibrant weak equivalences and fiber products.
In the closed model category setting this becomes a consequence of the axioms
CM1--CM5, see \cite{QuillenAnnals}.

In the proof of \cite{nCAT}
(modelled on that of \cite{Jardine}) the main step which is done first
(\cite{nCAT} Lemma 3.2) is to
prove this property of preservation of trivial cofibrations by coproducts.
. (On the other hand, note that with our definition
\ref{fibs} of fibrations, the preservation by fiber products is obvious).
\end{parag}

\begin{parag}
\label{perspective1}
We now try to put these properties in perspective in view of how we will use
them. If $A$ is any $n$-precat then applying {\em CM5} to the morphism
$A\rightarrow \ast$ we obtain a factorization
$$
A\rightarrow A' \rightarrow \ast
$$
with the first morphism a trivial cofibration, and the second morphism a
fibration. Thus $A'$ is a fibrant object. We call such a trivial cofibration to
a fibrant object $A\hookrightarrow A'$ a {\em fibrant replacement for $A$}.

In the constructions of \cite{Jardine}, \cite{nCAT} one
obtains the fibrant replacement by adding onto $A$ the pushouts by ``all
possible'' trivial cofibrations, making use of \ref{pushoutA}. The notion
of ``all possible'' has to be refined in order to avoid set-theoretical
problems:
actually one looks at $\omega$-bounded cofibrations. The number of them is
bounded by the maximum of $2^{\omega}$ or the cardinality of $A$.

When looking at morphisms into an $n$-category $C$ it is important that $C$ be
fibrant, for then  we obtain extension properties along trivial
cofibrations.  In particular, we will only define what it means for limits to
exist in $n$-categories $C$ which are fibrant.

When we finally get to our definition of  the $n+1$-category $nCAT$ below, it
will not be fibrant. Thus one of the main steps is to choose a fibrant
replacement $nCAT \hookrightarrow nCAT'$.
\end{parag}

\begin{parag}
\label{interval}
There is  a nice ``interval'' in our closed model category (in contrast with the
general situation envisioned by Quillen in \cite{Quillen}).
Let $\overline{I}$ denote the $1$-category with two objects $0,1$ and with
unique
morphisms going in either  direction between them, whose compositions are the
identity. Without changing notation, we can consider
$\overline{I}$ as an $n$-category (pull back by the obvious morphism $\Theta
^n\rightarrow \Theta ^1 = \Delta $).

\noindent
{\em Claim:}
Suppose $C$ is a fibrant $n$-category. Then two objects $x,y\in C_0$ are
equivalent (i.e. project to the same thing in $\tau _{\leq 0} C$
cf \ref{anotherapproach}) if and only if there exists a morphism
$\overline{I}\rightarrow C$ sending $0$ to $x$ and $1$ to $y$.

To prove this note that one direction is obvious: if there exists such a
morphism then by functoriality of $\tau _{\leq 0}$ $x$ and $y$ are equivalent
(because $\tau _{\leq 0}(\overline{I})=\ast$). For the other direction, suppose
$x$ and $y$ are equivalent. Use Proposition 6.5 of \cite{nCAT} which says that
there is an $n$-category $K$ with objects $0,1$ such that $K \rightarrow \ast$
is a weak equivalence, and there is a morphism $K\rightarrow C$ sending $0$ to
$x$ and $1$ to $y$.  Applying the factorization statement {\em CM5}
to the morphism
$$
K\cup ^{\{ 0 , 1\} } \overline{I} \rightarrow \ast
$$
we obtain a cofibration
$$
K\cup ^{\{ 0 , 1\} } \overline{I}\hookrightarrow A
$$
such that $A\rightarrow \ast$ is a weak equivalence. It follows from {\em CM2}
applied to $K\rightarrow A\rightarrow \ast$
that $K\rightarrow A$ is a weak equivalence, thus it is a trivial cofibration.
Now the fibrant property of $C$ implies that our morphism $K\rightarrow C$
extends to a morphism $A\rightarrow C$. This morphism restricted to
$\overline{I}\hookrightarrow A$ provides a morphism $\overline{I}\rightarrow C$
sending $0$ to $x$ and $1$ to $y$. This proves the other direction of the
claimed statement.
\end{parag}

\begin{parag}
\label{surje}
As a corollary of the above construction, suppose $f:A\rightarrow B$ is a
fibrant
morphism of fibrant  $n$-categories. Suppose that $a\in A_0$ and $b\in
B_0$ are objects such that $f(a)$ is equivalent to $b$ (i.e.
$f(a)$ is equal to $b$ in $\tau _{\leq 0} B$).  Then there is a different object
$a'\in A_0$ equivalent to $a$ such that $f(a')=b$ in $B_0$. To prove this, note
that the equivalence between $f(a)$ and $b$ corresponds by \ref{interval}
to a morphism $\overline{I}\rightarrow B$ sending $0$ to $f(a)$ and $1$ to
$b$. We have a lifting $a$ over $\{ 0\}$. The inclusion
$$
\{ 0\}  \subset \overline{I}
$$
is a trivial cofibration, so the fibrant property of $f$ means that there is a
lifting to a morphism $\overline{I} \rightarrow A$. The image of $1$ by
this map is an object $a'$ equivalent to $a$ and projecting to $b$.

A variant says that if $f:A\rightarrow B$ is a fibrant morphism between
fibrant $n$-categories and if $f$ is an equivalence then $f$ is surjective on
objects.  To obtain this note that
essential surjectivity of $f$ means that  every object $b$
is equivalent to some $f(a)$, then apply the previous statement.
\end{parag}

\begin{parag}
\label{internal1}
One of the main advantages to using a category of presheaves $nPC$ as
underlying category is that we obtain an internal $\underline{Hom}(A,B)$
between two $n$-precats. This represents a functor: a map
$$
E\rightarrow \underline{Hom}(A,B)
$$
is the same thing as a morphism $A\times E \rightarrow B$.

Of course for arbitrary $A$ and $B$, the internal $\underline{Hom}(A,B)$ will
not have any reasonable properties, for example it will not transform
equivalences of the $A$ or $B$ into equivalences. This situation is
rectified by imposing the hypothesis that $B$ should be fibrant.
\end{parag}

\begin{parag}
\label{internal2}
We describe some of the results saying that the
internal $\underline{Hom}(A,B)$ works nicely when $B$ is fibrant. The following
paragraphs are Theorem 7.1 and Lemma 7.2 of \cite{nCAT}.

Suppose $A$ is an $n$-precat  and $B$ is a fibrant $n$-precat. Then the
internal  $
\underline{Hom}(A,B)$ of
presheaves over $\Theta ^n$ is a fibrant $n$-category. Furthermore if
$B'\rightarrow B$ is a fibrant morphism then
$\underline{Hom} (A, B')\rightarrow \underline{Hom} (A, B)$ is fibrant.
Similarly if $A\hookrightarrow A'$ is a cofibration and if
$B$ is fibrant then $\underline{Hom}(A', B)\rightarrow \underline{Hom}(A,B)$ is
fibrant.

Suppose $A\rightarrow A'$ is a weak equivalence,  and  $B$
fibrant. Then
$$
\underline{Hom} (A', B)\rightarrow \underline{Hom} (A, B)
$$
is an
equivalence of $n$-categories.

If $B\rightarrow B'$ is an equivalence
of fibrant $n$-precats then $\underline{Hom}(A,B)\rightarrow
\underline{Hom}(A,B')$ is an equivalence.

Suppose $A\rightarrow B$ and $A\rightarrow C$ are cofibrations. Then
$$
\underline{Hom} (B\cup ^AC, D) = \underline{Hom} (B, D)
\times _{\underline{Hom}(A,D)}\underline{Hom}(C,D).
$$
\end{parag}

\begin{parag}
\label{homotopic1}
We can relate several different versions of the notion of two morphisms being
homotopic. Suppose $A$ and $B$ are $n$-precats with $B$ fibrant.
According to Quillen's definition \cite{Quillen}, two maps
$f: A\rightarrow B$ are {\em homotopic} if there is a diagram
$$
A\tworightarrows A' \rightarrow A
$$
such that all morphisms are weak equivalences, the first two morphisms are
cofibrations, and such that the compositions are the
identity of $A$, plus a morphism $A'\rightarrow B$ inducing $f$ and $g$ on
the two copies of $A$.

In our situation, if $B$ is fibrant then $\underline{Hom}(A,B)$ is a fibrant
$n$-category  whose objects are the morphisms $A\rightarrow B$. Two morphisms
are equivalent objects
in this $n$-category (cf \ref{essentialsurjectivity} above) if and only if they
are homotopic in Quillen's sense (this is \cite{nCAT} Lemma 7.3).
\end{parag}

\begin{parag}
\label{homotopic2}
In the above situation apply the claim of \ref{interval}. Two objects $f,g$ of
the fibrant  $n$-category  $\underline{Hom}(A,B)$ are equivalent if and only if
there is a morphism
$$
\overline{I}\rightarrow  \underline{Hom}(A,B)
$$
sending $0$ to $f$ and $1$ to $g$.  Such a morphism corresponds to a
map
$$
A\times \overline{I} \rightarrow B;
$$
so we can finish up by saying that two morphisms $f,g:A\rightarrow B$
are homotopic if and only if there exists a map
$$
A\times \overline{I}\rightarrow B
$$
restricting to $f$ on $A\times \{ 0\} $ and to $g$ on $A\times \{ 1\}$.
\end{parag}

\begin{parag}
\label{charequiv}
We obtain from CM4 the following characterization of fibrant weak equivalences.
A morphism $f: A\rightarrow B$ is a fibrant weak equivalence if and only if
it satisfies the lifting property for any cofibration $E'\hookrightarrow E$.
To prove this, note that CM4 shows that a fibrant weak equivalence has this
property. If $f$ has this property then it is fibrant (the case of
$E'\hookrightarrow E$ a trivial cofibration).  The morphisms
of $n-1$-categories
$$
A_{1/} (x,y)\rightarrow B_{1/}(f(x) , f(y))
$$
also have the same property (one can see this using the construction $\Upsilon
(E)$ below) and $f$ is surjective on objects (by the case $\emptyset
\hookrightarrow \ast$). Therefore $f$ is an equivalence.

We can give the following variant characterizing when a morphism is an
equivalence (not necessarily fibrant).  We say that a morphism $f: A\rightarrow
B$ between fibrant $n$-categories has the {\em homotopical lifting property for
$E'\hookrightarrow E$} if, given a morphism $v:E\rightarrow B$ and a lifting
$u':E'\rightarrow A$, there is a homotopy from $v$ to
a new morphism $v_1$, a lifting $u_1$ of $v_1$, and a homotopy
from $u'$ to $u'_1$ (the restriction of $u_1$ to $E'$) lifting the
homotopy from $v'$ to $v'_1$ (restriction of our first homotopy to $E'$).
In this definition
we can use any of the  equivalent notions of homotopy
\ref{homotopic1}, \ref{homotopic2} above.

{\em Claim:} A morphism $f: A\rightarrow B$ between two fibrant
$n$-categories is an equivalence if and only if it satisfies the
homotopical lifting property for all $E'\hookrightarrow E$.

To prove this, use CM5 to factor $A\rightarrow A'\rightarrow B$
with the first morphism a trivial cofibration and the second morphism fibrant.
Note that $A'$ is again fibrant. The statement being a homotopical one,
the same hypothesis holds for $A'\rightarrow B$. If we can prove that
$A'\rightarrow B$ is an equivalence then the composition with the trivial
cofibration $A\rightarrow A'$ will be a weak equivalence.
Thus we
may reduce to the case where $A\rightarrow B$ is a fibrant morphism. Now given
$E'\hookrightarrow E$ with $E\rightarrow B$ lifting to $E'\rightarrow A$,
choose homotopies
$$
E\times \overline{I}\rightarrow B
$$
and lifting
$$
E'\times \overline{I}\rightarrow A
$$
compatible with a lifting $E\times \{ 1\} \rightarrow A$
as in the definition of the homotopical lifting property. These give a lifting
$$
E'\times \overline{I}\cup ^{E'\times \{ 1\} }
E \times \{ 1\} \rightarrow A,
$$
and the morphism
$$
E'\times \overline{I}\cup ^{E'\times \{ 1\} }
E \times \{ 1\} \rightarrow E\times \overline{I}
$$
is a trivial cofibration, so by the fibrant property of $A\rightarrow B$
(which we are now assuming) there is a lifting
$$
E\times \overline{I} \rightarrow A.
$$
The restriction to $E\times \{ 0\}$ gives the desired lifting of the original
morphism $E\rightarrow B$, coinciding with the given lifting on $E'$.
This proves that $A\rightarrow B$ satisfies the lifting criterion given above
so it is a fibrant weak equivalence. This completes the proof of one direction
of the claim.
A similar argument (using CM2) gives the other direction.
\end{parag}

\bigskip

\subnumero{Families of $n$-categories}

\begin{parag}
\label{ncat1}
Using the internal $\underline{Hom}(A,B)$ of \ref{internal} between fibrant
$n$-categories, we define the $n+1$-category $nCAT$ of all fibrant
$n$-categories (cf \cite{nCAT} \S 7).
This is the ``right'' category of $n$-categories, and is not to be
confused with the first approximation $n-Cat$ as defined in
\ref{multisimplicial}
above.

The objects of $nCAT$ are the fibrant
$n$-categories.  Between any two objects we have an $n$-category of morphisms
$\underline{Hom}(A,B)$. Composition of morphisms gives a morphism
of $n$-categories
$$
\underline{Hom}(A,B)\times \underline{Hom}(B,C)\rightarrow \underline{Hom}(A,C),
$$
which is strictly associative and has a unit element, the identity morphism.
Using this we obtain an $n+1$-category $nCAT$: to be precise, if
$A_0,\ldots , A_p$ are objects then
$$
nCAT_{p/}(A_0,\ldots , A_p):= \underline{Hom}(A_0,A_1)\times
\ldots \times \underline{Hom}(A_{p-1},A_p)
$$
which organizes into a simplicial collection using the projections or, where
necessary, the composition morphisms.  The Segal maps are actually isomorphisms
here so this is an $n+1$-category.
\end{parag}

\begin{parag}
\label{ncat2}
Unfortunately, $nCAT$ is not a fibrant $n+1$-category, although  it does
have the
property that the $nCAT_{p/}$ are fibrant. Because of this, we must choose a
fibrant replacement
$$
nCAT \hookrightarrow nCAT'.
$$
\end{parag}

\begin{parag}
\label{families1}
Basic to the present paper is the notion of {\em family of $n$-categories
indexed by an $n+1$-category $A$}, which is defined using our fibrant
replacement (\ref{ncat2}) to be a morphism $A\rightarrow nCAT'$.

The {\em $n+1$-category of all families indexed by $A$} is the
$n+1$-category $\underline{Hom}(A, nCAT')$.
\end{parag}

\begin{parag}
\label{families2}
Suppose $\psi, \psi ' : A\rightarrow nCAT'$ are families. A {\em morphism}
from $\psi$ to $\psi '$ is an object of the $n$-category
$\underline{Hom}(A, nCAT')_{1/}(\psi , \psi ')$.  Let $I$ be the category
with two objects $0,1$ and a morphism from $0$ to $1$ (in our notations below
this will also be the same as what we will call $\Upsilon (\ast )$).
The set of
objects of $\underline{Hom}(A, nCAT')_{1/}(\psi , \psi ')$ is equal
to the set of
morphisms  $$
I\rightarrow  \underline{Hom}(A, nCAT')
$$
sending $0$ to $\psi$ and $1$ to $\psi '$. In view of the definition
of internal $\underline{Hom}$ this is the same thing as a morphism
$$
 A \times I \rightarrow nCAT'
$$
restricting on $A\times \{ 0\}$ to $\psi $ and on $A\times \{ 1\}$ to $\psi '$.
\end{parag}

\bigskip

\subnumero{The construction $\Upsilon$}

We will now introduce some of our main tools for the present paper. The basic
idea is that we often would like to talk about the basic $n$-category with two
objects (denoted $0$ and $1$) and with a given $n-1$-category $E$ of morphisms
from $0$ to $1$ (but no morphisms in the other direction and only identity
endomorphisms of $0$ and $1$).  We call this $\Upsilon (E)$.  To be more
precise we do this on the level of precats: if $E$ is an $n-1$-precat then we
obtain an $n$-precat $\Upsilon (E)$.  The main property of this construction is
that if $A$ is any $n$-category then a morphism of $n$-precats
$$
f:\Upsilon (E) \rightarrow A
$$
corresponds exactly to a choice of two objects $x=f(0)$ and $y=f(1)$ together
with a morphism of $n-1$-precats $E\rightarrow A_{1/}(x,y)$.

One can see $\Upsilon (E)$ as the universal $n$-precat $A$ with two objects
$x,y$ and a map $E\rightarrow A_{1/}(x,y)$.

\begin{parag}
\label{upsilon}
We also need more general things of the form $\Upsilon ^2(E,F)$ having objects
$0,1,2$ and similarly a $\Upsilon ^3$.
(These will not have quite so simple an interpretation as
universal objects.) Thus we present the definition in a
general way.

Suppose $E_1,\ldots , E_k$ are $n-1$-precats. Then we define the $n$-precat
$$
\Upsilon ^k(E_1,\ldots , E_k)
$$
in the following way. Its object set is the set with $k+1$ elements denoted
$$
\Upsilon ^k(E_1,\ldots , E_k)_0 = \{ 0,\ldots , k\} .
$$
Then
$$
\Upsilon ^k(E_1,\ldots , E_k)_{p/}(y_0, \ldots , y_p)
$$
is defined to be empty if any $y_i > y_j$ for $i<j$, equal to
$\ast$ if $y_0=\ldots = y_p$, and otherwise
$$
\Upsilon ^k(E_1,\ldots , E_k)_{p/}(y_0, \ldots , y_p):= E_{y_0} \times \ldots
\times E_{y_p}.
$$
\end{parag}

\begin{parag}
For example when $k=1$ (and we drop the superscript $k$ in this case)
$\Upsilon E$ is the $n$-precat with two objects $0,1$ and with $n-1$-precat of
morphisms from $0$ to $1$ equal to $E$. Similarly $\Upsilon ^2(E,F)$ has
objects $0,1,2$ and morphisms $E$ from $0$ to $1$, $F$ from $1$ to $2$ and
$E\times F$ from $0$ to $2$.  We picture $\Upsilon ^k(E_1,\ldots , E_k)$
as a $k$-gon (an edge for $k=1$, a triangle for $k=2$, a tetrahedron for $k=3$).
The edges are labeled with single $E_i$, or products $E_i \times \ldots , E_j$.
\end{parag}

\begin{parag}
There are inclusions of these $\Upsilon^k$ according to the faces of the
$k$-gon. The principal faces give inclusions
$$
\Upsilon ^{k-1}(E_1,\ldots , E_{k-1})\hookrightarrow \Upsilon ^k(E_1, \ldots ,
E_k),
$$
$$
\Upsilon ^{k-1}(E_2,\ldots , E_{k})\hookrightarrow \Upsilon ^k(E_1, \ldots ,
E_k),
$$
and
$$
\Upsilon ^{k-1}(E_1,\ldots , E_i\times E_{i+1}, \ldots  , E_{k})\hookrightarrow
\Upsilon ^k(E_1, \ldots , E_k).
$$
The inclusions of lower levels are deduced from these by induction. Note that
these faces $\Upsilon ^{k-1}$ intersect along appropriate $\Upsilon ^{k-2}$.
\end{parag}

\begin{remark}
\label{upsistar}
$\Upsilon (\ast )= I$ is the category with objects $0,1$ and with a unique
morphism from $0$ to $1$.  A map $\Upsilon (\ast )\rightarrow A$ is the same
thing as a pair of objects $x,y$ and a $1$-morphism from $x$ to $y$, i.e. an
object of $A_{1/}(x,y)$.
\end{remark}

\bigskip

Another way of constructing the $\Upsilon ^k$ is given in the following remarks
\ref{interpupsilon1}--\ref{interpupsilon3}.

\begin{parag}
\label{interpupsilon1}
For an $n-1$-precat $E$, denote by
$[p](E)$ the universal $n$-precat  $A$ with objects $x_0,\ldots ,x_p$ and with a
morphism $E\rightarrow A_{p/}(x_0,\ldots , x_p)$. This can be described
explicitly by saying that $[p](E)$ has objects $0,1,\ldots , p$, and for a
sequence of objects $i_1,\ldots , i_k$ the $n-1$-precat $[p](E)_{k/}(i_1,\ldots
, i_k)$ is empty if some $i_j < i_{j+1}$, is equal to $\ast$ if all $i_j$ are
equal, and is equal to $E$ if some $i_j < i_{j+1}$.
\end{parag}

\begin{parag}
\label{interpupsilon2}
One has $[1]E= \Upsilon (E)$.  The construction of the higher $\Upsilon ^k$ may
be described inductively as follows: we will construct $\Upsilon ^k(E_1,\ldots
, E_k)$ together with a morphism
$$
in:[k](E_1\times \ldots \times E_k)\rightarrow \Upsilon ^k(E_1,\ldots
, E_k).
$$
Suppose we have constructed these maps up to $k-1$.
Note that the first and last face morphisms coupled with the projections onto
the first and last $k-1$ factors give a map
$$
[k-1](E_1\times \ldots \times
E_{k}) \cup ^{[k-2](E_1\times \ldots \times
E_{k})}
[k-1](E_1\times \ldots \times
E_{k})
$$
$$
\stackrel{\alpha}{\rightarrow}
[k](E_1,\ldots , E_k),
$$
but on the other hand the projections onto subsets of factors of the product
$E_1\times \ldots \times E_k$ together with the maps $in$ in our inductive
construction for $k-1$ and $k-2$ give a map
$$
[k-1](E_1\times \ldots \times
E_{k}) \cup ^{[k-2](E_1\times \ldots \times
E_{k})}
[k-1](E_1\times \ldots \times
E_{k})
$$
$$
\stackrel{\beta}{\rightarrow}
\Upsilon ^{k-1}(E_1,\ldots , E_{k-1})\cup ^{
\Upsilon ^{k-2}(E_2,\ldots , E_{k-1})}
\Upsilon ^{k-1}(E_2,\ldots , E_{k}).
$$
Finally, $\Upsilon ^{k}(E_1,\ldots , E_k)$ is the coproduct of the maps
$\alpha$ and $\beta$.

We can think of this as saying that  $\Upsilon ^{k}(E_1,\ldots , E_k)$
is obtained by adding on the cell
$[k](E_1\times \ldots \times E_k)$ to the coproduct
$$
\Upsilon ^{k-1}(E_1,\ldots , E_{k-1})\cup ^{
\Upsilon ^{k-2}(E_2,\ldots , E_{k-1})}
\Upsilon ^{k-1}(E_2,\ldots , E_{k})
$$
of the earlier things we
have inductively constructed.
\end{parag}

\begin{parag}
\label{interpupsilon3}
The case $k=2$ is simpler to write down and is worth mentioning separately.
Recall that for $k=1$ we just had $\Upsilon (E)= [1](E)$. The next step is
$$
\Upsilon ^2(E,F) = [2](E\times F) \cup ^{[1](E\times F)\cup ^{\ast }
[1](E\times F)} ([1](E) \cup ^{\ast} [1](F)).
$$
\end{parag}

\begin{parag}
\label{trivinclusions}
One thing which we often will need to know below is when an inclusion from a
union of faces, into the whole $\Upsilon ^k$, is a trivial cofibration.
For $k=2$ the only inclusion which is a trivial cofibration is
$$
\Upsilon (E_1)\cup ^{\{ 1\} } \Upsilon (E_2) \hookrightarrow \Upsilon
^2(E_1,E_2). $$

For $k=3$ we denote our inclusions in shorthand notation where
$0,1,2,3$ refer to the vertices. To fix notations, the above inclusion for
$k=2$ would be noted
$$
(01) + (12) \subset (012).
$$
Now for $k=3$ the inclusions which are trivial cofibrations are:
$$
(01) + (12) + (23) \subset (0123)
$$
(which is the standard one, coming basically from the definition of
$n$-category);
and then some others which we obtain from this standard one by adding in
triangles on the right, keeping equivalence with $(01) + (12) + (23)$
according to the result for $k=2$:
$$
(01) + (123)\subset (0123),
$$
$$
(012) + (23) \subset (0123),
$$
$$
(012) + (123)\subset (0123),
$$
$$
(012) + (023)\subset (0123),
$$
$$
(013) + (123)\subset (0123),
$$
$$
(012) + (013) + (123)\subset (0123),
$$
$$
(012) + (023) + (123)\subset (0123).
$$
Our main examples of inclusions which are {\em not} trivial cofibrations
are when we leave out the first or the last faces:
$$
(012) + (023) + (013)\subset (0123) \;\;\;\; \mbox{not a t.c.};
$$
$$
(013) + (023) + (123)\subset (0123) \;\;\;\; \mbox{not a t.c.}.
$$
We call these the {\em left and right shells}. We shall meet both of them
and denote the left shell as
$$
(012) + (023) + (013) = Shell\Upsilon ^3(E_1,E_2,E_3),
$$
and the right shell as
$$
(013) + (023) + (123) = Shelr\Upsilon ^3(E_1,E_2,E_3).
$$
The main parts of our arguments for limits will consist of saying that
under certain circumstances we have an extension property for morphisms
with respect to these cofibrations which are not trivial.
\end{parag}

\begin{parag}
\label{trivcofibs2}
One of the main technical problems which will be encountered by the reader is
deciding when a morphism between $n$ or $n+1$-precats is a trivial cofibration:
we use this all the time in order to use the fibrant property of
the domain of morphisms we are trying to extend. It is not possible to give all
the details each time that this question occurs, as that would be much too long.
The general principles at work are: to be aware of the examples given in
\ref{trivinclusions}; to use the fact that the coproduct of a trivial
cofibration with something else again yields a trivial cofibration
(\ref{pushoutA}); and to use the fact that if a composable sequence of morphisms
$$
\cdot \stackrel{f}{\rightarrow} \cdot
\stackrel{g}{\rightarrow} \cdot
$$
has composition being a weak equivalence, and one of $f$ or $g$ being
a weak equivalence, then so is the other (\ref{explaincmc} CM2). And of
course to
use any available hypotheses that are in effect saying that certain
morphisms are
trivial cofibrations or equivalences. All of the cases where we need to
know that
something  is a trivial cofibration, can be obtained using these principles.
\end{parag}

\begin{parag}
\label{notationsr}
We will often be considering morphisms of the form
$$
f: \Upsilon ^k(E_1,\ldots , E_k) \rightarrow C.
$$
When we would like to restrict this to a face (or higher order face such
as an edge) then, denoting the face by $i_1,\ldots , i_j$ we denote the
restriction of $f$ to the face by
$$
r_{i_1\ldots i_j}(f).
$$
For example when $k=2$ the restriction of
$$
f: \Upsilon ^2(E,F)\rightarrow C
$$
to the edge $(02)$ (which is a $\Upsilon (E\times F)$) would be denoted
$r_{02}(f)$. The object $f(1)$ could also be denoted $r_1(f)$.

We make the same convention for restricting maps of the form
$$
A\times \Upsilon ^k(E_1,\ldots , E_k) \rightarrow C,
$$
to maps on $A$ times some face of the $\Upsilon ^k(E_1,\ldots , E_k)$.
\end{parag}

\bigskip

\subnumero{Inverting equivalences}

In preparation for \ref{obbyobequiv} we need the following
result. It says that a morphism which is an equivalence has an inverse
which is essentially unique, if the notion of ``inverse'' is defined in the
right way.
It is an $n$-category version of the theorem of \cite{flexible} which gives
a canonical inverse for a homotopy equivalence of spaces.

\begin{theorem}
\label{resttoIff}
For any fibrant $n$-category $C$ the morphism restriction from $\overline{I}$
to $I$:
$$
r:\underline{Hom}(\overline{I}, C)
\rightarrow \underline{Hom}(I,C)
$$
is fully faithful, so $\underline{Hom}(\overline{I}, C)$
is equivalent to the full sub-$n$-category of invertible elements of
$\underline{Hom}(I,C)$.
\end{theorem}
{\em Proof:}
We first construct some trivial cofibrations.

\begin{parag}
\label{ff1}
Recall (\ref{upsistar})
that $I= \Upsilon (\ast )$. The morphism
$$
\Upsilon (E)\cup ^{\{ 1\} } I\rightarrow
\Upsilon ^2(E,\ast )
$$
is a trivial cofibration (\ref{trivinclusions}), so by
\ref{pushoutA} the coproduct with
$$
\Upsilon (E)\cup ^{\{ 1\} } I\rightarrow
\Upsilon (E)\cup ^{\{ 1\} } \overline{I}
$$
gives a trivial cofibration
$$
\Upsilon (E)\cup ^{\{ 1\} } \overline{I}
\rightarrow
\Upsilon ^2(E,\ast )\cup ^{I} \overline{I} .
$$
The morphism
$$
\Upsilon (E) \rightarrow \Upsilon (E)\cup ^{\{ 1\} } \overline{I}
$$
is a weak equivalence (again by \ref{pushoutA} because it is pushout of the
trivial cofibration  $\ast  \rightarrow \overline{I}$). Therefore the composed
morphism
$$
i_{01}:\Upsilon (E)\rightarrow \Upsilon ^2(E,\ast )\cup ^{I} \overline{I}
$$
corresponding to the edge $(01)$
is an equivalence. Thus the projection
$$
\Upsilon ^2(E,\ast )\cup ^{I} \overline{I} \rightarrow \Upsilon (E)
$$
is an equivalence (by \ref{explaincmc} CM2). This in turn implies that the
morphism corresponding to the edge $(02)$
$$
i_{02}:\Upsilon (E)\rightarrow \Upsilon ^2(E,\ast )\cup ^{I} \overline{I}
$$
is a trivial cofibration.
\end{parag}

\begin{parag}
\label{ff2}
A similar argument shows that
$$
i_{02}:\Upsilon (E)\rightarrow \Upsilon ^2(\ast , E)\cup ^{I} \overline{I}
$$
is a trivial cofibration.
\end{parag}

\begin{parag}
\label{ff3}
Next, note that
$$
\Upsilon (E) \times I = \Upsilon ^2(E, \ast ) \cup ^{\Upsilon (E)}
\Upsilon ^2(\ast ,E)
$$
(the square decomposes as a union of two triangles). The morphisms in the
coproduct are both $i_{02}$. Thus if we attach $\overline{I}$ to each
of the intervals $I$ on the two opposite sides of this square,
the
result
$$
\Upsilon (E) \times I\cup ^{\{ 0,1\} \times I}(\{
0,1\} \times \overline{I})
$$
can be seen as a coproduct of the two objects
considered in \ref{ff1} and \ref{ff2} (we don't write this coproduct out).
Combining with the results of those paragraphs, the morphism from the diagonal
$$
\Upsilon (E)\rightarrow \Upsilon (E) \times I\cup ^{\{ 0,1\} \times I}(\{
0,1\} \times \overline{I})
$$
is an equivalence, which in turn implies that the projection
$$
\Upsilon (E) \times I\cup ^{\{ 0,1\} \times I}(\{
0,1\} \times \overline{I}) \rightarrow \Upsilon (E)
$$
is an equivalence  or, equally well, that the inclusion
$$
(\Upsilon (E) \times I)\cup ^{\{ 0,1\} \times I}(\{
0,1\} \times \overline{I}) \hookrightarrow
\Upsilon (E) \times \overline{I}
$$
is a trivial cofibration.
\end{parag}

\begin{parag}
\label{ff4}
Suppose now that $E'\subset E$. Let ${\bf G}$ denote the pushout of
$$
(\Upsilon (E) \times I)\cup ^{\{ 0,1\} \times I}(\{
0,1\} \times \overline{I})
$$
and $\Upsilon (E')\times \overline{I}$ along
$$
(\Upsilon (E') \times I)\cup ^{\{ 0,1\} \times I}(\{
0,1\} \times \overline{I}).
$$
Paragraph \ref{ff3} and the usual (\ref{pushoutA}) and '\ref{explaincmc} CM2)
imply that the
morphism
$$
{\bf G} \hookrightarrow \Upsilon (E)\times \overline{I}
$$
is a trivial cofibration. Note, however, the simpler expression
$$
{\bf G} = (\Upsilon (E) \times I)\cup ^{\Upsilon (E') \times I}
(\Upsilon (E')\times \overline{I}).
$$
\end{parag}

\begin{parag}
\label{ff5}
We are now ready to prove the theorem. Fix $u,v$ objects of
$\underline{Hom}(\overline{I},C)$. Suppose
$E'\hookrightarrow E$ is any cofibration, and suppose given a morphism
$$
E \rightarrow \underline{Hom}(I,C)_{1/}(r(u),r(v))
$$
provided with lifting
$$
E'\rightarrow \underline{Hom}(\overline{I},C)_{1/}(u,v).
$$
These correspond exactly to a morphism
$$
{\bf G} \rightarrow C,
$$
which since $C$ is fibrant extends along the trivial cofibration of
(\ref{ff4}) to a morphism
$$
\Upsilon (E)\times \overline{I}.
$$
This is exactly the lifting to a map
$$
E'\rightarrow \underline{Hom}(\overline{I},C)_{1/}(u,v)
$$
needed to establish the statement that the morphism
induced by $r$
$$
\underline{Hom}(\overline{I},C)_{1/}(u,v)
\rightarrow
\underline{Hom}(I,C)_{1/}(r(u),r(v))
$$
is an equivalence. This proves the theorem.
\eop
\end{parag}

\begin{corollary}
Suppose $f: U\rightarrow V$ is a morphism in a fibrant $n$-category $C$.
Then the $n$-category of morphisms $\overline{I}\rightarrow C$
restricting on $I\subset \overline{I}$ to $f$ is contractible.
\end{corollary}
{\em Proof:}
The $n$-category in question is just the fiber of the morphism in the theorem,
over the object $f\in \underline{Hom}(I,C)$.
\eop

The next corollary says that equivalences may be inverted with dependence on
parameters.

\begin{corollary}
\label{obbyobequiv}
Suppose $C$ is a fibrant $n$-category.
Suppose $\psi, \psi ' : A\rightarrow C$ are two morphisms and
suppose $f$ is a morphism from $\psi$ to $\psi '$. Suppose that for every
object $a\in A$ the induced morphism $f_a: \psi (a)\rightarrow \psi '(a)$
is an equivalence in $C$. Then $f$ is an equivalence considered as a
$1$-morphism in $\underline{Hom}(A,C)$.
\end{corollary}
{\em Proof:}
The morphism $f$ is a map
$$
f: A\times I \rightarrow C,
$$
which we can think of as a map
$$
f_1:A\rightarrow \underline{Hom}(I,C).
$$
From Theorem \ref{resttoIff} the morphism
$$
\underline{Hom}(\overline{I},C)\rightarrow \underline{Hom}(I,C)
$$
is a fibrant equivalence onto the full subcategory of invertible objects.
The hypothesis of the corollary says exactly that the morphism
$f_1$ lands in this full subcategory. Therefore it lifts to a morphism
$$
g: A\rightarrow \underline{Hom}(\overline{I}, C),
$$
in other words to
$$
A\times \overline{I} \rightarrow C
$$
or equally well
$$
\overline{I}\rightarrow \underline{Hom}(A,C).
$$
This shows that $f$ was an equivalence.
\eop

\bigskip

\numero{The definitions of direct and inverse limits}

One of the most useful tools in homotopy theory is the notion of homotopy limit
or ``holim''.  This can mean either direct or inverse limit and one of the two
is called a ``colimit'' but I don't know which one! So we'll call both
``limits'' and specify which one in context.  Our purpose is to define the
notions of inverse and direct limit in an $n$-category.

We always suppose that the target category $C$ is fibrant. When this is not
the case we first have to take a fibrant replacement (\ref{perspective1}).

\bigskip

\subnumero{Inverse limits}

Suppose $C$ is a fibrant $n$-category, and suppose $A$ is an $n$-category.
Suppose $\varphi : A \rightarrow C$ is a morphism. If $U \in C$ is an
object then we define
$$
Hom (U , \varphi ):= Hom (A, C)_{1/}(U_A, \varphi )
$$
where $U_A$ denotes the constant morphism with value
$U$.  If $V$ is another object of $C$ then we have a morphism
$$
C_{1/}(V , U )\rightarrow Hom (A, C)_{1/}(V_A,
U_A)
$$
and we use this to define
$$
Hom (V, U , \varphi ):= Hom (A, C)_{2/}
(V_A,U_A, \varphi )\times _{
Hom (A, C)_{1/}(V_A,
U_A)} C_{1/}(V , U )
$$
or more generally if $V ^0,\ldots , V ^p \in C_0$ we define
$$
Hom (V ^0, \ldots , V ^p, \varphi ):=
$$
$$
Hom (A, C)_{(p+1)/}
(V^0_A,\ldots , V ^p_A, \varphi )\times _{
Hom (A, C)_{p/}
(V^0_A,\ldots , V ^p_A)} C_{p/}(V ^0, \ldots ,
V ^p ).
$$
However we won't need this beyond $p=2$.

Notice now that since $C$ is fibrant, $Hom (A,C)$ is fibrant and in particular
an $n$-category, thus we get that the morphism
$$
Hom (V , U , \varphi )\rightarrow C_{1/}(V , U )
\times Hom (U , \varphi )
$$
is an equivalence.  On the other hand we have a projection
$$
Hom (V , U , \varphi )\rightarrow Hom (V , \varphi ).
$$
It is in this sense that we have a ``weak morphism'' from
$C_{1/}(V , U )
\times Hom (U , \varphi )$ to $Hom (V , \varphi )$.

\begin{definition}
\label{definverse}
We say that an object $U\in C_0$ together with element $f\in Hom (U , \varphi
)_0$ is an {\em inverse limit of $\varphi$} if for any $V\in C_0$ the
resulting weak morphism from $C_{1/}(V , U )$ to $Hom (V , \varphi )$
is an equivalence. To say this more precisely this means that the morphism
$$
Hom (V ,U , \varphi ) \times _{Hom (U , \varphi )} \{ f \}
\rightarrow Hom (V , \varphi )
$$
should be an equivalence. If such an inverse limit exists we say that {\em
$\varphi$ admits an inverse limit} (we will discuss uniqueness below).  If
any morphism $\varphi : A \rightarrow C$ from any $n$-category $A$ to $C$
admits an inverse limit then we say that {\em $C$ admits inverse limits}.
\end{definition}

\begin{parag}
\label{uniquenessinverse}
{\em Uniqueness:}
Suppose $f\in Hom (U , \varphi ) $ and $g \in Hom (V , \varphi )$
are two different inverse limits of $\varphi$.  Then the
inverse image of $g$ for the morphism
$$
Hom (V , U , \varphi ) \times _{Hom (U , \varphi )} \{ f \}
\rightarrow Hom (V , \varphi )
$$
is contractible. This gives a contractible $n$-category mapping to $Hom (V ,
U )$. We also have a contractible $n$-category mapping to $Hom (U
, V )$.  A similar argument with $p=3$ gives a contractible $n$-category
mapping to $Hom (V , U , V )$ which maps into the contractible things
for $V ,U$, for $U , V$ and for $V , V$. The image at the
end includes the identity. This shows that the composition of the morphisms in
the two directions is the identity. The same works in the other direction. This
shows that the essentially well defined morphisms $U \rightarrow V$ and
$V\rightarrow U$ are equivalences.
(The reader is challenged to find a nicer way of saying this!)
\end{parag}

\begin{parag}
\label{upsiloninverse1}
The condition of being an inverse limit may also be interpreted in terms of the
construction $\Upsilon$ described in the previous section.
To do this, start by noting that for an $n$-precat $E$ a morphism
$$
E\rightarrow Hom (U,\varphi )
$$
is the same thing as a morphism
$$
f:A\times \Upsilon (E)\rightarrow C
$$
such that $r_0(u)= U_A$ and $r_1(u)=\varphi $.
\end{parag}

\begin{parag}
\label{upsiloninverse2}
In view of the discussion
\ref{interpupsilon1}--\ref{interpupsilon3}, a morphism
$$
E\rightarrow Hom (V,U,\varphi )
$$
is the same thing as a morphism
$$
g:A\times [2](E)\rightarrow C
$$
with $r_0(g)=V_A$, $r_1(g)=U_A$ and $r_2(g)=\varphi$
and such that $r_{01}(g)$ comes from a morphism $\Upsilon (E)\rightarrow C$.
To see this, use the definition of $[2](E)$ by universal property
\ref{interpupsilon1}.

In a similar way using the description \ref{interpupsilon3}, a morphism
$$
E\rightarrow
Hom (V , U , \varphi ) \times _{Hom (U , \varphi )} \{ f \}
$$
is the same thing as a
morphism
$$
g: A\times \Upsilon ^2(E, \ast )\rightarrow C
$$
such that $r_2(g)= \varphi$ and $r_{01}(g)$ comes from a morphism
$g_{01}:\Upsilon (E)\rightarrow C$ with $r_0(g_{01})=V$ and
$r_1(g_{01})=U$.
\end{parag}

\begin{parag}
\label{upsiloninverse3}
Noting that the morphism
$$
Hom (V , U , \varphi ) \times _{Hom (U , \varphi )} \{ f \}
\rightarrow Hom (V , \varphi )
$$
is fibrant, it is an equivalence if and only if it satisfies the lifting
property for all cofibrations $E'\subset E$ (\ref{charequiv}).
\end{parag}

\begin{parag}
\label{upsiloninverse4}
Using the above
descriptions we can describe explicitly the lifting property
of the previous paragraph
and thus obtain the following characterization.
A morphism $f\in Hom (U,\varphi )$ is an inverse limit if and only if for every
morphism
$$
v:A\times \Upsilon (E) \rightarrow C
$$
with $r_0(v)= V_A$ for $V\in C_0$ and $r_1(v)= \varphi$, and for every
extension over $A\times \Upsilon (E')$ to a morphism
$$
w': A\times \Upsilon ^2(E', \ast )\rightarrow C
$$
with $r_{12}(w')= f$ and $r_{01}(w')$ coming from a morphism $z':\Upsilon
(E')\rightarrow C$ with $r_0(z')= V$ and $r_1(z')=U$, there exists a common
extension of these two: a morphism
$$
w: A\times \Upsilon ^2(E, \ast )\rightarrow C
$$
with $r_{12}(w)= f$ and $r_{01}(w)$ coming from a morphism $z:\Upsilon
(E')\rightarrow C$ with $r_0(z)= V$ and $r_1(z)=U$; such that the
restriction of $w$ to $A\times \Upsilon ^2(E', \ast )$ is equal to $w'$;
and such that $r_{02}(w)=v$.

This is the characterization we shall use in our proofs.
\end{parag}

\bigskip

\subnumero{Direct limits}

We obtain the notion of direct limit by ``reversing the arrows'' in the
above discussion.

Suppose $C$ is a fibrant $n$-category, and suppose $A$ is an $n$-category.
Suppose $\varphi : A \rightarrow C$ is a morphism. If $U \in C_0$ is an
object then we define
$$
Hom (\varphi , U  ):= Hom (A, C)_{1/}(\varphi , U_A)
$$
where again $U_A$ denotes the constant morphism with value
$U$.  If $V$ is another object of $C$ then we have a morphism
$$
C_{1/}(U ,V )\rightarrow Hom (A, C)_{1/}(
{U}_A, {V}_A)
$$
and we use this to define
$$
Hom (\varphi , U , V  ):= Hom (A, C)_{2/}
(\varphi , {U }_A,{V }_A, \varphi )\times _{
Hom (A, C)_{1/}({U }_A,
{V }_A)} C_{1/}( U ,V )
$$
or more generally if $V ^0,\ldots , V ^p \in C_0$ we define
$$
Hom (\varphi , V ^0, \ldots , V ^p):=
$$
$$
Hom (A, C)_{(p+1)/}
(\varphi , {V ^0}_A,\ldots , {V ^p}_A)\times _{
Hom (A, C)_{p/}
({V ^0}_A,\ldots , {V ^p}_A)} C_{p/}(V ^0, \ldots ,
V ^p ).
$$
Again we won't need this beyond $p=2$.

Notice now that since $C$ is fibrant, $Hom (A,C)$ is fibrant and in particular
an $n$-category, thus we get that the morphism
$$
Hom (\varphi , U , V )\rightarrow C_{1/}(U , V )
\times Hom (\varphi , U  )
$$
is an equivalence.  On the other hand we have a projection
$$
Hom (\varphi ,U,V)\rightarrow Hom (V , \varphi ).
$$
It is in this sense that we have a ``weak morphism'' from
$C_{1/}(U,V )
\times Hom ( \varphi ,U)$ to $Hom (\varphi ,V)$.

\begin{definition}
\label{defdirect}
We say that an element $f\in Hom (U , \varphi )_0$
is a {\em direct limit of $\varphi$} if for any $V \in C_0$ the resulting
weak morphism from $C_{1/}(U ,V )$ to $Hom (\varphi ,V )$ is
an equivalence. To say this more precisely this means that the morphism
$$
Hom (\varphi , U , V ) \times _{Hom (\varphi , U  )} \{ f
\}  \rightarrow Hom ( \varphi ,V )
$$
should be an equivalence. If such a direct limit exists we say that {\em
$\varphi$ admits an inverse limit}. Exactly the same discussion of uniqueness
as above
(\ref{uniquenessinverse})  holds here too.  If  any morphism $\varphi : A
\rightarrow C$ from any $n$-category $A$ to $C$ admits a direct limit then we
say that {\em $C$ admits direct limits}.
\end{definition}

\begin{parag}
\label{upsilondirect}
We have the following characterization analogue to \ref{upsiloninverse4}.
Again, this is the characterization which we shall use in the proofs.
It comes from considerations identical to
\ref{upsiloninverse1}--\ref{upsiloninverse3} which we omit here.

A morphism $f\in
Hom (\varphi ,U)$ is a direct limit if and only if for every morphism
$$
v:A\times \Upsilon (E) \rightarrow C
$$
with $r_0(v)= \varphi$ and $r_1(v)= V_A$ for $V\in C_0$, and for every
extension over $A\times \Upsilon (E')$ to a morphism
$$
w': A\times \Upsilon ^2(\ast , E')\rightarrow C
$$
with $r_{01}(w')= f$ and $r_{12}(w')$ coming from a morphism $z':\Upsilon
(E')\rightarrow C$ with $r_1(z')= U$ and $r_2(z')=V$, there exists a common
extension of these two: a morphism
$$
w: A\times \Upsilon ^2(\ast , E )\rightarrow C
$$
with $r_{01}(w)= f$ and $r_{12}(w)$ coming from a morphism $z:\Upsilon
(E')\rightarrow C$ with $r_1(z)= U$ and $r_2(z)=V$; such that the
restriction of $w$ to $A\times \Upsilon ^2(\ast ,E')$ is equal to $w'$;
and such that $r_{02}(w)=v$.
\end{parag}

\bigskip

\subnumero{Invariance properties}

\begin{proposition}
\label{invariance}
Suppose $f:A'\rightarrow A$ is an equivalence of $n$-categories and suppose
$C$ is a fibrant $n$-category. Suppose $\varphi : A\rightarrow C$ is a
morphism. Then the inverse (resp. direct) limit of $\varphi$ exists if and only
if the inverse (resp. direct) limit of $\varphi \circ f$ exists.

Suppose $\varphi$ and $\psi$ are morphisms from $A$ to $C$, and suppose
they are equivalent in $\underline{Hom}(A, C)$. Then the inverse (resp. direct)
limit of $\varphi$ exists if and only if the inverse (resp. direct) limit
of $\psi$ exists.

Finally suppose $g: C\rightarrow C'$ is an equivalence between fibrant
$n$-categories. Then the inverse (resp. direct) limit of $\varphi :
A\rightarrow C$ exists if and only if the inverse (resp. direct) limit of
$g\circ \varphi$ exists. In particular (combining with the previous paragraph)
$C$ admits inverse (resp. direct) limits if and only if $C'$ does.
\end{proposition}
{\em  Proof:}
There are several statements to prove so we divide the proof into several
paragraphs \ref{invariance1}--\ref{invariance9}.

\begin{parag}
\label{invariance1}
Suppose $f: A' \hookrightarrow A$ is a cofibrant equivalence
of $n$-categories. Suppose that $\varphi : A\rightarrow C$ is a morphism and
that  $u:\varphi \rightarrow U$ is a morphism from $\varphi$ to $U\in C$
which is  a direct limit. This corresponds to a diagram
$$
\epsilon : A \times \Upsilon (\ast )\rightarrow C
$$
and pullback by $f$ gives a diagram
$$
\epsilon ': A' \times \Upsilon (\ast )\rightarrow C.
$$
We claim that $\epsilon '$ is a direct limit (note that $\epsilon '$ is a
morphism from  $\varphi \circ f$ to $U$).  Suppose we are given
$$
u: A' \times \Upsilon (E) \rightarrow C
$$
and an extension over $E'\subset E$ to a diagram
$$
v_1:  A' \times \Upsilon ^2(\ast , E') \rightarrow C
$$
whose restriction to the edge $(01)$ is $\epsilon '$ and whose restriction to
the edge $(02)$ is $u$. Then we can first extend $v_1$ to a diagram
$$
A \times \Upsilon (\ast , E') \rightarrow C
$$
because
$$
A' \times \Upsilon ^2(\ast , E') \hookrightarrow
A \times \Upsilon ^2(\ast , E')
$$
is a trivial cofibration (and note also that we can assume that the extension
satisfies the relevant properties as in the definition of limit);
then we can
also extend our above morphism $u$ to a diagram
$$
A \times \Upsilon (E) \rightarrow C
$$
compatibly with the extension of $v_1$,
because the inclusion from the coproduct of
$A \times \Upsilon ( E')$ and
$A' \times \Upsilon (E)$ over $A' \times \Upsilon (E')$, into
$A \times \Upsilon ( E)$ is a trivial cofibration. Now we apply the limit
property of $\epsilon$ to conclude that there is an extension to a diagram
$$
v:
 A \times \Upsilon ^2(\ast , E) \rightarrow C.
$$
This restricts over $A'$ to a diagram of the form we would like, showing that
$\epsilon '$ is a direct limit.
\end{parag}

\begin{parag}
\label{invariance2}
Suppose that $f: A' \hookrightarrow A$ is a trivial cofibration and suppose
$\varphi: A\rightarrow C$ is a morphism to a fibrant $n$-category $C$,
and suppose now that we know that $\varphi \circ f$ has a limit
$$
\epsilon ': \varphi \circ f \rightarrow U
$$
for an object $U\in C$.
We claim that $\varphi$ has a limit.

The morphism $\epsilon '$ may be considered as
a diagram
$$
\epsilon ' : A' \times \Upsilon (\ast ) \rightarrow C.
$$
This extends along $A' \times \{ 0\}$ to
$$
\varphi : A \times \{ 0\} \rightarrow
C
$$
and it extends along $A'\times \{ 1\} $ to
$$
U_A: A \times \{ 1\} \rightarrow C.
$$
Putting these all together we obtain a morphism
$$
A\times \{ 0 \} \cup ^{A'\times \{ 0\} } A' \times \Upsilon (\ast )
\cup ^{A' \times \{ 1\} } A \times \{ 1\}
\rightarrow C.
$$
Since $A'\subset A$ is a trivial cofibration the morphism
$$
A\times \{ 0 \} \cup ^{A'\times \{ 0\} } A' \times \Upsilon (\ast )
\cup ^{A' \times \{ 1\} } A \times \{ 1\}
\rightarrow
A \times \Upsilon (\ast )
$$
is a trivial cofibration (applying \ref{explaincmc}, first part of CM4,
two times), so by the fibrant property of $C$ our morphism extends to a morphism
$$
\epsilon : A \times \Upsilon (\ast )\rightarrow C
$$
with the required properties of being constant along $A\times \{ 1\}$ and
restricting to $\varphi$ along $A\times \{ 0\}$. Thus we may write
$\epsilon : \varphi \rightarrow U$.

We claim that this map is a direct limit of $\varphi$. Given a diagram
$$
u:A\times \Upsilon (E)\rightarrow C
$$
going from $\varphi$ to a constant object $B$, the restriction $u'$ to
$A'\times \Upsilon (E)$ admits (by the hypothesis that $\epsilon '$ is a direct
limit) an extension to
$$
v': A'\times \Upsilon (\ast , E)\rightarrow C
$$
restricting along the edge $(01)$ to $\epsilon '$ and restricting along the edge
$(12)$ to the pullback of a diagram $\Upsilon (E)\rightarrow C$. Then using as
usual the fibrant property of $C$ and the fact that $A'\rightarrow A$ is a
trivial cofibration, we can extend $v'$ to a morphism
$$
v: A \times \Upsilon (\ast , E) \rightarrow C
$$
again
restricting along the edge $(01)$ to $\epsilon $, restricting along the edge
$(02)$ to our given diagram $u$, and restricting along the edge $(12)$ to the
pullback of a diagram $\Upsilon (E)\rightarrow C$.

If $E'\subset E$ and we are already given an extension $v_{E'}$ over
$A \times \Upsilon (\ast , E')$ then (as before, using the fibrant property of
$C$ applied to an appropriate cofibration) we can assume that our extension $v$
above restricts to $v_{E'}$.  This completes the proof that $\epsilon$ is a
direct limit, and hence the proof of the statement claimed for
\ref{invariance2}.
\end{parag}

\begin{parag}
\label{invariance3}
Now suppose $p:A' \rightarrow A$ is a trivial fibration.  Then there
exists a section $s:A\rightarrow A'$ (with $ps = 1_A$). Note that $s$ is a
trivial cofibration.  If $\varphi : A \rightarrow C$ is a morphism then
$$
(\varphi \circ p)\circ s = \varphi
$$
so applying the previous two paragraphs \ref{invariance1} and \ref{invariance2}
to the morphism $s$ we conclude that $\varphi$ admits a limit if and only if
$\varphi \circ p$ admits a limit.
\end{parag}

\begin{parag}
\label{invariance4}
Now suppose $f: A' \rightarrow A$ is any
equivalence between $n$-categories.
Decomposing $f= p\circ j$ into a composition
of  a trivial cofibration followed by a trivial fibration and applying
\ref{invariance1}, \ref{invariance2} to $j$ and \ref{invariance3} to $p$ we
conclude that a functor $\varphi : A\rightarrow C$ admits a direct
limit if and
only if $\varphi \circ f$ admits a direct limit. This proves the first
paragraph of Proposition \ref{invariance} for direct limits.
\end{parag}

\begin{parag}
\label{invariance5}
The proof of the first paragraph of \ref{invariance} for inverse limits is
exactly the same as the above.
\end{parag}

\begin{parag}
\label{invariance6}
Now we prove the second paragraph of the proposition. If $f,g: A\rightarrow C$
are two morphisms which are equivalent in $\underline{Hom}(A,C)$ (i.e. they are
homotopic) then there exists  a morphism
$\varphi : A\times \overline{I}\rightarrow C$ restricting to $f$ on $A\times \{
0\}$ and restricting to $g$ on $A\times \{ 1\}$. Applying the first paragraph
of the proposition (for either direct or inverse limits) to the two inclusions
$A\rightarrow A\times \overline{I}$ we find that $f$ admits a limit if and only
if $\varphi$ admits a limit and similarly $g$ admits a limit
if and only if $\varphi$ does---therefore $f$ admits a limit if and only if
$g$ admits a limit.
\end{parag}

\begin{parag}
\label{lifting}
Suppose $C\rightarrow C'$ is an equivalence between fibrant
$n$-categories. In general if $F'\subset F$ is any cofibration and if
$F'\rightarrow C$ is a morphism, then there exists an extension to
$F\rightarrow C$ if and only if the composed  morphism $F'\rightarrow C'$
extends over $F$. To see this, look at the (exactly commutative) diagram
$$
\begin{array}{ccc}
\underline{Hom}(F,C) & \rightarrow & \underline{Hom}(F',C) \\
\downarrow && \downarrow \\
\underline{Hom}(F,C') & \rightarrow & \underline{Hom}(F',C')
\end{array} .
$$
The horizontal arrows are fibrations and the vertical arrows are equivalences.
If an element $a\in \underline{Hom}(F',C)$ maps to something $b$ which is hit
from $c\in \underline{Hom}(F,C')$ then  there is $d\in
\underline{Hom}(F,C) $ mapping to something equivalent to $c$; thus the image
$e$ of $d$ in $\underline{Hom}(F',C)$ maps to something equivalent to
$b$. This implies (since the right vertical arrow is  an equivalence)
that $e$ is equivalent to $a$. Since the top morphism is fibrant, there
is another element $d'\in
\underline{Hom}(F,C) $ which maps directly to $a$.
\end{parag}

\begin{parag}
\label{invariance7}
Suppose still that $C\rightarrow C'$ is an equivalence between fibrant
$n$-categories.
Using the general lifting principle \ref{lifting} and the fact that the
property of being a limit is expressed in terms of extending morphisms
across certain cofibrtions $F'\subset F$, we conclude that a functor
$A\rightarrow C$ has a (direct or inverse) limit if and only if the composition
$A\rightarrow C'$  does. This proves the first sentence of the last paragraph of
the proposition.
\end{parag}

\begin{parag}
\label{invariance8}
If $f: C\rightarrow C'$ is an equivalence between fibrant $n$-categories and if
$C'$ admits direct (resp. inverse) limits then any functor $A\rightarrow C$
admits a direct (resp. inverse) limit by \ref{invariance7}.
\end{parag}

\begin{parag}
\label{invariance9}
Suppose on the
other hand that we know that $C$ admits direct (resp. inverse) limits. Suppose
that $\varphi : A \rightarrow C'$ is a functor. Since $f$ induces an
equivalence from $\underline{Hom}(A,C)$ to $\underline{Hom}(A,C')$ there is a
morphism $\psi : A \rightarrow C$ such that $f\circ \psi $ is equivalent to
$\varphi$ in $\underline{Hom}(A,C')$. By the second paragraph of the
proposition (proved in \ref{invariance6} above), $\varphi$ admits direct (resp.
inverse) limits if and only if $f \circ \psi$ does. Then by the first part of
the last paragraph proved in \ref{invariance7} above, $f\circ \psi$ admits
direct (resp. inverse) limits if and only if $\psi$ does. Now by hypothesis
$\psi$ has a direct (resp. inverse) limit, so $\varphi$ does too. This shows
that $C'$ admits direct (resp. inverse) limits.
\end{parag}

We have now completed the proof of Proposition \ref{invariance}.
\eop

\begin{parag}
\label{variance1}
We now start to look at variance properties in other situations.
Suppose
$h:C\rightarrow C'$ is a morphism between fibrant $n$-categories, and suppose
$\varphi : A \rightarrow C$ is a morphism. If
$$
u: \varphi \rightarrow U
$$
is a
direct limit then $h(u): \varphi \circ h \rightarrow h(U)$ is a morphism.
Suppose that $\varphi \circ h$ admits a direct limit
$$
v: \varphi \circ h \rightarrow V.
$$
Then by the limit property there is a factorization i.e. a diagram
$$
[v, w]:\varphi \circ h \rightarrow V \rightarrow h(U)
$$
whose third edge $(02)$ is $h(u)$.
We say that {\em the morphism $h$ commutes with the direct limit of $\varphi$}
if the direct limit of $\varphi \circ h$ exists and if the factorization
morphism $w: V\rightarrow h(U)$ is an equivalence.

Suppose that $C$ and $C'$ admit direct limits. We say that {\em the morphism $h$
commutes with direct limits} if $h$ commutes with the direct limit of any
$\varphi : A \rightarrow C$ in the previous sense.
\end{parag}

\begin{parag}
\label{variance2}
We have similar definitions for inverse limits, which we repeat for the record.
Suppose again that
$h:C\rightarrow C'$ is a morphism between fibrant $n$-categories, and suppose
$\varphi : A \rightarrow C$ is a morphism. If
$$
u: U\rightarrow \varphi
$$
is an inverse limit then $h(u): h(U)\rightarrow \varphi \circ h $ is a
morphism. Suppose that $\varphi \circ h$ admits an inverse limit
$$
v: V\rightarrow \varphi \circ h .
$$
Then by the limit property there is a factorization i.e. a diagram
$$
[w,v]:h(U)\rightarrow V\rightarrow \varphi \circ h
$$
whose third edge $(02)$ is $h(u)$.
We say that {\em the morphism $h$ commutes with the inverse limit of $\varphi$}
if the inverse limit of $\varphi \circ h$ exists and if the factorization
morphism $w:h(U)\rightarrow V$ is an equivalence.

Suppose that $C$ and $C'$ admit inverse limits. We say that {\em the
morphism $h$
commutes with inverse limits} if $h$ commutes with the inverse limit of any
$\varphi : A \rightarrow C$ in the above sense.
\end{parag}

\bigskip

\subnumero{Behavior under certain precat inverse limits of $C$}

We will now study certain situations of what happens when we take fiber
products or other inverse limits (here we mean inverse limits in the category
of $n$-precats) of the target $n$-category $C$. We study what happens to
inverse limits in $C$. We could also say the same things about direct limits in
$C$ but the inverse limit case is the one we need, so we state things there and
leave it to the reader to make the corresponding statements for direct limits.

\begin{lemma}
\label{dirprod}
Suppose $\{ C_i\} _{i\in S}$ is a collection of fibrant $n$-categories indexed
by a set $S$. Let $C= \prod _{i\in S} C_i$ and suppose $\varphi = \{ \varphi
_i\}$ is a morphism from $A$ to $C$. Suppose that the $\varphi _i$ admit
inverse limits
$$
u_i : U_i \rightarrow \varphi _i
$$
in $C_i$. Then $U= \{ U_i\}$ is an object of $C$ and we have a morphism
$$
u: U\rightarrow \varphi
$$
composed of the factors $u_i$. This morphism is an inverse limit of $\varphi$
in $C$.
\end{lemma}
{\em Proof:}
The property that $u$ be an inverse limit consists of a collection of extension
properties that have to be satisfied. The morphisms $u_i$ admit the
corresponding extensions and putting these together we get the required
extensions for $u$.
\eop

\begin{lemma}
\label{fiprod}
Suppose $f:C\rightarrow D$ and $g:E\rightarrow D$ are morphisms of
fibrant $n$-categories with $f$ fibrant.  Suppose that $\varphi : A \rightarrow
C\times _DE$ is a morphism such that the component morphisms $\varphi _C:
A\rightarrow C$,  $\varphi _D:A\rightarrow D$ and $\varphi _E: A\rightarrow E$
have inverse limits $\lambda _C$, $\lambda _D$ and $\lambda _E$ respectively.
Suppose furthermore that $f$ and $g$ preserve these inverse limits, which means
that the projections of $\lambda _C$ and $\lambda _E$ into $D$ are equivalent
(as objects with morphisms to $\varphi _D$) to $\lambda _D$.  Then we may (by
changing the $\lambda _C, \lambda _D, \lambda _E$ by equivalences) assume that
$\lambda _C$ and $\lambda _E$ project to $\lambda _D$; and the resulting object
$\lambda \in C\times _DE$ is an inverse limit of $\varphi$.
\end{lemma}
{\em Proof:}
Set $\lambda '_E:=\lambda _E$ and let $\lambda '_D:= g(\lambda _E)$ be the
projection to $D$. Note that by hypothesis $\lambda '_D$ is an inverse limit
of $\varphi _D$.  Now $\lambda _C$ (considered as an object with
morphism to $\varphi _C$) projects in $D$ to something equivalent to
$\lambda _D$ and hence equivalent to $\lambda '_D$ (equivalence of the diagrams
including the morphism to $\varphi _D$). Since $f$ is a fibrant morphism, we
can modify $\lambda _C$ by an equivalence, to obtain $\lambda ' _C$ projecting
directly to $\lambda ' _D$.  Note that the equivalent $\lambda '_C$ is again an
inverse limit of $\varphi _C$.  Together these give an element $\lambda \in
C\times DE$ with a map
$$
u: \lambda \rightarrow \varphi ,
$$
and we claim that $u$ is an inverse limit.  Suppose
$F'\subset F$ is a cofibration of $n-1$-precats and suppose
$$
v: V \stackrel{F}{\rightarrow} \varphi
$$
is any $F$-morphism (i.e. a diagram
$$
A\times \Upsilon (F) \rightarrow C\times _DE
$$
restricting on $A\times \{ 0\}$ to the constant $V_A$ and restricting on
$A\times \{ 1\}$ to $\varphi$) provided with an extension over $F'$ to
a diagram
$$
w': A\times \Upsilon ^2(F, \ast )\rightarrow C\times _DE
$$
restricting on $(02)$ to $v'$ (the restriction of $v$ to $F'$) and on $(12)$ to
$u$.  We look for an extension of $w'$ to a diagram
$$
w: A\times \Upsilon ^2(F, \ast )\rightarrow C\times _DE
$$
restricting on $(02)$ to $v$ and on $(12)$ to $u$.
Denoting with subscripts the
components in $C$, $D$ and $E$, we have that the pairs $(v_C,w'_C)$ and $(v_E,
w'_E)$ admit extensions $w_C$ and $w_E$ respectively. The projections of these
extensions in $D$ give diagrams which we denote
$$
w_{C/D}, w_{E/D} : A\times \Upsilon ^2(F, \ast )\rightarrow D,
$$
both restricting on $(02)$ to $v_D$ and on $(12)$ to $u_D$, and extending
$w'_D$.  Applying again the limit property for $u_D$ to the cofibration
$$
F\times \{ 0\} \cup ^{F' \times \{ 0\} } F' \times \overline{I}
\cup ^{ F' \times \{ 1\} } F \times \{ 1\} \hookrightarrow F\times \overline{I}
$$
we find that there is a diagram
$$
z_D: A \times \Upsilon  ^2(F\times \overline{I} , \ast ) \rightarrow D
$$
giving a homotopy between $w_{C/D}$ and $w_{E/D}$.

Notice that
this is a homotopy
in Quillen's sense \cite{Quillen} because the diagram
$$
\Upsilon ^2(F, \ast )\tworightarrows \Upsilon ^2(F\times \overline{I}, \ast )
\rightarrow \Upsilon ^2(F, \ast )
$$
is of the form used by Quillen cf \ref{homotopic1}. Such a homotopy can
be changed into one of the more classical form
$$
A\times \Upsilon ^2(F, \ast )\times \overline{I} \rightarrow D
$$
because the relation of homotopy in Quillen's sense is the same as the relation
of equivalence of morphisms, which in turn is the same as existence of a
homotopy for the ``interval'' $\overline{I}$. In fact we don't use this remark
here but we use it with $D$ replaced by $C$, below.

Now apply the lifting property for the morphism $C\rightarrow D$, for the
above map $z_D$, with respect to the trivial cofibration
$$
A\times \Upsilon ^2(
F\times \{ 0\} \cup ^{F' \times \{ 0\} } F' \times \overline{I}, \ast )
$$
$$
\rightarrow
A\times \Upsilon ^2(F\times \overline{I} , \ast ).
$$
We get a morphism
$$
z_C: A \times \Upsilon  ^2(F\times \overline{I} , \ast ) \rightarrow C
$$
providing a homotopy between $w_C$ and a new morphism $w^n_C$ which
projects into $D$ to $w_{E/D}$.  The new $w^n_C$ is again a solution of the
required extension problem (or to be more precise we can impose conditions
on our lifting $z_C$ to insure that this is the case. The fact that it projects
to $w_{E/D}$ means that the pair $w=(w^n_C, w_E)$ is a solution of the required
extension problem to show that $u$ is an inverse limit. This completes the
proof.
\eop

\begin{lemma}
\label{invlim}
Suppose $C_i$ is a collection of fibrant $n$-categories for $i=0,1,2,\ldots$
and suppose that $f_i : C_i \rightarrow C_{i-1}$ are fibrant morphisms. Let $C$
be the inverse limit of this system of $n$-precats. Then $C$ is a fibrant
$n$-category. Suppose that we have $\varphi : A \rightarrow C$ projecting to the
$\varphi _i : A \rightarrow C_i$ and suppose that the $\varphi _i$ admit
inverse limits $u_i : U_i \rightarrow \varphi _i$. Suppose finally that
the $f_i$ commute with the inverse limits of the $\varphi _i$. Then $\varphi$
admits an inverse limit and the projections $C\rightarrow C_i$ commute with the
inverse limit of $\varphi$.
\end{lemma}
{\em Proof:}
The fibrant property of $C$ may be directly checked by producing liftings
of trivial cofibrations.

First we construct  a morphism $u: U\rightarrow \varphi$ projecting in each
$C_i$ to an inverse limit of $\varphi _i$.  To do this, note by \ref{invariance}
that it suffices to have $u$ project to a morphism equivalent to $u_i$. On the
other hand, by the hypothesis that the $f_i$ commute with the inverse limits
$u_i$, we have that $f_i (u_i)$ is an inverse limit of $\varphi _{i-1}$.
In particular $f_i(u_i)$ is equivalent to $u_{i-1}$ as a diagram from $A\times
\Upsilon (\ast )$ to $C_{i-1}$.  The morphism from such diagrams in $C_i$,
to such diagrams in $C_{i-1}$, is fibrant (since it comes from $f_i$ which is
fibrant). Therefore we can change $u_i$ to an equivalent diagram with
$f_i(u_i)=u_{i-1}$. Do this successively for $i=1,2,\ldots$, yielding a system
of morphisms $u_i$ with $f_i(u_i)=u_{i-1}$. These now form a morphism
$$
u: U\rightarrow \varphi .
$$
We claim that $u$ is an inverse limit of $\varphi$. Suppose
$E'\subset E$ is an inclusion of $n-1$-precats and suppose
$$
w:W\stackrel{E}{\rightarrow}\varphi
$$
is an $E$-morphism (i.e. a diagram of the form
$$
A\times \Upsilon (E)\rightarrow C
$$
being constant equal to $W$ on $A\times \{ 0\}$),
provided over $E'$ with an extension to a diagram
$$
[v',u]:W\stackrel{E'}{\rightarrow} U \rightarrow \varphi
$$
(i.e. a morphism
$$
A\times \Upsilon ^2 (E,\ast ) \rightarrow C
$$
restricting
to $u$ on the second edge $(12)$ and restricting to $w|_{E'}$ on the third edge
$(02)$).     We would like to extend this to a diagram $[v,u]$ giving
$w$ on the third edge.  Let $w_i$ (resp. $v'_i$) be the projections of these
diagrams in $C_i$.  These admit extensions $v_i$. The projection of $v_i$ to
$C_{i-1}$ is an extension of the desired sort for $w_{i-1}$ and $v'_{i-1}$. The
extensions
$v_i$ are unique up to equivalence---which means a diagram
$$
A\times \Upsilon ^2 (E,\ast )\times \overline{I} \rightarrow C_i
$$
satisfying appropriate boundary conditions---and from this and our usual sort
of argument constructing  a trivial cofibration (this occurs several times
below) then making use of the fibrant property of $f_i$, we conclude that $v_i$
may be modified by an equivalence so that it projects to $v_{i-1}$. As before,
do this successively for $i=1,2,\ldots $ to obtain a system of extensions $v_i$
with $f_i(v_i)=v_{i-1}$. This system corresponds to an extension $v$ of the
desired sort for $w$ and $v'$.  This shows that the morphism $u$ is an inverse
limit.

Note from our construction the projection of the inverse limit $u$ to
$C_i$ is an inverse limit $u_i$ for $\varphi _i$ so the projections
$C\rightarrow C_i$ commute with the inverse limit of $\varphi$.
\eop

The application of the above results which we have in mind is the following
theorem.

\begin{theorem}
\label{variation}
Suppose $C$ is a fibrant $n$-category and $B$ is an $n$-precat. Then if $C$
admits inverse (resp. direct) limits, so does the fibrant
$n$-category $\underline{Hom}(B,C)$. The morphisms of functoriality
for $B'\rightarrow B$ commute with inverse (resp. direct) limits.
\end{theorem}
{\em Proof:}
Suppose $\varphi : A\rightarrow \underline{Hom}(B, C)$. We will construct an
inverse limit $\lambda \in \underline{Hom}(B, C)$
such that for any $b\in B$ the restriction $\lambda (b)$ is equivalent (via the
natural morphism) to the inverse limit of $\varphi (b): A\rightarrow C$.
This condition implies that the restriction morphism for any $B'\rightarrow B$
commutes with the inverse limit.  In effect, there is a morphism from the
inverse limit over $B$  (pulled back to $B'$) to the inverse limit over $B'$,
and this morphism is an equivalence over every object in $B'$ by the condition,
which implies that it is an equivalence by Lemma \ref{obbyobequiv}.

\begin{parag}
\label{pfvar1}
The first remark is that if
$B\subset B'$ and $B\subset B''$ are cofibrations of $n$-precats such that
$\underline{Hom}(B, C)$, $\underline{Hom}(B', C)$, and
$\underline{Hom}(B'', C)$ admit inverse limits complying with the above
condition, then  $\underline{Hom}(B'\cup ^B B'', C)$ admits inverse limits again
complying with the above condition. To see this we apply  Lemma \ref{fiprod}.
The only thing that we need to know is that the restriction of the inverse
limits is again equivalent to the inverse limit in $\underline{Hom}(B, C)$.
This follows from the fact that there is a morphism which, thanks to the
condition given at the start of the proof, is an equivalence for each object of
$B$---therefore it is an equivalence.
\end{parag}

\begin{parag}
\label{pfvar2}
The next remark is that weak equivalences of $n$-precats $B$ are turned into
equivalences of the $\underline{Hom}(B,C)$.  The morphisms
$$
h(1, M') \cup ^{\ast} \ldots \cup ^{\ast }h(1,M')\rightarrow h(m,M')
$$
are weak equivalences. By the previous remarks if we know the theorem (with the
condition of the first paragraph of the proof) for $h(1,M')$ then we get it for
any $h(m,M')$ (again with the condition of the first paragraph).

As pointed out at the start of the  proof, morphisms of restriction between any
$B$'s for which we know that the limits exist (and satisfying the condition of
the first paragraph), commute with the limits.  In particular when we apply
\ref{fiprod} and \ref{invlim} the hypotheses about commutation with the limits
will hold.

Now any $n$-precat may be expressed as a direct union of pushouts of the
$h(m,M')$. The  pushouts in question may be organized into a countable direct
union of pushouts each of which is adding a disjoint  direct sum; and the
addition will be of a direct sum of things of the form $h(m,M')$ added along
their boundary.  Taking $\underline{Hom}$ into $C$ transforms this expression
into an inverse limit indexed by the natural numbers, of fiber products of
terms which are direct products of things of the form $\underline{Hom}(h(m,M'),
C)$.  Applying Lemmas \ref{dirprod}, \ref{fiprod} and \ref{invlim} we find that
if we know the existence of limits for $h(1,M')$ (hence $h(m,M')$)---as
always with the additional condition of the first paragraph of the proof---then
we get existence of limits for $\underline{Hom}(B,C)$ for any $n$-precat $B$.
\end{parag}

\begin{parag}
\label{pfvar3}
Therefore it suffices to prove the theorem for $B= h(1,M')$. This is
more generally of the form $\Upsilon E$ (in this case $E=h(M')$).   Thus it now
suffices to prove the theorem for $B=\Upsilon E$.

Suppose we have a morphism $\varphi : A\times \Upsilon E \rightarrow C$.
Let $\varphi (0)$ (resp. $\varphi (1)$) denote the restriction of $\varphi$ to
$A\times \{ 0\} $ (resp.  $A\times \{ 1\} $).
Let $(\lambda (0),\epsilon (0))$ and $(\lambda (1),\epsilon (1))$ denote inverse
limits of $\varphi (0)$ and $\varphi (1)$.
The morphism $\varphi$ corresponds to a morphism
$$
E\rightarrow \underline{Hom}(A,C)_{1/}(\varphi (0), \varphi (1)).
$$
We can lift this together with $\epsilon (0)$ to a morphism
$$
E\rightarrow \underline{Hom}(A,C)_{2/}(\lambda (0), \varphi (0), \varphi (1)).
$$
The resulting $E\rightarrow \underline{Hom}(A,C)_{1/}(\lambda (0), \varphi (1))$
exends to
$$
E\rightarrow \underline{Hom}(A,C)_{2/}(\lambda (0)_A, \lambda (1)_A,
\varphi (1))
$$
projecting to $\epsilon (1)$ on the second edge,
by the limit property for $\epsilon (1)$.
The first edge of this comes from a morphism $\lambda : \Upsilon E \rightarrow
C$. Noting that the product $\Upsilon E
\times I$ is a pushout of two triangles the above morphisms glue
together to give a morphism
$$
\Upsilon E \times I \rightarrow \underline{Hom}(A,C),
$$
in other words we get a morphism $\epsilon$ from $\lambda$ to $\varphi$
considered as families over $A$ with values in $\underline{Hom}(\Upsilon E, C)$.
\end{parag}

\begin{parag}
\label{pfvar4}
To finish the proof we just have to prove that $\epsilon$ is an inverse limit
of $\varphi$.  Suppose we have $\mu \in \underline{Hom}(\Upsilon E, C)$
and suppose given a morphism
$f: \Upsilon F \times A \times \Upsilon E\rightarrow C$ restricting
over $0^F$ to $\mu _A$ and over $1^F$ to $\varphi$ (here $0^F$ and $1^F$
denote the endpoints of $\Upsilon F$ and we will use similar notation for
$E$).  This extends to morphisms
$$
f'(0^E): \Upsilon ^2(F, \ast )\times A \rightarrow C,
$$
$$
f'(1^E): \Upsilon ^2(F, \ast )\times A \rightarrow C,
$$
by the limit properties of $\lambda (0)$ and $\lambda (1)$ (the above morphisms
restricting to $\epsilon (0)$ and $\epsilon (1)$ on the second edges).

Now we try to extend to a morphism on all of
$\Upsilon ^2(F, \ast )\times A \times \Upsilon E\rightarrow C$. For this we use
notation of the form $(i,j)$ for the objects of $\Upsilon ^2(F,\ast )\times
\Upsilon E$, where $i=0,1,2$ (objects in $\Upsilon ^2(F, \ast )$ and $j= 0,1$
(objects in $\Upsilon E$).  We are already given maps defined over
the triangles
$(012, 0)$ and $(012,1)$ (these are $f'(0^E)$ and $f'(1^E)$), as well as
overthe
squares $(02, 01)$ (our given map $f$) and  $(12, 01)$ (the map $\epsilon$).
First extend using the fibrant property of $C$ to a map on the tetrahedron
$$
(0,0)(1,0)(2,0)(2,1).
$$
Then extend again using the fibrant property of
$C$ to a map on the tetrahedron
$$
(0,0)(1,0)(1,1)(2,1).
$$
Here note that on the face $(0,0)(1,0)(1,1)$ the map is chosen first as
coming from a map $\Upsilon ^2(F,E)\rightarrow C$.  Finally we have to
find an extension over the tetrahedron
$$
(0,0)(0,1)(1,1)(2,1).
$$
Again we require that the map on the first face
$(0,0)(0,1)(1,1)$ come from a map
$$
\Upsilon ^2(E,F)\rightarrow C.
$$
Our problem at this stage is that the map is already specified on all of
the other faces, so we can't do this using the fibrant property of $C$ (the
face that is missing is not of the right kind).  Instead we have to use the
limit property of $\epsilon (1^E)$.

The limit condition on $\lambda (1^E)$ means that the morphism
$$
\underline{Hom}^{\epsilon (1^E)}(\mu (0^E), \mu (1^E), \lambda (1^E), \varphi )
\rightarrow
\underline{Hom}(\mu (0^E), \mu (1^E), \varphi )
$$
is an equivalence.  The morphisms
$$
\underline{Hom}^{\epsilon (1^E)}(\mu (0^E), \lambda (1^E), \varphi )
\rightarrow
\underline{Hom}(\mu (0^E), \varphi )
$$
and
$$
\underline{Hom}^{\epsilon (1^E)}(\mu (1^E), \lambda (1^E), \varphi )
\rightarrow
\underline{Hom}(\mu (1^E), \varphi )
$$
are equivalences too.  This implies (in view of the fact that the edges
containing $\epsilon (1^E)$ are fixed) that the morphism
$$
\underline{Hom}^{\epsilon (1^E)}(\mu (0^E), \mu (1^E), \lambda (1^E), \varphi )
\rightarrow
$$
$$
\underline{Hom}^{\epsilon (1^E)}(\mu (0^E), \lambda (1^E), \varphi )
\times _{
\underline{Hom}(\mu (0^E), \varphi )}
$$
$$
\underline{Hom}(\mu (0^E), \mu (1^E), \varphi )\times _{
\underline{Hom}(\mu (1^E), \varphi )}
$$
$$
\underline{Hom}^{\epsilon (1^E)}(\mu (1^E), \lambda (1^E), \varphi )
$$
is an equivalence. This exactly implies that the restriction to the shell that
we are interested in is an equivalence. The fact that this equivalence is a
fibration implies that it is surjective on objects, giving finally the
extension that we need.

In the relative case where we are already given an extension over $F'\subset F$,
one can choose our extension in a compatible way (adding on the part concerning
$F'$ in the above argument doesn't change the properties of the relevant
morphisms being trivial cofibrations).
\eop
\end{parag}

\begin{corollary}
\label{constant}
Suppose $\varphi : A\rightarrow C$ is a morphism from an $n$-precat to a
fibrant $n$-category $C$ admitting inverse limits, and suppose that $B$ is an
$n$-precat. Let  $\varphi _B: A\rightarrow \underline{Hom}(B,C)$ denote the
morphism constant along $B$. Suppose that $\varphi$ admits an inverse limit
$u:U\rightarrow \varphi$. Then $u_B: U_B \rightarrow \varphi _B$ (the pullback
of $u$ along $B\rightarrow \ast$) is an inverse limit of $\varphi _B$.
\end{corollary}
{\em Proof:}
This is just commutativity for pullbacks for the morphism
$B\rightarrow \ast$.
\eop

\begin{parag}
\label{usevar}
We can use the result of the previous theorem to obtain the variation of the
limit depending on the family. Suppose $C$ is a fibrant $n$-category in which
inverse limits exist, and suppose $A$ is an $n$-precat. Let $B=
\underline{Hom}(A, C)$. We have a tautological morphism
$$
\zeta : A\rightarrow \underline{Hom}(B,C).
$$
By the previous theorem, limits exist in $\underline{Hom}(B,C)$. Thus we obtain
the limit of $\zeta$ which is an element of $\underline{Hom}(B,C)$: it is a
morphism $\lambda $ from $B=\underline{Hom}(A,C)$ to $C$, which is the morphism
which to $\varphi \in \underline{Hom}(A,C)$ associates $\lambda (\varphi )$
which is the limit of $\varphi$.

The same remark holds for direct limits.
\end{parag}

\begin{theorem}
\label{commute}
Suppose $A$, $B$ and $C$ are $n$-categories. Suppose $F: A\times B \rightarrow
C$ is a functor.  Then letting $\psi : A\rightarrow \underline{Hom}(B,C)$
denote the corresponding functor, suppose that $\psi$ admits an inverse limit
$\lambda \in \underline{Hom}(B,C)$.  Suppose now that $\lambda$ (considered as
a morphism $B\rightarrow C$) admits an inverse limit $\mu \in C$.  Then $\mu$
is an inverse limit of $F: A\times B\rightarrow C$.  In particular if the
intermediate limits exist going in the other direction then the composed limits
are canonically equivalent.  Thus if $C$ admits inverse limits then
inverse limits commute with each other.
\end{theorem}
{\em Proof:}
The general proof is left to the reader. In the case $C=nCAT'$ this will be
easy to see from our explicit construction of the limits.
\eop

\numero{Limits in $nCAT'$}

Let $nCAT \hookrightarrow nCAT'$ be a trivial cofibration to a fibrant
$n+1$-category.

\begin{theorem}
\label{inverse}
The fibrant $n+1$-category $nCAT'$ admits inverse limits.
\end{theorem}

The rest of this section is devoted to the proof. As a preliminary remark
notice that by \ref{invariance} the statement doesn't depend on which
choice of $nCAT'$ we make.  We also remark, in the realm of set-theoretic
niceties, that the statement means that $nCAT'$ (an $n+1$-category
composed of classes) admits inverse limits indexed by any $n+1$-category
composed of sets.  To be more precise our proof will show that if we restrict
to a subcategory of $nCAT'$ of $n$-categories represented in a certain set of
fixed cardinality $\alpha$, then the inverse limit indexed by $A$ exists if
$\alpha$ is infinite and at least equal to
$2^{\# A}$, also at least equal to what is needed for making the fibrant
replacement $nCAT\subset nCAT'$.

\bigskip

\subnumero{Construction of the limit}

\begin{parag}
\label{subscriptA}
Suppose $A$ is an $n+1$-category. If $B$ is a fibrant $n$-category we have
denoted by $B_A$ the  constant morphism $A\rightarrow nCAT$ with value $B$,
considered as  a morphism $A\rightarrow nCAT'$.
\end{parag}

\begin{parag}
\label{lambda}
We now give the construction of the inverse limit.  Suppose $\varphi : A
\rightarrow nCAT'$ is a morphism.  We define an $n$-category
$$
\lambda := \underline{Hom}(A, nCAT')_{1/}(\ast _A, \varphi ).
$$
This has the universal property that for any $n$-precat $B$,
$$
Hom (B, \lambda )= Hom ^{\ast _A, \varphi }(A\times \Upsilon B, nCAT').
$$
The notation on the right means the fiber of the map
$$
(r_0,r_1): Hom (A\times \Upsilon B, nCAT')\rightarrow
Hom (A, nCAT')\times Hom(A, nCAT')
$$
over $(\ast _A, \varphi )$.

The problem below will be to prove that $\lambda$ is an inverse limit of
$\varphi$.
\end{parag}

\bigskip

\subnumero{Diagrams}

\begin{parag}
\label{diag1}
Suppose $C$ is an $n+1$-precat. Then for $n$-precats $E_1,\ldots , E_k$
we define
$$
Diag (E_1,\ldots , \underline{E_i},\ldots , E_k; C)
$$
to be the $n$-precat which represents the functor
$$
F \mapsto Hom (\Upsilon ^k (E_1,\ldots , E_i\times F,\ldots , E_k), C).
$$
We establish some properties.
\end{parag}

\begin{parag}
\label{diag1.5}
The first remark is that
$Diag (E_1,\ldots , \underline{E_i},\ldots , E_k; C)$ decomposes as a
disjoint union over all pairs $(a,b)$ where
$$
a: \Upsilon ^{i-1}(E_1,\ldots , E_{i-1})\rightarrow C
$$
and
$$
b: \Upsilon ^{k-1-i}(E_{i+1}, \ldots , E_k)\rightarrow C
$$
are the restrictions to the first and last faces separated by the $i$-th edge.
Employ the notation
$$
Diag^{a,b} (E_1,\ldots , \underline{E_i},\ldots , E_k; C)
$$
for the subobject restricting to a given $a$ and $b$. If we don't wish to
specify $b$ for example, then denote this by the superscript $Diag^{a,\cdot}$.

In particular note that we can decompose into a disjoint union over the
$k+1$-tuples of objects which are the images of the vertices $0,\ldots , k$
(these objects are all specified either as a part of $a$ or as a part of $b$).
\end{parag}

\begin{parag}
\label{diag1.6}
In case $C=nCAT$ we have
$$
Diag ^{a,b}(E_1,\ldots , \underline{E_i},\ldots , E_k; nCAT)
= \underline{Hom}(U_{i-1} \times E_i, U_i)
$$
where $U_j$ are the fibrant $n$-categories which are the images of the vertices
$$
j\in \Upsilon ^k(E_1,\ldots , E_k)
$$
by the maps $a$ (if $j\leq i-1$) or $b$
(if $j\geq i$).
\end{parag}

\begin{parag}
When we are only interested in the set of objects, it doesn't matter which
$E_i$ is underlined and we denote by
$$
Diag (E_1,\ldots , E_k; C) = Hom (\Upsilon ^k (E_1,\ldots , E_k), C)
$$
this set of objects.  We can put a superscript $Diag ^{a,b}$ here if we want
(with the obvious meaning as above).  The edge $i$ dividing between $a$ and $b$
should be understood from the data of $a$ and $b$.
\end{parag}

\begin{parag}
\label{quasifib}
We need a way of understanding the statement that $nCAT'$ is a fibrant
replacement for $nCAT$. In order to do this we will  use the
following property of $nCAT$ which shows that in some sense it is close to
being fibrant.

We say that an $n+1$-category $C$ is {\em quasifibrant} if for
any sequence of objects $x_0, \ldots , x_p$
the morphism
$$
C_{p/}(x_0,\ldots , x_p)\rightarrow
C_{(p-1)/}(x_0,\ldots , x_{p-1})\times _{C_{(p-1)/}(x_1,\ldots , x_p)}
C_{(p-2)/}(x_1,\ldots , x_{p-1})
$$
is a fibration of $n$-categories. Note inductively that the morphisms involved
in the fiber product here are  themselves fibrations, and we get that the
projections
$$
C_{p/}(x_0,\ldots , x_p)\rightarrow
C_{(p-1)/}(x_0,\ldots , x_{p-1})
$$
and
$$
C_{p/}(x_0,\ldots , x_p)\rightarrow
C_{(p-1)/}(x_1,\ldots , x_{p})
$$
are fibrations.
\end{parag}

\begin{parag}
\label{quasifib1.4}
The condition that $C$ is an $n+1$-category implies that the morphism
in the definition of quasifibrant, is an equivalence whenever $p\geq 2$.
Thus if $C$ is quasifibrant, the morphism in question is actually a fibrant
equivalence.
\end{parag}

\begin{parag}
\label{quasifib1.5}
If $C'$ is a fibrant $n+1$-category then it is quasifibrant. This is because
the morphisms (in the notation of \ref{interpupsilon1})
$$
[p-1](E)\cup ^{[p-2](E)}[p-1](E)\rightarrow [p](E)
$$
are trivial cofibrations.
\end{parag}

\begin{parag}
\label{quasifib2}
The $n+1$-category $nCAT$ is easily seen to be quasifibrant: the morphisms
in question are actually isomorphisms for $p\geq 2$ and for $p=1$
they are just projections from the $\underline{Hom}(A_0,A_1)$---which are
fibrant---to $\ast$.
\end{parag}

\medskip

We now have two claims which allow us to pass between something quasifibrant
such as $nCAT$ and its fibrant completion.

\begin{parag}
\label{diag2}
First of all, if $C$ is quasi-fibrant (\ref{quasifib}) then
$Diag (E_1,\ldots , \underline{E_i},\ldots , E_k; C)$ is fibrant.
Furthermore in this case for cofibrations $E'_j\hookrightarrow E_j$
the morphism
$$
Diag (E'_1,\ldots , \underline{E'_i},\ldots , E'_k; C)
\rightarrow
Diag (E_1,\ldots , \underline{E_i},\ldots , E_k; C)
$$
is fibrant.
\end{parag}

\begin{parag}
\label{diag3}
Secondly, if $C$ quasifibrant (\ref{quasifib})  and if
$C\rightarrow C'$ is an equivalence to a fibrant $C'$ then the morphism
$$
Diag ^{a,b}(E_1,\ldots , \underline{E_i},\ldots , E_k; C)
\rightarrow
Diag ^{a',b'}(E_1,\ldots , \underline{E_i},\ldots , E_k; C')
$$
is an equivalence of fibrant $n$-categories. Here $a,b$ are fixed as in
\ref{diag1.5}, and $a',b'$ denote the images in $C'$.

{\em Caution:} it is essential to restrict to the components for a fixed
$a,b$ coming from $C$.
\end{parag}

\begin{parag}
\label{diag3.5}
Before getting to the
proofs of \ref{diag2} and \ref{diag3},
we discuss diagrams in a quasifibrant $C$. A morphism
$$
u:\Upsilon ^k(E_1,\ldots , E_k)\rightarrow C
$$
may be described inductively as triple $u=(\tilde{u}, u^-, u^+)$
where
$$
u^-= \Upsilon ^{k-1}(E_2,\ldots , E_k)\rightarrow C
$$
and
$$
u^+:\Upsilon ^{k-1}(E_1,\ldots , E_{k-1})\rightarrow C,
$$
are morphisms which agree on $\Upsilon ^{k-1}(E_1,\ldots , E_{k-1})$,
and where
$$
\tilde{u}: E_1\times \ldots \times E_k \rightarrow C_{k/}(x_0,\ldots , x_k)
$$
is a lifting of the morphism $(\tilde{u}^- ,\tilde{u}^+)$
(these are the components of $u^-$ and $u^+$ analogous to the component
$\tilde{u}$ of $u$)
along the morphism
$$
C_{k/}(x_0,\ldots , x_k)\rightarrow
C_{(k-1)/}(x_0,\ldots , x_{k-1})\times _{C_{(k-1)/}(x_1,\ldots , x_k)}
C_{(k-2)/}(x_1,\ldots , x_{k-1}).
$$
\end{parag}

\begin{parag}
\label{diag3.6}
If $C$ is quasifibrant then the morphisms involved in the previous description
are fibrations. We obtain the following result: that if $E'_i\subset E_i$ are
trivial cofibrations and $C$ is quasifibrant then any diagram
$$
\Upsilon ^k(E'_1,\ldots , E'_k)\rightarrow C
$$
extends to a diagram
$$
\Upsilon ^k(E'_1,\ldots , E'_k)\rightarrow C.
$$
We can prove this by induction on $k$, and we are reduced exactly to the
lifting property for the trivial cofibration
$$
E'_1\times \ldots \times E'_k \hookrightarrow
E_1\times \ldots E_k
$$
along the morphism
$$
C_{k/}(x_0,\ldots , x_k)\rightarrow
C_{(k-1)/}(x_0,\ldots , x_{k-1})\times _{C_{(k-1)/}(x_1,\ldots , x_k)}
C_{(k-2)/}(x_1,\ldots , x_{k-1}).
$$
This morphism being fibrant by hypothesis, the lifting property holds.
\end{parag}

\begin{parag}
\label{diag3.7}
Suppose $C$ is quasifibrant. Then for $a,b$ fixed as in \ref{diag1.5} the
morphisms
$$
Diag ^{a,b}(E_1,\ldots , \underline{E_i},\ldots , E_k; C)
\rightarrow
Diag ^{a(i-1), b(i)}(\underline{E_i}; C)
$$
are fibrant weak equivalences, where $a(i-1)$ and $b(i)$ are the images
by $a$ and $b$ of the $i-1$-st and $i$-th vertices.

To prove this we use the description \ref{diag3.5}, inductively reducing
$k$. The
remark \ref{quasifib1.4} says that for any $k\geq 2$ the choice of lifting
$\tilde{u}$ doesn't change the equivalence type of the $Diag$ $n$-category.
This reduces down to the case $k=1$, which gives exactly that the
restriction to
the $i$-th edge is an equivalence (the restrictions to the other edges are
fixed and don't contribute anything because we fix $a,b$).
\end{parag}

\begin{parag}
\label{diag4}
{\em Proofs of \ref{diag2} and \ref{diag3}:}
The statements of \ref{diag2} are direct consequences of the lifting property
\ref{diag3.6}.

To prove \ref{diag3}, in view of \ref{diag3.7} it suffices to consider the case
$k=1$. Now
$$
Diag ^{U,V}(\underline{E}; C) = \underline{Hom}(E, C_{1/}(U,V)).
$$
Therefore if $C\rightarrow C'$ is any morphism of quasifibrant $n+1$-categories
which is ``fully faithful'' i.e. induces equivalences of fibrant $n$-categories
$$
C_{1/}(U,V)\rightarrow C'_{1/}(U,V),
$$
then
$$
Diag ^{U,V}(\underline{E}; C)\rightarrow
Diag ^{U,V}(\underline{E}; C')
$$
are equivalences by \ref{internal}.

(Note by \ref{quasifib1.5} that
the equivalence $C\rightarrow C'$ to a fibrant $C'$ that occurs in the
hypothesis
of \ref{diag3} is, in particular, a fully faithful morphism
of quasifibrant $n+1$-categories.)

This proves \ref{diag3} for $k=1$ and hence by \ref{diag3.7} for any $k$.
\end{parag}

\begin{parag}
\label{diag5}
The hypotheses on $C\rightarrow C'$ used in \ref{diag2} and \ref{diag3}
are satisfied by $nCAT \rightarrow nCAT'$, cf \ref{quasifib2}.
Therefore we may apply the results \ref{diag2} and \ref{diag3} to $nCAT
\rightarrow nCAT'$.

Fix $a,b$ as in \ref{diag1.5} for the following $Diag$'s, and suppose
that the restriction of $a$ to $\Upsilon (B)$ is equal to $1_B$.
From \ref{diag2} and \ref{diag3},
the morphism
$$
Diag ^{a,b}(B,E_2,\ldots , \underline{E_i}, \ldots , E_k; nCAT)
\rightarrow
Diag ^{a,b}(B,E_2,\ldots , \underline{E_i}, \ldots , E_k; nCAT')
$$
is an equivalence between fibrant $n$-categories.
\end{parag}

\bigskip

\subnumero{Some extensions}

\begin{parag}
\label{k23}
In the following preliminary statements we fix $k\geq 2$. We will only use
these statements for $k=2,3$.
\end{parag}

\begin{parag}
\label{def1B}
We now describe what will be our main technical tool.  Suppose $B$ is a fibrant
$n$-category.  We have a natural morphism
$$
1_B \in Hom ^{\ast , B}(\Upsilon B, nCAT')
$$
coming from the identity morphism $\ast \times B \rightarrow B$ in $nCAT$
(which is then considered as a morphism in $nCAT'$).
\end{parag}

\begin{parag}
\label{shl}
For any $E_1,E_2, \ldots , E_k$ let
$$
Shell\Upsilon ^k(E_1, \ldots , E_k):=
$$
$$
\Upsilon ^{k-1}(E_1,\ldots , E_{k-1})
\cup
\bigcup _{i=1}^{k-1}\Upsilon ^{k-1}(\ldots ,
E_i\times E_{i+1}, \ldots )
$$
(thus
it consists of all of the ``faces'' except  the one $\Upsilon
^{k-1}(E_2,\ldots , E_k)$).  We have a cofibration
$$
Shell\Upsilon ^k(E_1,\ldots , E_k)\rightarrow \Upsilon ^k(E_1,\ldots , E_k).
$$
\end{parag}

\begin{parag}
Now set $E_1=B$ and let $\underline{Hom}^{1_B}(\Upsilon ^k(B, E_2,\ldots , E_k),
nCAT')$ denote the fiber of
$$
\underline{Hom}(\Upsilon ^k(1_B, E_2,\ldots , E_k),nCAT')\rightarrow
\underline{Hom}(\Upsilon B, nCAT')
$$
over $1_B$ and let
$\underline{Hom}^{1_B}(Shell\Upsilon ^k(B, E_2,\ldots , E_k),
nCAT')$ denote the fiber of
$$
\underline{Hom}(Shell\Upsilon ^k(1_B, E_2,\ldots , E_k),nCAT')\rightarrow
\underline{Hom}(\Upsilon B, nCAT')
$$
over $1_B$.
\end{parag}

\begin{lemma}
\label{extension1}
Suppose $B$ is a fibrant $n$-category and $E$ an $n$-precat,
and $U$ an object of $nCAT'$ (it is also an object of $nCAT$).
The morphism
$$
Diag ^{1_B,U}(B,\underline{E}; nCAT')\rightarrow Diag ^{\ast
,U}(\underline{B\times
E}, nCAT')
$$
is an equivalence of $n$-categories.
\end{lemma}
{\em Proof:}
In view of \ref{diag2} and \ref{diag3} it
suffices to prove the same thing for diagrams in $nCAT$. In this case, use the
calculation of \ref{diag1.6}: both sides become equal to
$\underline{Hom}(B\times E, U)$.
\eop

\begin{remark}
\label{itsfibrant}
Since $nCAT'$ is a fibrant $n+1$-category, the morphism
$$
Diag ^{1_B,U}(B,\underline{E}; nCAT')\rightarrow Diag ^{\ast
,U}(\underline{B\times
E}, nCAT')
$$
is fibrant.  One checks directly the lifting property for a trivial cofibration
$F'\hookrightarrow F$, using the fibrant property of $nCAT'$.
\end{remark}

\begin{corollary}
\label{extension2}
The morphism
$$
Hom ^{1_B}(\Upsilon ^2(B,E),nCAT')\rightarrow Hom ^{\ast}(\Upsilon (B\times E),
nCAT')
$$
is surjective.
\end{corollary}
{\em Proof:}
We can fix an object $U$ for the image of the last vertex.
The morphism
$$
Diag ^{1_B,U}(B,\underline{E}; nCAT')\rightarrow Diag ^{\ast
,U}(\underline{B\times
E}, nCAT')
$$
is fibrant by the above remark \ref{itsfibrant}, and it is an equivalence by
\ref{extension1}. This implies that it is surjective on objects (\ref{surje}).
\eop

\begin{corollary}
\label{extension3}
Suppose $E'\subset E$ is a cofibration of $n$-precats. Suppose we are given
an object of
$$
Hom ^{1_B}(\Upsilon ^2(B, E'), nCAT').
$$
and an extension over the shell to an object of
$$
Hom ^{1_B}(Shell \Upsilon ^2(B,E), nCAT').
$$
Then these two have a common extension to an element of
$$
Hom ^{1_B}(\Upsilon ^2(B, E), nCAT').
$$
\end{corollary}
{\em Proof:}
Again we can fix $U$.
By Lemma \ref{extension1} and remark \ref{itsfibrant}, the morphism
$$
Diag ^{1_B,U}(B,\underline{\ast}; nCAT')\rightarrow Diag ^{\ast
,U}(\underline{B}, nCAT')
$$
is a trivial fibration.  Therefore it has the lifting property with respect to
any cofibration $E'\subset E$. This lifting property gives exactly what
we want to show---this is because a morphism
$$
E\rightarrow Diag ^{1_B,U}(B,\underline{\ast}; nCAT')
$$
is the same thing as an object of
$$
Diag ^{1_B,U}(B,E; nCAT')
$$
or equivalently an element of $Hom ^{1_B}(\Upsilon ^2(B, E), nCAT')$.
\eop

Now
we treat a similar type of extension problem for shells with $k=3$.

\begin{parag}
\label{extension4}
Now suppose we have an object $b\in Diag (F; nCAT')$.
Let
$$
Diag ^{1_B, b}_{\rm Shell}(B,\underline{E}, F; nCAT')
$$
be the $n$-precat representing the functor
$$
G\mapsto Hom ^{1_B, b} (Shell \Upsilon ^3(B,E\times G,F), nCAT')
$$
where the superscript on the $Hom$ has the obvious meaning that we look only at
morphisms restricting to $1_B$ on the edge $01$ and to $b$ on the edge $23$.

The shell $Shell \Upsilon ^3(B,E,F\times G)$ has three faces. We call the faces
$(013)$ and $(023)$ the {\em last faces} and the face $(012)$ the {\em first
face}.  Restriction to the last faces (which meet along the edge $(03)$) gives a
map
$$
Diag ^{1_B, b}_{\rm Shell}(B,\underline{E}, F; nCAT')
\rightarrow
$$
$$
Diag ^{1_B, b_3}(B, \underline{E\times F}; nCAT') \times _{Diag
^{\ast , b_3}
(\underline{B\times E \times F}; nCAT')} Diag ^{\ast , b}(\underline{B\times
E}, F; nCAT'),
$$
where $b_3$ denotes the object image of $3$ under the map $b$;
a similar
definition will hold for $b_2$ below---and recall that
the image of $0$ under the map $1_B$ is $\ast$.
\end{parag}

\begin{parag}
\label{extension5}
{\em Claim:} that the map at the end of the previous paragraph is a fibrant
equivalence.

Call the object on the right in this morphism $D$.  Restriction to the
edge $(02)$ is a map
$$
D\rightarrow Diag ^{\ast , b_2} (\underline{B\times E}; nCAT').
$$
We have an isomorphism
$$
Diag ^{1_B, b}_{\rm Shell}(B,\underline{E}, F; nCAT')
\stackrel{\cong}{\rightarrow} D \times
_{Diag ^{\ast , b_2} (\underline{B\times E}; nCAT')}
Diag ^{1_B, b_2}(B, \underline{E}; nCAT').
$$
However, the second morphism in this fiber product is
$$
Diag ^{1_B, b_2}(B, \underline{E}; nCAT')
\rightarrow
Diag ^{\ast , b_2} (\underline{B\times E}; nCAT')
$$
which is a fibrant equivalence by Lemma \ref{extension1}.
It follows that the morphism
$$
Diag ^{1_B, b}_{\rm Shell}(B,\underline{E}, F; nCAT')
\rightarrow D
$$
is a weak equivalence (note also that it is fibrant).
This proves the claim.
\end{parag}

\begin{corollary}
\label{extension6}
The morphism
$$
Diag ^{1_B, b}(B,\underline{E}, F; nCAT') \rightarrow
Diag ^{1_B, b}_{\rm Shell}(B,\underline{E}, F; nCAT')
$$
is a fibrant  equivalence.
\end{corollary}
{\em Proof:}
It is fibrant because $nCAT'$ is fibrant.
In view of the claim \ref{extension5} it suffices to note that
the map
$$
Diag ^{1_B, b}(B,\underline{E}, F; nCAT') \rightarrow
D
$$
is an equivalence, and this by the fibrant property of $nCAT'$
(the union of the faces $(013)$ and $(023)$ is one of the admissible ones in
our list of \ref{trivinclusions}).
\eop

\begin{corollary}
\label{extension7}
Suppose $E'\subset E$ is a cofibration of $n$-precats. Then
for any morphism
$$
\Upsilon ^3(B, E', F)\rightarrow nCAT'
$$
sending the edge $(01)$ to $1_B$, and any extension of this over the shell to a
morphism
$$
Shell \Upsilon ^3(B, E, F)\rightarrow nCAT'
$$
again restricting to $1_B$ on the edge $(01)$, there exists a common extension
to a morphism
$$
\Upsilon ^3(B, E, F)\rightarrow nCAT'.
$$
\end{corollary}
{\em Proof:}
By the previous Corollary \ref{extension6}  the morphism
$$
Diag ^{1_B, b}(B,\underline{\ast}, F; nCAT') \rightarrow
Diag ^{1_B, b}_{\rm Shell}(B,\underline{\ast}, F; nCAT'),
$$
is a fibrant equivalence. Therefore it satisfies the lifting property for any
cofibration $E' \hookrightarrow E$, and as before (\ref{extension3}) a map from
$E$ into
$Diag ^{1_B, b}(B,\underline{\ast}, F; nCAT')$ is the same thing as an object of
$Diag ^{1_B, b}(B,\underline{E}, F; nCAT')$ (and the same things for $E'$ and
for  $Diag ^{1_B, b}_{\rm Shell}$). This gives the required statement.
\eop

\bigskip

\subnumero{Extension properties for internal $\underline{Hom}$}

Now we take the above extension properties and recast them in terms of internal
$\underline{Hom}$. This is because we will need them for products of our
precats $\Upsilon$ with an arbitrary $A$. Note that there is a difference
between the internal $\underline{Hom}$ refered to in this section (which are
$n+1$-categories) and the $Diag$ $n$-categories above.

We state the following lemma for any value of $k$
 but we will only use $k=2$ and $k=3$; and we give the proofs only in these
cases, leaving it to the reader to fill in the combinatorial details for
arbitrary $k$.

\begin{lemma}
\label{main}
For any $n$-precats $E_2,\ldots , E_k$,
the morphism
$$
\underline{Hom}^{1_B}(\Upsilon ^k(B, E_2,\ldots , E_k),
nCAT') \rightarrow \underline{Hom}^{1_B}(Shell\Upsilon ^k(B, E_2,\ldots , E_k),
nCAT')
$$
is an equivalence of $n+1$-categories.
\end{lemma}
{\em Proof:}
The morphism in question is fibrant---cf \ref{internal}.

The proof is divided into several paragraphs. In \ref{main1}---\ref{main3}
we give the proof for $k=2$.
Then
in \ref{main4} we give the proof for $k=3$.

\begin{parag}
\label{main1}
We begin the proof for $k=2$.
Corollary \ref{extension2} implies that the morphism in question
$$
\underline{Hom}^{1_B}(\Upsilon ^2(B, E_2),
nCAT')\rightarrow
\underline{Hom}^{\ast}(\Upsilon (B\times E_2), nCAT')
$$
is surjective on objects.
\end{parag}

\begin{parag}
\label{main2}
Now we have to prove that our morphism induces equivalences between the
morphism $n$-categories. Suppose
$$
f,g: \Upsilon ^2(B, E_2)\rightarrow nCAT'
$$
are two morphisms (with the appropriate behavior on $(01)$).  Then the
$n$-category of morphisms between them represents the functor
$$
F\mapsto Hom ^{f,g; 1_B}(\Upsilon F \times \Upsilon ^2(B, E_2), nCAT')
$$
where the superscript means morphisms restricting to $f$ and $g$ over $0,1\in
\Upsilon F$ and restricting to $1_B$ over $\Upsilon F \times \Upsilon B$.
This maps (by restricting to the edge $02$) to the functor
$$
F\mapsto Hom ^{f,g; \ast}(\Upsilon F \times \Upsilon (B\times E_2), nCAT').
$$
We would like to prove that this restriction map of functors is an equivalence.
In order to prove this it suffices to prove that it has the lifting property
for any cofibrations $F'\subset F$.  Thus we suppose that we have a morphism
$$
\eta : \Upsilon F \times \Upsilon (B\times E_2)\rightarrow nCAT'
$$
(restricting appropriately to $f$ and $g$ and to $\ast$), as well as a morphism
$$
\zeta ' : \Upsilon F' \times \Upsilon ^2(B, E_2)\rightarrow nCAT',
$$
restricting appropriately to $f$, $g$ and $1_B$. We would like to extend this
latter to a map defined on $F$ and compatible with the previous one.
This extension will complete the proof for $k=2$.
\end{parag}

\begin{parag}
\label{main3}
We now prove the extension statement claimed above. As in \ref{pfvar4} we
consider the diagram as the product of an interval $(01)$ and a triangle $(012)$
and we denote the points by $(i,j)$ for $i=0,1$ and $j=0,1,2$.
More generally for example $(ab,cd)$ denotes the square which is the edge
$(ab)$ crossed with the edge $(cd)$. We are provided with maps on the end
triangles $f$ on $(0, 012)$ and $g$ on $(1,012)$ as well as $\eta $
on the top square
$(01, 02)$. We fix the map on the square $(01,01)$ (which is $\Upsilon (F)
\times \Upsilon (B)$ pullback of $1_B$, and call this again $1_B$. We are also
provided with a map $\zeta '$ defined on the whole diagram with respect to $F'$
and we would like to extend this all to $\zeta$ defined on the whole diagram.

Note that we can write $\Upsilon (F) \times \Upsilon ^2(B,E)$
as the coproduct of three tetrahedra which we denote
$$
(0,0)\;\; (1,0) \;\; (1,1) \;\; (1,2) \;\; i.e. \;\; \Upsilon ^3(F,B,E),
$$
$$
(0,0)\;\; (0,1) \;\; (0,2) \;\; (1,2) \;\; i.e. \;\; \Upsilon ^3(B,E, F),
$$
$$
(0,0)\;\; (0,1) \;\; (1,1) \;\; (1,2) \;\; i.e. \;\; \Upsilon ^3(B,F,E).
$$

The first step is to use the fibrant property of $nCAT'$ to extend
our given
morphisms $g$, the restriction of $1_B$ to the triangle $(0,0)$,
$(1,0)$, $(1,1)$, and the restriction of $\eta $ to the triangle $(0,0)$,
$(1,0)$, $(1,2)$, to a map on the tetrahedron
$$
(0,0)\;\; (1,0) \;\; (1,1) \;\; (1,2).
$$
We can do this in a way which extends the map $\zeta '$.

Next we again use the fibrant property of $nCAT'$ to extend across the
tetrahedron
$$
(0,0)\;\; (0,1) \;\; (0,2) \;\; (1,2).
$$
Note that we are provided with the map $f$ on the triangle
$(0,012)$, and the restriction of the $\eta$ on the triangle
$(0,0), (0,2), (1,2)$.
We can find our extension again in a way extending the given map $\zeta '$.

Finally we come to the tetrahedron
$$
(0,0)\;\; (0,1) \;\; (1,1) \;\; (1,2),
$$
which is of the form $\Upsilon ^3(B, F,E)$.
Here we are given maps on all of the faces except the last one,
i.e. on the shell of this tetrahedron, and we would like to extend
it. The given maps are the pullback of $1_B$ on the first face, and the maps
coming from the two previous paragraphs on the other two faces. Furthermore we
already have a map $\zeta '$ over the tetrahedron $\Upsilon ^3(B, F',E)$. The
given map on the shell restricts on the first edge to $1_B$, so this is an
extension problem of the type which we have already treated in Corollary
\ref{extension7} above. (N.B. the notations $E$ and $F$ are interchanged between
\ref{extension7} and the present situation.) Thus Corollary \ref{extension7}
provides the extension we are looking for, and we have finished making our
extension across the three tetrahedra.  This completes the proof of ``fully
faithfulness'' so the morphism in the lemma is an equivalence in the case $k=2$.
\end{parag}

\begin{parag}
\label{main4}
Here is the proof for $k=3$. First of all, the morphism
$$
\underline{Hom}^{1_B}(\Upsilon ^3(B, E_2,E_3),
nCAT')\rightarrow
$$
$$
\underline{Hom}^{1_B}(\Upsilon ^{2}(B, E_2),
nCAT') \times_{\underline{Hom}(\Upsilon  (B\times E_2),
nCAT')}
\underline{Hom}(\Upsilon ^2(B\times E_2, E_3),
nCAT')
$$
is an equivalence, by the fibrant property for $nCAT'$.

By the case $k=2$ (\ref{main1}--\ref{main3}) applied to the face $012$,
the morphism
$$
\underline{Hom}^{1_B}(Shell\Upsilon ^3(B, E_2,E_3),
nCAT')\rightarrow
$$
$$
\underline{Hom}^{1_B}(\Upsilon ^{2}(B, E_2),
nCAT') \times_{\underline{Hom}(\Upsilon  (B\times E_2),
nCAT')}
\underline{Hom}(\Upsilon ^2(B\times E_2, E_3),
nCAT')
$$
is an equivalence. This implies that the morphism
$$
\underline{Hom}^{1_B}(\Upsilon ^3(B, E_2,E_3),
nCAT')\rightarrow
\underline{Hom}^{1_B}(Shell\Upsilon ^3(B, E_2,E_3),
nCAT')
$$
is an equivalence.
This completes the case $k=3$.
\end{parag}

This completes the proof of the lemma (as far as we are going).
\eop

\begin{remark}
\label{main5}
One might think that we have a simple argument for the case $k=3$,
and the only difficult part of the argument for $k=2$ was the part where we
used $k=3$. However one cannot simplify the proof: the simple argument for
$k=3$ is based upon the use of internal $\underline{Hom}$ and
to get $k=2$ for internal $\underline{Hom}$ we need a statement like the case
of $k=3$---the statement which in the above proof is provided by Corollary
\ref{extension7}. This is why we were obliged to do all of the stuff in the
previous subsection.
\end{remark}

We will only use the subsequent corollaries in the cases $k=2$ and $k=3$, so the
proof we have given of \ref{main} is sufficient.  Again the reader is invited
to treat the case of any $k$.

\begin{corollary}
\label{main6}
Suppose $A$ is an $n+1$-precat and suppose $ E_i$
are $n$-precats for $i=2,\ldots , k$.  Suppose we are given
a morphism
$$
A\times Shell\Upsilon ^k(B,E_2,\ldots , E_k)\rightarrow nCAT'
$$
restricting to $1_B$ on $A\times \Upsilon B$.
Then there is an extension  to a morphism
$$
A\times \Upsilon ^k(B,E_2,\ldots , E_k)\rightarrow nCAT'.
$$
\end{corollary}
{\em Proof:}
The restriction morphism on the $\underline{Hom}$ is fibrant and
an equivalence by the previous lemma, therefore it is surjective on objects.
\eop

What we really need to know is a relative version of this for cofibrations
$E'_i\subset E_i$.

\begin{corollary}
\label{main7}
Suppose $A$ is an $n$-category and suppose $E'_i\subset E_i$
are cofibrations of $n$-precats for $i=2,\ldots , k$.  Suppose we are given
a morphism
$$
A\times Shell\Upsilon ^k(B,E_2,\ldots , E_k)\rightarrow nCAT'
$$
restricting to $1_B$ on $A\times \Upsilon B$,
together with a filling-in
$$
A\times \Upsilon ^k(B,E'_2,\ldots , E'_k)\rightarrow nCAT' ,
$$
then there is an extension of all of this to a morphism
$$
A\times \Upsilon ^k(B,E_2,\ldots , E_k)\rightarrow nCAT'.
$$
\end{corollary}
{\em Proof:}
Let $H_E$ denote the $\underline{Hom}$ for the full $\Upsilon$ and let
$H^{Sh}_E$ denote the $\underline{Hom}$ for $Shell\Upsilon$. The morphism
$$
H_E\rightarrow H^{Sh}_E \times _{H^{Sh}_{E'}}H_{E'}
$$
is an equivalence (as is seen by applying the lemma for both $E$ and $E'$)
and it
is fibrant (since it comes from $\underline{Hom}$ applied to a cofibration).
Therefore it is surjective on objects, which exactly means that we have the
above extension property.
\eop

\bigskip

\subnumero{Proof of Theorem \ref{inverse}}

\begin{parag}
\label{epsilon}
Recall that $\lambda$ was defined in \ref{lambda}.
We first apply the above statements to find our morphism
$\epsilon
:\lambda _A\rightarrow \varphi$.  The universal property of $\lambda$
(\ref{lambda}) applied to the identity map
$\lambda \rightarrow \lambda$ gives a
morphism
$$
\eta : A\times \Upsilon (\lambda )
\rightarrow nCAT'
$$
sending $A\times \{ 0\} $ to $\ast _A$ and sending
$A\times \{ 1\}$ to $\varphi$.  By Corollary \ref{main6} (for $k=2$) there is a
morphism
$$
\epsilon ^{(2)} :A\times \Upsilon ^2(\lambda , \ast )
\rightarrow nCAT'
$$
such that
$$
r_{02}(\epsilon ^{(2)}) = \eta
$$
and
$$
r_{01}(\epsilon ^{(2)})=1_{\lambda} .
$$
Note that $r_{12}(\epsilon ^{(2)})$ is a morphism from $A\times \Upsilon \ast
= A\times I$ into $nCAT'$ restricting to $\lambda _A$ and $\varphi$, which by
 definition means a morphism $\lambda \rightarrow \varphi$. Call this morphism
$\epsilon$.
\end{parag}

\begin{parag}
\label{proofclaim}
{\em Claim:} that for any fibrant $n$-category $B$ and any morphism
$$
f:A\times \Upsilon E \rightarrow nCAT'
$$
with
$$
r_0(f)= B_A,\;\;\; r_1(f)=\varphi
$$
there is a morphism
$$
f':\Upsilon ^2(E,\ast ) \rightarrow nCAT'
$$
with
$$
r_{02}(f')= f,\;\;\; r_{12}(f')=\epsilon .
$$
This almost gives the required property to show that $\lambda
\stackrel{\epsilon}{\rightarrow }\varphi$ is an inverse limit. Technically
speaking we also will have to show the above claim in the relative situation of
$E'\subset E$. This we will do below (\ref{proofF}--\ref{proofEnd}) after first
going through the argument in the absolute case (\ref{proofA}--\ref{proofE}).
\end{parag}

\begin{parag}
\label{proofA}
The basic idea is to use what we know up until now to construct a morphism
$$
F: A\times \Upsilon ^3(B, E,\ast )\rightarrow nCAT'
$$
with
$$
r_{01}(F)= 1_B,\;\;\; r_{13}(F)= f, \;\;\; r_{23}(F) = \epsilon .
$$
Setting $f'= r_{123}(F)$ we will obtain the morphism asked for in the
previous paragraph. In order to follow the construction the reader is urged to
draw a tetrahedron with vertices labeled $0,1,2,3$, putting respectively
$B$, $E$, $\ast$, $B\times E$, $E$, $B\times E$ along the edges $01$, $12$,
$23$, $02$, $13$, $03$; then putting in $\ast _A$, $B_A$, $\lambda _A$ and
$\varphi$ at the vertices $0,1,2,3$ respectively.  And finally putting in
$1_B$ along edge $01$, $f$ along edge $13$ and $\epsilon$ along edge $23$.

Our strategy is to fill in all of the faces except $123$, then call upon
Corollary \ref{main6}  to fill in the tetrahedron thus getting face $123$.
\end{parag}

\begin{parag}
\label{proofB}
The first step is the face $013$. This we fill in using simply the fact that
$nCAT'$ is a fibrant $n+1$-category. The edges $01$ and $13$ are specified so
we can fill in  to a morphism $A\times \Upsilon ^2(B,E)\rightarrow nCAT'$
(restricting to $1_B$ and $f$ on the edges $01$ and $13$).  Now the restriction
of this morphism to edge $03$ provides  a morphism $g:A\times \Upsilon (B\times
E)\rightarrow nCAT'$ restricting to $\ast _A$ and $\varphi$.
\end{parag}

\begin{parag}
\label{proofC}
The next step is to notice that by the universal property
\ref{lambda} of $\lambda$ there is
a morphism $B\times E\rightarrow \lambda$ such that $g$ is deduced from $\eta$
by pullback via $\Upsilon (E\times B)\rightarrow \Upsilon (\lambda )$.
This same morphism yields
$$
\Upsilon ^2(E\times B, \ast )\rightarrow \Upsilon ^2(\lambda , \ast ),
$$
and we can use this to pull back the above morphism $\epsilon ^{(2)}$.
This gives a morphism
$$
h: A \times \Upsilon ^2(E\times B, \ast ) \rightarrow nCAT'
$$
where (adopting exceptionally for obvious reasons here the notations
$0$, $2$ and $3$ for the vertices of this $\Upsilon ^2$)
$$
r_{03}(h)= g, \;\;\; r_{23}(h)=\epsilon .
$$
This treats the face $023$.
\end{parag}

\begin{parag}
\label{proofD}
Finally, for the face $012$ we have a morphism
$$
r_{02}(h):A\times \Upsilon (E\times B)\rightarrow nCAT'
$$
restricting to $\ast _A$ and $\lambda _A$.  By Corollary \ref{main6}
applied with
$k=2$ (for the map
$$
A\times Shell\Upsilon ^2(B, E)\rightarrow nCAT'
$$
given by $1_B$ and $h$) we get a morphism
$$
m: \Upsilon ^2(B,E)\rightarrow nCAT'
$$
with $r_{01}(m)= 1_B$ and $r_{02}(m)= r_{02}(h)$.
\end{parag}

\begin{parag}
\label{proofE}
Putting all of these together
we obtain a morphism
$$
F' : A\times Shell \Upsilon ^3(B,E, \ast )\rightarrow nCAT'
$$
restricting to $1_B$ on edge $01$ and restricting to $f$ on edge $13$ and
$\epsilon$ on edge $23$.  Corollary \ref{main6} applied with $k=3$ gives an
extension over the tetrahedron to a morphism
$$
F : A\times \Upsilon ^3(B,E, \ast )\rightarrow nCAT'
$$
again restricting to $1_B$ on edge $01$ and restricting to $f$ on edge $13$ and
$\epsilon$ on edge $23$. The restriction to the last face $r_{123}$ yields
the filling-in desired.
\end{parag}

\begin{parag}
\label{proofF}
We now treat the case where $E'\subset E$ is a cofibration and where we already
have a filling-in of the face $123$ for $E'$.  We would like to obtain a
filling-in of this face for $E$.  Basically the only difficulty is that
we don't  yet know that the filling-in of face $123$ for $E'$ comes from
a filling-in of the whole tetrahedron compatible with the above process.
In particular this causes a problem at the step where we fill in face
$023$.
\end{parag}

\begin{parag}
\label{proofG}
Before getting started we use the fibrant property of $nCAT'$ to obtain a
morphism
$$
A\times \Upsilon ^3(B,E', \ast )\rightarrow nCAT'
$$
restricting to our given morphism on the face $123$, and restricting to $1_B$
on the edge $01$. Actually we would like to insure that the restriction to
the face $012$ comes from a morphism
$$
\Upsilon ^2(B,E')\rightarrow nCAT'
$$
by pulling back along the projection $A\rightarrow \ast$.  In order to do this
notice that the restriction of the given map to the edge $12$ comes
from $\Upsilon (E')\rightarrow nCAT'$.  Thus we can first extend this map
combined with $1_B$ to a morphism $\Upsilon ^2(B,E')\rightarrow nCAT'$.
Now the morphism
$$
(A\times \Upsilon ^2(B,E'))\cup ^{A\times \Upsilon E'}
A\times \Upsilon ^2(E', \ast )\rightarrow \Upsilon ^3(B,E', \ast )
$$
is a trivial cofibration so we can extend from here to obtain
$A\times \Upsilon ^3(B,E', \ast )\rightarrow nCAT'$ with restriction to the
face $012$ coming from $\Upsilon ^2(B,E')\rightarrow nCAT'$. This is our point
of departure for the rest of the argument.
\end{parag}

\begin{parag}
\label{proofH}
The first step following the previous outline is to fill in the face $013$.  We
note that the morphism  $\Upsilon B \cup ^{\{ 1\} } \Upsilon E \rightarrow
\Upsilon
^2(B,E)$ is a trivial cofibration. Thus also the morphism
$$
(\Upsilon B \cup ^{\{ 1\} } \Upsilon E)\cup ^{\Upsilon B \cup ^{\{ 1\} }
\Upsilon E'}
\Upsilon ^2(B,E') \rightarrow \Upsilon ^2(B,E)
$$
is a trivial cofibration, so given the edges $01$ and $13$ (for $E$) with
filling-in over the face $013$ with respect  to $E'$, we can fill in $013$
with respect to $E$.
\end{parag}

\subsubnumero{The face $023$}

Now we treat the face $023$.
Let
$$
g: A\times \Upsilon (E\times B)\rightarrow nCAT'
$$
be the restriction of the map obtained in \ref{proofH} to the edge $03$.  Let
$g'$ denote its restriction to $A\times \Upsilon (E'\times B)$.
The map  given in \ref{proofG} restricts on $(023)$ to a morphism
$$
h': A \times \Upsilon ^2(E'\times B, \ast ) \rightarrow nCAT'
$$
where (using as above the notations
$0$, $2$ and $3$ for the vertices of this $\Upsilon ^2$)
$$
r_{03}(h')= g', \;\;\; r_{23}(h')=\epsilon .
$$
Let $a'=r_{02}(h')$. It is a morphism
$$
a'_A:A\times \Upsilon (E'\times B)\rightarrow nCAT'
$$
with $r_0(a'_A)= \ast _A$ and $r_2(a'_A)= \lambda _A$.  By hypothesis on our
map over the full tetrahedron for $E'$ (cf \ref{proofF}), $a'_A$ comes from a
map  $$
a':\Upsilon
(E'\times B)\rightarrow nCAT'
$$
again with values $\ast$ and $\lambda$ on the
endpoints.  This map corresponds to
$$
E'\times B\rightarrow nCAT'_{1/}(\ast ,
\lambda ).
$$
The morphism $nCAT\rightarrow nCAT'$ is an equivalence so our
morphism is equivalent to a different morphism $b':E'\times B\rightarrow
nCAT_{1/}(\ast , \lambda )$.   These two resulting morphisms $\Upsilon (E'\times
B)\rightarrow nCAT'$ are equivalent
so  by \ref{homotopic1}, \ref{homotopic2} there is a
morphism
$$
\overline{I}\times \Upsilon (E'\times B)\rightarrow nCAT'
$$
sending the
endpoints $0,1\in \overline{I}$ to $a'$ and $b'$.

Using this different morphism $b'$ (which is now the same thing as a map
$E'\times B\rightarrow \lambda )$ we pull back our standard
$$
\eta \in \underline{Hom}(A\times \Upsilon ^2(\lambda , \ast ), nCAT')
$$
to get a morphism
$$
A\times \Upsilon ^2(E'\times B, \ast )\rightarrow nCAT'
$$
restricting on the edges to $b'$ and $\epsilon$ respectively.

Now we have a map from
$$
\left( A\times \Upsilon ^2(E'\times B, \ast ) \right) \cup
\left( A\times \overline{I}\times [\Upsilon (E'\times B)\cup \Upsilon
(\ast )]\right)  \cup
\left( A \times \Upsilon ^2(E'\times B, \ast )\right)
$$
to $nCAT'$,
where  the first term is glued to the second term along $1\in \overline{I}$
and the last term is glued to the second term along $0\in \overline{I}$
(we have omitted in the notation the $n+1$-precats along which the glueing takes
place, the reader may fill them in as an exercise!).

The morphism from the above domain to
$$
A\times \overline{I} \times \Upsilon ^2(E'\times B, \ast )
$$
is a trivial cofibration, so since $nCAT'$ is fibrant there exists an extension
of the above to a morphism
$$
A\times \overline{I} \times \Upsilon ^2(E'\times B, \ast )\rightarrow nCAT'.
$$
This morphism is a standard one coming from $b': E'\times B \rightarrow \lambda$
on the end $1\in \overline{I}$, and it is our given $h'$ on the
end $0\in \overline{I}$.

We now go to the edge $03$ of the triangle $023$. We are also given an
extension
of $g'$ to $g: A\times \Upsilon (E\times B)\rightarrow nCAT'$ along
the edge $03$ of the triangle and $0$ of the interval $\overline{I}$.
Thus, using the face $(03)\times \overline{I}$,
we have a morphism
$$
(A\times \Upsilon (E\times B))\cup ^{A\times \Upsilon (E'\times B)}
(A\times \overline{I}\times \Upsilon (E'\times B))
\rightarrow nCAT'.
$$
Fill this in along the trivial cofibration
$$
(A\times \Upsilon (E\times B))\cup ^{A\times \Upsilon (E'\times B)}
(A\times \overline{I}\times \Upsilon (E'\times B))
\hookrightarrow
A\times \overline{I}\times \Upsilon (E\times B),
$$
to give on the whole a
morphism
$$
A\times \overline{I}\times \Upsilon (E\times B) \cup ^{
A\times \overline{I}\times \Upsilon (E'\times B)}
(A\times \overline{I}\times \Upsilon ^2(E'\times B , \ast )
\rightarrow nCAT',
$$
where the morphism $\Upsilon (E'\times B)\rightarrow \Upsilon ^2(E'\times
B,\ast )$ in question is the one coming from the edge $03$.

Next we extend down along the triangle $023$ times the end $1\in \overline{I}$.
To do this, notice that our extension from the previous paragraph gives
an extension of the morphism $b': \Upsilon (E'\times B)\rightarrow nCAT'$
to a morphism $b: \Upsilon (E\times B)\rightarrow nCAT'$. By the universal
property of $\lambda$ this corresponds to an extension  $E\times B\rightarrow
\lambda$.  Now the morphism that we already have on the end $1\in
\overline{I}$ comes by pulling back the standard $\eta : A\times \Upsilon
^2(\lambda , \ast )\rightarrow nCAT'$ via the map $E'\times B\rightarrow
\lambda$ so our extension allows us to pull back $\eta$ to get a map
$b: A\times \Upsilon ^2(E\times B, \ast )\rightarrow nCAT'$ extending the
previous $b'$.

Now we have our map
$$
A\times \overline{I} \times \Upsilon ^2(E'\times B , \ast )\rightarrow nCAT'
$$
which is provided with an extension from $E'$ to $E$, over the faces
$(03)\times \overline{I}$ and $023 \times 1$ of the product of the triangle
with the interval.
Another small step is to notice that along the face $(02)\times \overline{I}$
the morphism is pulled back along $A\rightarrow \ast$ from a morphism
$\overline{I} \times \Upsilon (E'\times B )\rightarrow nCAT'$. On the other
hand at the edge $(02)\times \{ 1\}$ the extension from $E'$ to $E$ again
comes from a morphism $\Upsilon (E\times B)\rightarrow nCAT'$.  We get a
morphism
$$
\overline{I} \times \Upsilon (E'\times B) \times ^{\{ 1\} \times \Upsilon
(E'\times B)} \Upsilon (E\times B) \rightarrow nCAT',
$$
which can be extended along the trivial cofibration
$$
\overline{I} \times \Upsilon (E'\times B) \times ^{\{ 1\} \times \Upsilon
(E'\times B)} \Upsilon (E\times B) \hookrightarrow
\overline{I}\times \Upsilon (E\times B)
$$
to give a map
$$
\overline{I}\times \Upsilon (E\times B) \rightarrow nCAT'.
$$
Similarly we note that the map on the face $(23) \times \overline{I}$
is pulled back from our map $\epsilon : \Upsilon (\ast )\rightarrow nCAT'$.

All together on the triangular icosahedron $(023) \times \overline{I}$
we have a morphism defined for $E'$ plus, along the faces
$$
(03)\times \overline{I}, \;\; (02)\times \overline{I}, \;\; (23)\times
\overline{I}, \;\; (023)\times \{ 1\}
$$
extensions from $E'$ to $E$ (all compatible on intersections of the faces and
having the required properties along $02$ and $23$).
The inclusion of this
$n+1$-precat (which we will call ${\bf G}$ for ``gory'' instead of writing it
out) into
$$
A\times \overline{I} \times \Upsilon ^2(E\times B , \ast )
$$
is a trivial cofibration.  Indeed ${\bf G}$ comes by attaching to the
end
$$
A \times \{
1\} \times \Upsilon ^2(E\times B , \ast ),
$$
something of the form
$$
A\times \partial \Upsilon ^2(E\times B, \ast )
\cup ^{A\times \partial \Upsilon ^2(E'\times B, \ast )}
A\times \Upsilon ^2(E'\times B,
\ast )
$$
where $\partial \Upsilon ^2(E\times B, \ast )$ denotes
the coproduct of the three ``edges'' $\Upsilon (E\times B)$ (two times)
and $\Upsilon (\ast )$.
The
inclusion of the end $A \times \{ 1\} \times \Upsilon ^2(E\times B , \ast )$
into ${\bf G}$ is an equivalence, as is the inclusion of this end into
the full product
$$ A\times \overline{I} \times \Upsilon ^2(E\times B , \ast ),
$$
which proves that the map in question (from ${\bf G}$ to
the above full product) is a weak
equivalence (and it is obviously a cofibration).

Now we again make use of  the fibrant property to extend our map
from ${\bf G}$ to a morphism
$$
A\times \overline{I} \times \Upsilon ^2(E\times B , \ast )\rightarrow
nCAT'.
$$
When restricted to $A\times \{ 0\} \times \Upsilon ^2(E\times B ,
\ast )$ this gives the extension $h$ desired in order to complete our treatment
of the face $023$.

\begin{parag}
For the face $012$ the argument is the same as in the previous case but we
apply Corollary \ref{main7} rather than \ref{main6} in view of our relative
situation $E'\subset E$.
\end{parag}

\begin{parag}
\label{proofEnd}
{\em End of the proof of \ref{inverse}}

We have constructed a morphism
$$
F' : A\times Shell \Upsilon ^3(B,E, \ast )\rightarrow nCAT'
$$
restricting to $1_B$ on edge $01$ and restricting to $f$ on edge $13$ and
$\epsilon$ on edge $23$.  Furthermore, by construction it restricts to our
already-given morphism over $E'$.
Corollary \ref{main7} applied with $k=3$ gives an extension over
the tetrahedron to a morphism
$$
F : A\times \Upsilon ^3(B,E, \ast )\rightarrow nCAT'
$$
again restricting to $1_B$ on edge $01$ and restricting to $f$ on edge $13$ and
$\epsilon$ on edge $23$, and restricting to the already-given morphism over
$E'$. The restriction to the last face $r_{123}(F)$ yields the filling-in
desired. This completes the proof that
$\lambda
\stackrel{\epsilon}{\rightarrow }\varphi$ is an inverse limit,
finishing
the proof of Theorem \ref{inverse}.
\eop
\end{parag}

\begin{corollary}
\label{commutencat}
If $F:A\times B\rightarrow nCAT'$ is a functor from the product of two
$n+1$-categories, then taking the inverse limits first in one direction and
then in the other, is independent of which direction is chosen first.
\end{corollary}
{\em Proof:}
This is a consequence of Theorem \ref{commute} but can also be seen directly
from the construction \ref{lambda} of the limit.
\eop

\numero{Direct limits}

\begin{theorem}
\label{direct}
The $n+1$-category $nCAT'$ admits direct limits.
\end{theorem}

One should probably be able to construct these direct limits
in much the same way as in the topological case, roughly speaking by
replacing a family by an equivalent one in which the morphisms are cofibrations
(some type of telescope construction) and then taking the direct
limit of $n$-precats in the usual sense.  This seems a bit complicated to put
into practice so we will avoid doing so by a trick.

\begin{parag}
\label{arg1cat}
{\em The argument for a $1$-category:}
Consider the following argument which shows that if $C$ is a category in which
all inverse limits exist and in which projectors are effective, then $C$ admits
direct limits. For a functor $\psi :A\rightarrow C$ let $D$ be the category
whose
objects are pairs $(c,u)$ where $c$ is an object of $C$ and $u: \psi \rightarrow
c$ is a morphism. There is a forgetful functor $f:D\rightarrow C$. Let $\delta
\in C$ be the inverse limit of $f$.  Then for any $a\in A$ there is a unique
morphism  $\psi (a)\rightarrow f$. By the inverse limit property this yields
a morphism $\psi (a)\rightarrow \delta$ and uniqueness implies that it is
functorial in $a$. Thus we get a morphism $v: \psi \rightarrow \delta$ and
$(\delta , v)$ is in $D$. As an object of $D$, $\delta$ has a morphism
$$
p:=v(\delta ): \delta \rightarrow \delta.
$$
This is itself a morphism in $D$, so we get $p\circ p  = p$ from naturality
of $v$. Thus $p$ is a projector. Let $t$ be the direct factor of $\delta$
given by $p$. Composition $\psi \rightarrow \delta \rightarrow t$ gives
a map $\psi \rightarrow t$ and
we get a factorization $\psi \rightarrow t \rightarrow \delta$. Now
$t$ is seen to be an initial object of $D$, hence $\psi \rightarrow t$ is
a direct limit.
\end{parag}

\begin{parag}
\label{problem}
The only problem with this argument is a set-theoretic one.
Namely, when one speaks of
``limits'' it is presupposed that the indexing category $A$ is small, i.e. is a
set of some cardinality rather than a class. However our category $C$ is likely
to be a class. Thus, in the above argument, $D$ is not small and we are not
allowed to take the inverse limit over $D$.
\end{parag}

\begin{parag}
\label{fixup}
Let's see how to fix this up in the case $C=Set$ is the category of sets.
Suppose we have a functor $\psi :
A\rightarrow Set$ from a small category $A$. Let $\alpha$ be a cardinal number
bigger than $|A|$ and bigger than the cardinal of any set in the image of
$\psi$.  Let $D_{\alpha}$ be the category of pairs $(c,u)$ as above where $c$ is
contained in a fixed set of cardinality $\alpha$.  Note that $D_{\alpha}$ has
cardinality $\leq 2^{\alpha}$.  Let $(\delta , v)$ be as above.  The only hitch
is that (since we know an expression of $\delta$ as a subset of certain types of
functions on $D_{\alpha}$ with values in the parametrizing sets which themselves
have cardinality $\leq \alpha$) the cardinality of $\delta$ seems {\em a priori}
only to be bounded by $2^{2^{\alpha}}$.  Let $\delta ' \subset \delta $ be the
smallest subset through which the map $v: \psi \rightarrow \delta$ factors.
Note that the cardinality of $\delta '$ cannot be bigger than the sum of the
cardinals of the $\psi (a)$ over $a\in A$, in particular $\delta '$ has
cardinal $\leq \alpha$.  But now the universal property of $\delta$ implies
that $\delta = \delta '$, for it is easy to see that $\delta '\rightarrow f$ is
again an inverse limit.  Thus by actually counting we see that the
cardinality of $\delta$ is really $\leq \alpha$ and up to isomorphism we may
assume that $(\delta , v)\in D_{\alpha}$.  This argument actually shows that
the cardinality of $\delta$ is bounded independantly of the choice of $\alpha$.
Thus $\delta$ satisfies the universal property of a direct limit for morphisms
to a set of any cardinality, so $\delta$ is the direct limit of $A$.

More generally, in the situation of \ref{arg1cat} if we can define the
$D_{\alpha}$ and if we know for some reason that every object $B\in D$ admits a
map $B'\rightarrow D$ from an object $D'\in D_{\alpha}$ then we can fix up the
argument.
\end{parag}

We would like to do the same thing for limits in $nCAT'$, namely show that
direct limits exist just using a general argument working from the existence of
inverse limits. In order to do this we first need to discuss cardinality
questions for $n$-categories.

\bigskip

\subnumero{Cardinality}

Suppose $A$ is an $n$-category. We define the {\em cardinal of $A$}, denoted
$\# A$
in the following way. Choose for every $y\in \pi _0(A)=\tau _{\leq 0}(A)$
(the set of equivalence classes of objects) a
lifting to an object $\tilde{y}\in A_0$. Then
$$
\# A := \sum _{y,z\in \pi _0(A)}\# A_{1/}(\tilde{y},\tilde{z}).
$$
The sum of cardinals is of course the cardinal of the disjoint union of
representing sets. This definition is recursive, as what goes into the
formula is the cardinal of the $n-1$-category
$A_{1/}(\tilde{y},\tilde{z})$.  At the start we define the cardinal of a
$0$-category (i.e. a set) in the usual way.

\begin{lemma}
The above definition of $\# A$ doesn't depend on choice of representatives. If
$A\rightarrow B$ is an equivalence of $n$-categories then $\# A = \# B$.
\end{lemma}
{\em Proof:}
Left to the reader.
\eop

An easier and more obvious notion is the {\em precardinality of $A$}. If $A$ is
any $n$-precat we define (with the notations of \cite{nCAT})
$$
\# ^{\rm pre}A:= \sum _{M\in \Theta ^n} \# (A_M).
$$
For infinite cardinalities the precardinal of $A$ is also the maximum of the
cardinalities of the sets $A_M$. In any case note that the precardinality
is infinite unless $A$ is empty.

\begin{remark}
 Let $A\mapsto Cat(A)$ denote the operation of replacing an
$n$-precat by the associated $n$-category. Then the precardinal of $Cat(A)$ is
bounded by the maximum of $\omega$ and
the precardinal of $A$. Similarly by the argument of (\cite{nCAT} \S 6, proof
of CM5(1)), for any $n$-precat $A$ there is a replacement by a fibrant
$n$-category $A\hookrightarrow A'$ with
$$
\# ^{\rm pre}A ' \leq  max (\omega , \# ^{\rm pre}A).
$$
Actually since $ \# ^{\rm pre}A\geq \omega$ we can write more simply that
$\# ^{\rm pre}A ' = \# ^{\rm pre}A$.
\end{remark}

Note trivially that
$$
\# A \leq \# ^{\rm pre}A .
$$
The following lemma gives a converse up to equivalence.

\begin{lemma}
Suppose $A$ is an $n$-category with $\# A\leq  \alpha$ for an infinite cardinal
$\alpha$. Then $A$ is equivalent to an $n$-category $A'$ of precardinality
$\leq \alpha$.
\end{lemma}
{\em Proof:}
Left to the reader.
\eop

\bigskip

\subnumero{A criterion for direct limits in $nCAT'$}

Before getting to the application of the theory of cardinality
we give a criterion which simplifies the problem of finding direct limits
in $nCAT'$.

\begin{parag}
\label{criterion1}
For this section we need another type of universal morphism. Suppose
$E$ and $B$ are $n$-precats, with $B$ fibrant. Then $\underline{Hom}(B,E)$
is fibrant and we have a canonical morphism
$$
\underline{Hom}(B,E)\times E \rightarrow B.
$$
This may be interpreted as an object
$$
\nu \in Diag ^{\underline{Hom}(B,E), B}(E; nCAT)
$$
which yields by composition with $nCAT \rightarrow nCAT'$
the element which we denote by the same symbol
$$
\nu \in Diag ^{\underline{Hom}(B,E), B}(E; nCAT')
$$
\end{parag}

\begin{parag}
\label{criterion2}
The element $\nu$ has the following universal property: for any $n$-precat
$F$ the  morphism
$$
Diag ^{U , \nu }(\underline{F}, E; nCAT')
\rightarrow
Diag ^{U , B }(\underline{F\times E}; nCAT')
$$
is a fibrant equivalence of fibrant $n$-categories.

To prove this note that the fibrant property comes from the fact that $nCAT'$
is fibrant.
Note that both sides are fibrant by \ref{diag2}. The fact that it is an
equivalence may be checked using diagrams in $nCAT$ rather than diagrams in
$nCAT'$, according to \ref{diag3}. Using \ref{diag1.6} for diagrams in $nCAT$,
both sides are equal to
$$
\underline{Hom}(U\times F\times E, B)
$$
where $U$ is the image of the first object $0\in \Upsilon ^2(F,E)$.
This shows that the morphism is an equivalence.
\end{parag}

\begin{parag}
\label{criterion3}
As a corollary of the above, given a morphism
$$
f:\Upsilon(F\times E) \rightarrow nCAT'
$$
with image of the last vertex equal to $B$, there is an extension
to a morphism
$$
g:\Upsilon ^2(F , E)\rightarrow nCAT'
$$
such that $r_{02}(g)=f$ and $r_{12}(f) = \nu$.
Similarly if $E'\subset E$ and we are already given the extension $g'$
for $E'$ then we can assume that $g$ is compatible with $g'$.
\end{parag}

\begin{parag}
\label{criterion3.5}
We also have a version of this universal property for shell-extension
in higher degree.  This concerns the {\em right shell} $Shelr \Upsilon ^k$
(cf \ref{shl}).  Suppose we are given a morphism
$$
f: Shelr \Upsilon ^k(F_1,\ldots , F_{k-1}, E)\rightarrow nCAT'
$$
such that $f$ restricts on the last edge to $\nu$.  Then there is a
filling-in to a morphism
$$
g: \Upsilon ^k(F_1,\ldots , F_{k-1}, E)\rightarrow nCAT'.
$$
If $g'$ is already given over $F'_1,\ldots , F'_{k-1}, E'$
then we can assume that $g$ is compatible with $g'$.
This is the analogue of \ref{extension7} and the like.
\end{parag}

\begin{parag}
\label{criterion4}
The above property also works in a family.
Given a morphism
$$
f:A\times \Upsilon(F\times E) \rightarrow nCAT'
$$
sending the last vertex to the constant object $B_A$, there is an extension
to a morphism
$$
g:A\times \Upsilon ^2(F , E)\rightarrow nCAT'
$$
such that $r_{02}(g)=f$ and $r_{12}(f) = \nu _A$ is the morphism pulled
back from
$\nu$. Again if an extension $g'$ is already given on $E'\subset E$ then $g$ may
be chosen compatibly with $g'$.

Similarly there is a shell-extension property as in \ref{criterion3.5}
in a family.

For the proof
one has to go through a procedure analogous to the passage from diagrams
to internal $\underline{Hom}$ in \ref{main1}--\ref{main7}. This discussion of
the universal morphism $\nu$ is parallel to the discussion of the discussion of
the universal $1_B$, but with ``arrows reversed''.
\end{parag}

We now come to our simplified criterion for limits in $nCAT'$.

\begin{parag}
\label{caution}
{\em Caution:} Note that the following  lemma only applies as such to limits
taken in $nCAT'$ and not in general to limits in an arbitrary $n+1$-category
$\Cc$. The proof uses in an essential way the fact that the morphism
objects for the ``category''  $nCAT'$ are $n$-categories which are also
basically the same thing as the objects of $nCAT'$. Of course it is possible
that the same techniques of proof might work in a limited other range of
circumstances which are closely related to these.
\end{parag}

\begin{lemma}
\label{easydirect}
Suppose $A$ is an $n+1$-precat and $\psi : A \rightarrow
nCAT'$ is a morphism. Suppose that $\epsilon : \psi \rightarrow \delta$ is a
morphism to an object $\delta \in nCAT'$ having the following weak limit-like
property: for any other morphism $f:\psi \rightarrow \mu$ there exists a
morphism $g:\delta \rightarrow \mu$  such that the composition $g\epsilon $
(well defined up to homotopy) is homotopic to $f$; and furthermore that such a
factorization is unique up to a (not necessarily unique) homotopy of the
factorization. Then $\psi
\stackrel{\epsilon}{\rightarrow} \delta $ is a direct limit.
\end{lemma}
{\em Proof:}
First we explain more precisely what the existence and uniqueness of the
factorization mean.  Given an element $f\in Hom (\psi , \mu )$ there exists an
element  $g'\in Hom ^{\epsilon }
(\psi , \delta , \mu )$ projecting via $r_{02}$ to a
morphism equivalent to $f$. This equivalence may be measured in the $n$-category
$Hom (\psi , \mu )$. Note that since $nCAT'$ is fibrant the projection
$$
Hom ^{\epsilon }
(\psi , \delta , \mu ) \rightarrow Hom (\psi , \mu )
$$
is fibrant, so if an object equivalent to $f$ is in the image then $f$ is in
the image. Thus we can restate the criterion as saying simply that there exists
an element $g'$ projecting via $r_{02}$ to $f$.

Suppose given two such factorizations $g'_1$ and $g'_2$.
By ``homotopy of the factorization'' we mean a homotopy between $r_{12}(g'_1)$
and $r_{12}(g'_2)$ such that the resulting homotopy between $f$ and itself
(this homotopy being well defined up to $2$-homotopy) is $2$-homotopic to the
identity $1_f$. Again using the fibrant condition of $nCAT'$ we obtain that
this condition implies the simpler statement that there exists
a morphism
$$
A\times \Upsilon ^2(\ast , \ast ) \times \overline{I}\rightarrow nCAT'
$$
restricting to $g'_1$ and $g'_2$ on the two endpoints $0,1\in \overline{I}$;
restricting to the pullback $\epsilon$ on the edge $(01)$ of the $\Upsilon ^2$,
this edge being $A\times \Upsilon (\ast )\times \overline{I}$; and restricting
to the pullback of $f$ on the edge $(02)$ which is $A\times \Upsilon
(\ast ) \times \overline{I}$.

\begin{parag}
\label{existen}
{\em Simple factorization}
We start by showing the simple version of the factorization property
necessary to show that $\epsilon$ is an inverse limit; we will treat the
relative case for $E'\subset E$ below.So for now, suppose that we are given a
morphism
$$
u:A \times \Upsilon (E) \rightarrow nCAT'
$$
restricting to $\psi$ on $A\times \{ 0\} $ and restricting to a constant object
$B\in nCAT'$ (i.e. to the pullback $B_A$) on $A\times \{ 1\}$. We would like
to extend this to a morphism
$$
v: A\times \Upsilon ^2(\ast , E) \rightarrow nCAT'
$$
restricting to our given morphism on the edge $(02)$, and restricting to
$\epsilon$ on the edge $(01)$.  Our given morphism corresponds by
\ref{criterion4} to a morphism $w:A\times \Upsilon (\ast )\rightarrow nCAT'$
restricting to $\psi$ on $A\times \{ 0\}$ and restricting to the constant object
$\underline{Hom}(E,B)$ (pulled back to $A$) on $A\times \{ 1\}$. More precisely
there is a morphism
$$
w': A\times \Upsilon ^2(\ast , E)\rightarrow nCAT '
$$
restricting to $u$ over the edge $(02)$, and restricting to the universal
morphism $\nu$ (cf \ref{criterion1}) over the edge $(12)$. The restriction to
the edge $(01)$ is the morphism $w$.

Now $w$ is an element of $Hom (\psi , B)$, so by hypothesis there is a diagram
$$
g: A\times \Upsilon ^2(\ast , \ast )\rightarrow nCAT'
$$
sending the  edge $(01)$ to $\epsilon$ and sending the edge $(02)$ to
$w$. Putting this together with the diagram $w'$ and using the fibrant property
of $nCAT'$ (i.e. composing these together) we obtain existence of a diagram
$$
A\times \Upsilon ^3(\ast , \ast , E)\rightarrow nCAT'
$$
restricting to $g$ on the face $(012)$ and restricting to $w'$ on the face
$(023)$. The face $(013)$  yields a diagram
$$
A\times \Upsilon ^2(\ast , E)\rightarrow nCAT'
$$
restricting to $\epsilon$ on the first edge and restricting to our original
morphism $u$ on the edge $(03)$: this is the morphism $v$ we are looking for.
\end{parag}

\begin{parag}
\label{uniquen}
{\em Uniqueness of these factorizations}
The homotopy uniqueness property for factorization of morphisms implies a
similar property for the factorizations of $E$-morphisms obtained in the
previous paragraph.
Suppose that we are given a morphism
$$
u:A \times \Upsilon (E) \rightarrow nCAT'
$$
as above, and suppose that we are given two extensions
$$
v_1, v_2: A\times \Upsilon ^2(\ast , E) \rightarrow nCAT'
$$
restricting to our given morphism on the edge $(02)$, and restricting to
$\epsilon$ on the edge $(01)$. We can complete the $v_i$ to diagrams
$$
z_i: A\times \Upsilon ^3(\ast , \ast , E)\rightarrow nCAT'
$$
restricting to $v_i$ on the faces $(013)$ and restricting to the universal
morphism $\nu$ of \ref{criterion1} on the edge $(23)$. To do this, use the
universal property of $\nu$ (cf \ref{criterion4}) to fill in the faces $(023)$
and $(123)$; then we have a map defined on the shell and by the universal
property of $\nu$ which gives shell extension (\ref{criterion3.5},
\ref{criterion4}) we can extend to the whole tetrahedron.

Note furthermore that we can assume that the restrictions to the faces $(023)$
are the same for $z_1$ and $z_2$ (because we have chosen these faces using only
the map $u$ and not refering to the $v_i$).
Call these $r_{023}(z)$. In particular the restrictions to
$(02)$ give the same map $w:\psi \rightarrow \underline{Hom}(E,B)$.
Now the restrictions of the above diagrams $z_i$ to the faces $(012)$ give two
different factorizations of this map $w$
so by hypothesis there is a homotopy between these factorizations:
it is a morphism
$$
A \times \Upsilon ^2(\ast , \ast ) \times \overline{I}
\rightarrow nCAT'
$$
restricting to $r_{012}(z_i)$ on the endpoints $i=0,1$ of $\overline{I}$,
restricting to the pullback of $\epsilon$ along $(01)\times \overline{I}$
and restricting to the pullback of our morphism $w$ along $(02)\times
\overline{I}$. We can attach this homotopy to the constant homotopy
which is the pullback of $r_{023}(z)$ from
$A\times \Upsilon ^2(\ast , E)$ to $A\times \Upsilon ^2(\ast , E)\times
\overline{I}$. We obtain a homotopy defined on the union of the faces $(012)$
and $(023)$ and going between $z_0$ and $z_1$. Using the fact that the
inclusion of this union of faces into the tetrahedron is a trivial cofibration
(see the list \ref{trivinclusions} above) we get that the inclusion (written
in an obvious shorthand notation where $(0123)$ stands for $A \times \Upsilon
^3(\ast , \ast , E)$ and $(012 + 023)$ for the union of the two faces)
$$
(0123) \times \{ 0\} \cup ^{(012 + 023)\times \{ 0\} }
(012 + 023)\times \overline{I}
\cup ^{(012 + 023)\times \{ 1\} }
(0123) \times \{ 1\}
\hookrightarrow (0123)\times \overline{I}
$$
is a trivial cofibration.  We have a map from the left side into $nCAT'$
so it extends to a map
$$
A \times \Upsilon ^3(\ast , \ast
, E) \times \overline{I} \rightarrow nCAT'.
$$
The restriction of this map to the face $(013)$ is a homotopy
$$
A \times \Upsilon ^2(\ast
, E) \times \overline{I} \rightarrow nCAT'
$$
between our factorizations $v_1$ and $v_2$.
\end{parag}

\subsubnumero{The  relative case}
To actually prove the lemma, we need to obtain a factorization property as
above in the relative situation $E'\subset E$ where we already have the
factorization over $E'$ and we would like to extend to $E$. This is where we
use the homotopy uniqueness of factorization which was in the hypothesis of the
lemma (we use it in the form given in the previous paragraph \ref{uniquen}).
Suppose we are given
$$
v': A\times \Upsilon ^2(\ast , E') \rightarrow nCAT'
$$
restricting to $\epsilon$ on the first edge,
and suppose we are given
$$
u:A \times \Upsilon (E) \rightarrow nCAT'
$$
restricting to $\psi$ on $A\times \{ 0\} $ and restricting to a constant object
$B\in nCAT'$ (i.e. to the pullback $B_A$) on $A\times \{ 1\}$.
Suppose that the restriction of $u$ to $A\times \Upsilon (E')$ is equal to
the restriction of $v'$ to the edge $(02)$. By \ref{existen}
there exists an extension
$$
v_0: A\times \Upsilon ^2(\ast , E) \rightarrow nCAT'
$$
which restricts to $\epsilon$ on the first edge and to $u$ on the edge $(02)$.
Let $v'_0$ denote the restriction of $v_0$ to
$A\times \Upsilon ^2(\ast , E')$.  By the uniqueness statement \ref{uniquen}
there exists a homotopy
$$
A\rightarrow \Upsilon ^2(\ast , E') \times \overline{I} \rightarrow nCAT'
$$
between $v'_0$ and $v'$, constant along edges $(01)$ and $(02)$.
Let $D$ be the coproduct of $\Upsilon ^2(\ast , E')$ and $\Upsilon (E)$
with the latter attached along the edge $(02)$ (i.e. the coproduct is taken
over the copy of $\Upsilon (E')\subset \Upsilon ^2(\ast , E')$ corresponding to
the edge $(02)$).  Our homotopy glues with the constant map $u$
to give a morphism
$$
A \times  D \times \overline{I} \rightarrow nCAT',
$$
and this glues with $u$ to obtain
$$
A\times \Upsilon ^2(\ast , E)\times \{ 0\} \cup ^{A \times  D \times
\{ 0\} } A\times D \times I \rightarrow nCAT'.
$$
The inclusion
$$
A \times  D \times
\{ 0\} \hookrightarrow A\times D \times I
$$
is a trivial cofibration so the inclusion
$$
A\times \Upsilon ^2(\ast , E)\times \{ 0\} \cup ^{A \times  D \times
\{ 0\} } A\times D \times I
\hookrightarrow
A\times \Upsilon ^2(\ast , E)\times \overline{I}
$$
is a trivial cofibration and by the fibrant property of $nCAT'$ there exists an
extension of the above morphism to a morphism
$$
A\times \Upsilon ^2(\ast , E)\times \overline{I}
\rightarrow nCAT'.
$$
The value of this over $1\in \overline{I}$ is the extension
$$
v: A\times \Upsilon ^2(\ast , E) \rightarrow nCAT'
$$
we are looking
for: it restricts to $\epsilon$ on the edge $(01)$, it restricts to $u$ on the
edge $(02)$, and it restricts to $v'$ over
$A\times \Upsilon ^2(\ast , E')$.  This completes the proof of Lemma
\ref{easydirect}.
\eop

\begin{remark}
\label{easyinverse}
One also obtains a criterion similar to \ref{easydirect} for inverse limits in
$nCAT'$. The proof is the same as above but using the universal diagram
$$
B \stackrel{E}{\rightarrow} B\times E
$$
in the place of $\underline{Hom}(E,B)\stackrel{E}{\rightarrow} B$.
We did not choose to use this in the proof of \ref{inverse} because
it didn't seem to make any substantial savings (and in fact probably would
have complicated the notation in many places).
\end{remark}

\begin{parag}
\label{alpha}
We now improve the above criterion with a view toward applying this in the
argument \ref{fixup} given above for the case of sets. Fix a functor $\psi :
A\rightarrow nCAT'$ of $n+1$-categories. Suppose $\alpha$ is an infinite
cardinal number such that $$
\# A \leq \alpha
$$
and
$$
\# \psi (a) \leq \alpha
$$
for all $a\in A$.
\end{parag}

\begin{parag}
\label{factorization}
Suppose $B\in nCAT'$ and suppose $u: \psi \rightarrow B$ is a morphism.
Then we claim that there is $B'\in nCAT'$ with $\#^{\rm pre} B' \leq
\alpha$, and
with a factorization $\psi \rightarrow B' \rightarrow B$.

To prove this notice that $u$ is a morphism
$$
u:A\times I \rightarrow nCAT',
$$
and the image of $u$ is contained in an $\alpha$-bounded set of additions of
trivial cofibrations to $nCAT$ (recall that $nCAT'$ was constructed by
adding pushouts along trivial cofibrations to $nCAT$ \ref{perspective1}).

We can take $B'\subset B$ to be a sub-precat containing all of the objects
necessary for the morphisms involved in the trivial cofibrations which are
added in the previous paragraph, as well as the morphisms involved in $u$
given that we fix $\psi$, all of the objects necessary for the structural
morphisms of a precat, and finally add on what is necessary to get $B'$ fibrant.
This has cardinality $\# ^{\rm pre}B' \leq \alpha$.
\end{parag}

\bigskip

\subnumero{A construction}

\begin{parag}
\label{hypothesis}
{\em Hypothesis---}
With the above notations, suppose we have an object $U\in nCAT'$ together with a
morphism $a: \psi \rightarrow U$  provided with the following data:
\newline
(A)---for every morphism $\psi \rightarrow V$ where  $\# ^{\rm pre}
V\leq \alpha$,
a factorization which we call the {\em official factorization}
$$
\psi \stackrel{a}{\rightarrow} U \rightarrow V
$$
(in other words a diagram
$$
A\times \Upsilon ^2(\ast , \ast ) \rightarrow nCAT'
$$
restricting to $a$ on the edge $(01)$ and restricting to our given morphism
on the edge $(02)$);
\newline
(B)---for every diagram
$$
\psi \rightarrow V\rightarrow V',
$$
a
completion of this and the official factorization diagrams
$$
\psi \rightarrow U\rightarrow V,\;\; \psi \rightarrow U\rightarrow V'
$$
to a diagram (the {\em official commutativity diagram})
$$
\psi \rightarrow U \rightarrow V \rightarrow V'
$$
(which again means a morphism
$$
A\times \Upsilon ^3(\ast , \ast ,\ast )\rightarrow
nCAT'
$$
restricting to our given diagrams on the faces $(023)$, $(012)$, and
$(013)$).
\end{parag}

\begin{parag}
\label{construction1}
Keep the above hypothesis \ref{hypothesis}.
Let
$$
[b,i]:\psi \stackrel{b}{\rightarrow} U' \stackrel{i}{\rightarrow} U
$$
be a
factorization of the morphism $a$ as above (\ref{factorization}) with  $\# U'
\leq \alpha$.
This means a diagram whose third edge $(02)$ is equal to $a$.

Unfortunately at
this point we have no control over the choice of $U'$, so the ``real'' $U'$
which we would like to choose to satisfy the criterion of \ref{easydirect} may
be a direct factor of this $U'$. To explain this notice that by hypothesis
\ref{hypothesis} (A) there is a morphism
$$
q:U\rightarrow U'
$$
giving a factorization
$$
[a,q]:\psi \rightarrow U \rightarrow U'.
$$
Let $b$ be the edge $(02)$ of this diagram, so we can write
$b\sim qa$.

Using the fibrant property of $nCAT'$ we can glue the diagrams $[b,i]$ and
$[a,q]$ together to give a diagram
$$
[b,i,q]: \psi \rightarrow U' \rightarrow U \rightarrow U',
$$
in other words  a morphism
$$
A\times \Upsilon ^3(\ast , \ast , \ast ) \rightarrow nCAT'
$$
restricting to $[b,i]$ on the face $(012)$ and restricting to
$[a,q]$ on the face $(023)$ (and satisfying the usual condition that the
restriction to the face $(123)$ be constant in the $A$ direction).
Denote by $p$
the restriction  to the edge $(13)$, and denote by $[b,p]$ the restriction to
the face $(013)$. Thus
$$
[b,p]: \psi \rightarrow U' \rightarrow U'
$$
is a diagram whose restrictions to the edges $(01)$  and $(02)$ are both equal
to the morphism $b$.

Restriction to the face $(123)$ is a diagram $[i,q]$ with third edge equal to
$p$, in other words we can write $p \sim q\circ i $.

The official commutativity diagram for $[b,p]$ is a diagram of the form
$$
[a,q,p]: \psi \rightarrow U \rightarrow U' \rightarrow U'.
$$
The restriction of this diagram to the face $(023)$ is the diagram $[b,p]$.
The restrictions to $(012)$ and $(013)$ are both equal (by hypothesis (B))
to the official factorization diagram $[a,q]$. In particular, the face $(123)$
gives a diagram
$$
[q,p]: U\rightarrow U' \rightarrow U'
$$
whose third edge (which we should here denote $(13)$) is again the morphism $p$.
Homotopically we get an equation
$$
p \circ q \sim q.
$$

In view of the fact that $p\sim  q \circ i$  we get
$$
p\circ p \sim p.
$$
This equation says that, up to homotopy, $p$ is a projector. It is the
projector onto the answer that we are looking for.
\end{parag}

\begin{parag}
\label{construction2}
{\em Construction---}
Continuing with hypothesis \ref{hypothesis} and the
notations of \ref{construction1}, we will construct the object corresponding to
the ``image'' of the homotopy projector $p$.  To do this we will take the
``mapping telescope'' of the sequence
$$
U' \stackrel{p}{\rightarrow } U' \stackrel{p}{\rightarrow }
U' \stackrel{p}{\rightarrow}\ldots .
$$
In the present setting of $n$-categories we do this as follows
(which is basically just the mapping telescope in the closed model category
structure of \cite{nCAT}). Recall that $\overline{I}$ is the $1$-category with
two objects $0,1$ and two morphisms inverse to each other between the
objects. We
consider it as an $n$-category. Glue together  the $n$-precats $U' \times
\overline{I}$, one for each natural number, by attaching $U'\times \{ 1\}$ in
the $i-1$-st copy to $U' \times \{ 0\}$ in the $i$-th copy via the map $p:
U'\times \{ 1\} \rightarrow U' \times \{ 0\}$. Denote by $T'$ the resulting
$n$-precat and by $T'\hookrightarrow T$  a fibrant replacement. Inclusion of $U'
\times \{ 0\}$ in the first copy gives a morphism $$
j: U' \rightarrow T.
$$
On the other hand, using the projection $p$
in each variable and the homotopy $p\circ p\sim p$ gives a morphism
$$
r : T\rightarrow U'
$$
(which comes by extension from a map $r':T'\rightarrow U'$)
and we have $rj=p$.
\end{parag}

\begin{parag}
\label{construction3}
{\em Claim:}
The morphism $jr : T \rightarrow T$ is homotopic to the identity,
via a homotopy compatible with the homotopy $p\circ p \sim p$.

This is by a classical construction that works in any closed model category
with ``interval object'' such as $\overline{I}$.  As a sketch of proof, let
$T^{m}$ denote the subobject of $T'$ obtained by taking only the first
$m$ copies of $U'\times \overline{I}$. Then $T^m$ retracts to the last
copy of $U'$, so the restriction  of $r'$ to $T^m$ is
homotopic (via this retraction) to $p$. On the other hand, the inclusion $T^m
\hookrightarrow T^{m+1}$ is also homotopic to $p$ (via the
retractions to the end copies of $U'$).
Thus we may choose a homotopy (in Quillen's sense cf \ref{homotopic1})
between the
restriction of $r'$ to $T^m$, and
the inclusion $T^m \hookrightarrow T^{m+1}$.  We can make this into a homotopy
between
$$
r', 1_{T^m} : T^m \tworightarrows T,
$$
and since $T$ is fibrant we can do this with a homotopy using the interval
$\overline{I}$. Again because $T$ is fibrant we can assume that these
homotopies are compatible for all $m$, so they glue together to give a homotopy
between the two maps
$$
r', 1_{T'} : T'\tworightarrows T.
$$
Then extend from $T'$ to $T$.
\end{parag}

\begin{parag}
\label{construction4}
We wrap things up by pointing out how $T$ fits in with the situation of
\ref{construction1}.
Consider the sequence of morphisms
$$
\psi \rightarrow U' \rightarrow T \rightarrow U' \rightarrow T.
$$
The composition of the first two gives a morphism $jb:\psi \rightarrow T$.
The composition of the first three morphisms is equal to
$rjb\sim pb$ which has a homotopy to the usual
morphism $b:\psi \rightarrow U'$. Thus the morphism $b$ factors through $T$.
Finally from our claim \ref{construction3} the composition of the last two
arrows is homotopic to the identity on $T$.

Our original morphism $\psi \rightarrow U$ factors through $U'$ hence it
factors through $T$: the composition
$$
\psi \rightarrow T \rightarrow U' \rightarrow U
$$
is equal to the original morphism $a: \psi \rightarrow U$.

We have the morphism
$$
jq: U\rightarrow T
$$
providing a factorization
$$
\psi \rightarrow U \rightarrow T.
$$
The composition $T\rightarrow U \rightarrow T$ is homotopic to the identity
on $T$ by claim \ref{construction3}.
\end{parag}

\begin{lemma}
\label{construction5}
Under hypothesis \ref{hypothesis} and with the above notations,
the morphism $\psi \rightarrow T$ has the unique homotopy factorization property
of \ref{easydirect} with respect to any morphism $\psi \rightarrow B$ (without
bound on the cardinality of $B$).
\end{lemma}
{\em Proof:}
This is  really only a
statement about $1$-categories. We can consider the $1$-category $M$ which is
the truncation of the $n+1$-category of objects under $\psi$ (cf \ref{lazy}
below). Our objects $U,U',T$ and so on togther with maps from $\psi$ may be
considered as objects in the category $M$. The result of \ref{construction1}
says that $p: U'\rightarrow U'$ is a projector in the category $M$, and in
\ref{construction2}, \ref{construction3} and \ref{construction4} we show
that the
object $T$ corresponding to this projector exists. The criterion of
\ref{easydirect} asks simply that $T$ be an initial object in $M$.

What we know from hypothesis \ref{hypothesis} is that $T$ is provided with a
collection of morphisms $T\rightarrow B'$ to every $\alpha$-bounded object of
$M$, in such a way that these form a natural transformation from the constant
functor $T$ to the identity functor $M_{\alpha}\rightarrow M$ (where
$M_{\alpha}$ is the full subcategory of objects having cardinality bounded by
$\alpha$).

The fact that $T$ is the object corresponding to the projector $p$ (and
that $p$ was the projector defined by the natural transformation for $U'$) means
that the value of this natural transformation on $T$ itself is the identity.

Suppose $\psi \rightarrow B$ is an object of $M$. Then there is a factorization
through
$\psi \rightarrow B'\rightarrow B$ with $\# ^{\rm pre}B' \leq \alpha$.
This just says that every object of $M$ has a morphism from an object in
$M_{\alpha}$. It is worth mentioning that if $B\in M$ and if $B'\rightarrow B$
and $B'' \rightarrow B$ are two morphisms from objects in $M_{\alpha}$ then
they both factor through a common morphism $B'''\rightarrow B$ from an object
in $M_{\alpha}$.

Using the above formal properties, we show that $T$ is an initial object of
$M$ to prove the lemma.
Suppose $B$ is an object of $M$. There exists a morphism $B'\rightarrow B$
from an object of $M_{\alpha}$ so applying our natural transformation, there
exists a morphism $T\rightarrow B'$ and hence a morphism $T\rightarrow B$.
Suppose that $T\tworightarrows B$ is a pair of morphisms. These factor through
a common object of $M_{\alpha}$
$$
T\tworightarrows B' \rightarrow B,
$$
and applying our natural transformation we obtain that the compositions
of the two morphisms
$$
T\rightarrow T\tworightarrows B'
$$
are equal to the given morphism  $T\rightarrow B'$; however, since the
natural transformation $T\rightarrow T$ is the identity, this implies that our
two morphisms $T\tworightarrows B'$ were equal and hence that the two original
morphisms $T\tworightarrows B$ were equal.
\eop

\begin{corollary}
\label{construction6}
In the situation of Lemma \ref{construction5} the map $\psi \rightarrow T$
is a direct limit.
\end{corollary}
{\em Proof:}
By \ref{construction5} it satisfies the condition of \ref{easydirect}
so by the latter, it is a direct limit.
\eop

\bigskip

\subnumero{Proof of Theorem \ref{direct}}

\begin{parag}
\label{lazy}
{\em Objects under $\psi$:}
In order to replicate the proof that was given above for the category of sets,
we need to know what the category of ``objects under $\psi$'' is.
Suppose $C$ is
an $n+1$-category and $A$ another $n+1$-category and suppose $\psi :
A\rightarrow
C$ is a morphism.
We define the $n+1$-category $\psi /C$ of {\em
objects under $\psi$}
to be the category of morphisms
$$
(A\times I)\cup ^{A\times \{ 1\} }\{
1\}\rightarrow C
$$
restricting to $\psi$ on $A\times \{ 0\}$.  In other words,
it is the fiber of the morphism
$$
\underline{Hom}((A\times I)\cup ^{A\times \{ 1\} }\{
1\},C)\rightarrow \underline{Hom}(A,C)
$$
over $\psi$.
\end{parag}

\begin{parag}
Let $(\psi /C)_{\alpha}$ denote the category of objects under $\psi$
which are (set-theoretically speaking) contained in a given fixed set of
cardinality $\alpha$.  It has cardinality $\leq 2^{\alpha}$. It is a full
subcategory of $\psi/C$.
\end{parag}

\begin{parag}
\label{reeasydirect}
We can restate the criterion of \ref{easydirect} in terms of the above
definition. Let $\tau _{\leq 1}(\psi /C)$ denote the $1$-truncation of the
category of objects under $\psi$ defined in \ref{lazy}. It is a $1$-category.
The criterion says that if $u: \psi \rightarrow U$ is an initial object in
this category then it is (the image under the truncation operation of) a
direct limit of $\psi$.
Definition \ref{lazy} and the present remark were used in the proof of
\ref{construction5} already, where we denoted $\psi /C$ by $M$.
\end{parag}

{\em Proof of Theorem \ref{direct}:}
Suppose $\psi : A\rightarrow nCAT'$. Fix a cardinal $\alpha$ bounding $(A,\psi
)$ as above. Let $M_{\alpha}:= (\psi /C)_{\alpha}$
denote the $n+1$-category of objects
under $\psi$,
of cardinality bounded by $\alpha$.  Let $U$ be the inverse
limit of the forgetful functor $f:M_{\alpha}\rightarrow nCAT'$,
given by Theorem \ref{inverse}.
By Corollary \ref{constant}, the pullback of $U$ to a constant family $U_A$ over
$A$ is again an inverse limit of the functor
$$
f_A: A\times M_{\alpha}\rightarrow nCAT'
$$
($f$ pulled back along the second
projection to $M_{\alpha}$).

We have a morphism of families over $A\times M_{\alpha}$, from $\psi$ to $f_A$,
which thus factorizes into
$$
\psi \rightarrow U_A \rightarrow f_A.
$$

The morphism  $\psi \rightarrow U_A$
is automatically provided with the data
required for Hypothesis \ref{hypothesis}.

Apply the above construction \ref{construction1}--\ref{construction4} to obtain
$\psi \rightarrow T_A$, and Lemma \ref{construction5} and Corollary
\ref{construction6} show that  $\psi \rightarrow T_A$ is a direct limit of
$\psi$.
\eop

\numero{Applications}

We will discuss several different possible applications for the notions of
inverse and direct limit in $n$-categories in general, and of the existence of
limits in $nCAT'$ in particular. Many of these applications are only
proposed as conjectural. Only in the first section do we give full proofs.

The conjectures are for the most part supposed to be possible to do with the
present techniques, except possibly \ref{mts}.

\bigskip

\subnumero{Coproducts and fiber products}

\begin{parag}
\label{fiberproducts}
Taking $A$ to be the category with three objects
$a$, $b$ and $c$ and morphisms $a\rightarrow b$ and $c\rightarrow b$,
a functor $A\rightarrow nCAT$ is just a triple of $n$-categories
$X,Y,Z$ with maps $u:X\rightarrow Y$ and $v:Z\rightarrow Y$. The inverse limit
of the projection into $nCAT'$ is the {\em homotopy fiber product}
denoted $X\times ^{\rm ho} _YZ$.
\end{parag}

\begin{lemma}
\label{calcinverse}
Suppose $A$ is as above and $\varphi : A\rightarrow nCAT$ is a morphism
corresponding to a pair of maps $u:X\rightarrow Y$ and $v:Z\rightarrow Y$ of
$n$-categories such that $u$ is fibrant. Then the usual fiber product $X\times
_YZ$ is a limit of $\varphi$ so we can write
$$
X\times ^{\rm ho} _YZ = X\times  _YZ.
$$
\end{lemma}
{\em Proof:}
One way to prove this is to use our explicit construction of the inverse limit
(\ref{lambda}).  The second way is to show that $U:=X\times  _YZ$
satisfies the required universal property as follows.  First of all note
that the commutative square
$$
\begin{array}{ccc}
U&\rightarrow & X \\
\downarrow && \downarrow \\
Z &\rightarrow & Y
\end{array}
$$
corresponds to a map $I\times I\rightarrow nCAT$ which we can project into
$I\times I\rightarrow nCAT'$. Then combine this with the projection
$$
A \times I \rightarrow I \times I
$$
which sends $A\times \{ 0\}$ to $(0,0)$ and sends $A\times \{ 1\}$
to the copy of $A\subset I\times I$ corresponding to the sides
$(1,01)$ and $(01,1)$ of the square. We get a map
$A\times I \rightarrow nCAT'$ having the required constancy property to
give an element $\epsilon \in Hom (U, \varphi )$.  This is the map which
we claim is an inverse limit.

In passing note that since $X,Y,Z$ are elements of $nCAT$ they are by
definition fibrant, and since by hypothesis the map $X\rightarrow Y$ is
fibrant, the map $U\rightarrow Z$ is fibrant too, and so $U$ is fibrant.

We now fix a fibrant $n$-category $V$ and study the functor which to an
$n$-precat $F$ associates the set of morphisms
$$
g:A\times \Upsilon (F)\rightarrow nCAT'
$$
with $r_0(g)=V_A$ and $r_1(g)= \varphi $.  This is of course just the
functor represented by $Hom (V, \varphi )$.  Recalling that
$\underline{Hom}(V,U)$ is the morphism set in $nCAT$, we obtain by composition
with $\epsilon $ a morphism
$$
C_{\epsilon}: \underline{Hom}(V,U)\rightarrow Hom (V,\varphi ).
$$
In this case, since $\epsilon$ comes from $nCAT$ in which the composition at
the first stage is strict, the morphism $C_{\epsilon}$ is
strictly well defined rather than being a weak morphism as usual in the
notion of limit. We would like to show that $C_{\epsilon}$ is an equivalence
(which would prove the lemma).

A morphism $g:A\times \Upsilon (F)\rightarrow nCAT'$ decomposes as a pair of
morphisms $(g_1,g_2)$ with
$$
g_i : I \times \Upsilon (F) \rightarrow nCAT';
$$
in turn these decompose as pairs $g_i^+$ and $g_i^-$
where
$$
g_i ^+ : \Upsilon (\ast , F)\rightarrow nCAT',
$$
$$
g_i^- : \Upsilon (F,\ast )\rightarrow nCAT' .
$$
(Decompose the square $I \times \Upsilon (F)$ into two triangles,
drawing the edge $I$ vertically with vertex $0$ on top.)
The conditions on everything to correspond to a morphism $g$ are that
$$
r_{12}(g_i^+) = r_{01}(g_i^-)
$$
and
$$
r_{02}(g_1^-)= r_{02}(g_2^-).
$$
The endpoint conditions on $g$ correspond to the conditions
$$
r_{12}(g_1^-)= u,\;\;\;\;
r_{12}(g_2^-) = v,
$$
and
$$
r_{01}(g_i^+)=1_V.
$$
Putting these all together we see that our functor of $F$ is
of the form a fiber product of four diagram $n$-categories
\ref{diag1}.  More precisely, put
$$
M_u :=Diag ^{1_V, X}(\ast ,\underline{\ast}; nCAT')
\times _{Diag ^{V, X}(\underline{\ast} ; nCAT')}
Diag ^{V, u}(\underline{\ast}, \ast ; nCAT')
$$
where the morphisms in the fiber product are  $r_{12}$ then $r_{01}$;
and define $M_v$ similarly.
Then
$$
Hom (V, \varphi )= M_u \times _{Diag ^{V,Y}(\underline{\ast} ; nCAT')}
M_v,
$$
where here the morphisms in the fiber product are the restrictions
$r_{02}$ on the second factors of the $M$.

Refer now to the calculation of \ref{diag1.6} in view of
the comparison result \ref{diag3} (applied to $nCAT\rightarrow nCAT'$).

By this calculation the restriction morphism
$$
r_{12}:
Diag ^{1_V, X}(\underline{\ast}, \ast ; nCAT') \rightarrow
Diag ^{V,Z}(\underline{\ast} ; nCAT')
$$
is a fibrant equivalence. Therefore the second projections are equivalences
$$
M_u \rightarrow Diag ^{V,u}(\ast , \underline{\ast}; nCAT')
$$
and similarly for $v$.
Using these second projections in each of the factors $M$ we get an equivalence
$$
Hom (V, \varphi )\rightarrow
$$
$$
Diag ^{V,u}( \underline{\ast},\ast ; nCAT')
\times _{Diag ^{V,Y}(\underline{\ast} ; nCAT')}
Diag ^{V,u}(\underline{\ast},\ast ; nCAT')
$$
where the morphisms in the fiber product are $r_{02}$. There is a morphism
from the same fiber product taken with respect to $nCAT$, into here.
In the case of the fiber product taken with respect to $nCAT$ the calculation
of \ref{diag1.6} gives directly that it is equal to
$$
\underline{Hom}(V,X)\times _{\underline{Hom}(V,Y)}\underline{Hom}(V,Z)
$$
which is just $\underline{Hom}(V,U)$.  The morphism
$$
\underline{Hom}(V,X) =
Diag ^{V,u}( \underline{\ast},\ast ; nCAT)
\rightarrow
Diag ^{V,u}( \underline{\ast},\ast ; nCAT')
$$
is an equivalence by \ref{diag3}, and similarly for the other factors
in the fiber product.

Now we are in the general situation that we
have  equivalences of fibrant $n$-precats
$P\rightarrow P'$, $Q\rightarrow Q'$ and $R\rightarrow R'$ compatible with
diagrams
$$
P\rightarrow Q\leftarrow R, \;\;\; P'\rightarrow Q'\leftarrow R'.
$$
If we know that the morphisms $P\rightarrow Q$ and $P'\rightarrow Q'$
are fibrant then we can conclude that  these induce an equivalence
$$
P\times _QR\rightarrow P'\times _{Q'} R'.
$$
Prove this in several steps using \ref{pushoutA} and \cite{nCAT} Theorem 6.7:
$$
P'\times _{Q'}R' \stackrel{\sim}{\rightarrow} P' \times _{Q'}R = (P'\times
_{Q'}Q)\times _Q R
$$
and
$$
P'\times _{Q'}Q \stackrel{\sim}{\rightarrow} P'
$$
so
$$
P \stackrel{\sim}{\rightarrow} P' \times _{Q'}Q
$$
giving  finally
$$
P \times _Q R \stackrel{\sim}{\rightarrow} (P'\times _{Q'}Q)\times _Q R;
$$
then apply (\ref{explaincmc}, CM2).

Applying this general fact to the previous situation gives that the morphism
$$
\underline{Hom}(V,X)\times _{\underline{Hom}(V,Y)}\underline{Hom}(V,Z)
\rightarrow
$$
$$
Diag ^{V,u}( \underline{\ast},\ast ; nCAT')
\times _{Diag ^{V,Y}(\underline{\ast} ; nCAT')}
Diag ^{V,u}(\underline{\ast},\ast ; nCAT')
$$
is an equivalence. By (\ref{explaincmc}, CM2) this implies that
$$
C_{\epsilon} : \underline{Hom}(U,V)\rightarrow Hom (V,\varphi )
$$
is an equivalence.
\eop

Lemma \ref{calcinverse} basically says that for calculating homotopy fiber
products we can forget about the whole limit machinery and go back to our usual
way of assuming that one of the morphisms is fibrant.

\begin{parag}
\label{coproducts}
Taking $A$ to be the opposite of the category in the
previous paragraph, a functor $A\rightarrow nCAT$ is a triple $U,V,W$ with
morphisms $f:V\rightarrow U$ and $g:V\rightarrow W$. The direct limit is
the {\em
homotopy pushout} of $U$ and $W$ over $V$, denoted
$U\cup ^V_{\rm ho}W$.
\end{parag}

\begin{lemma}
\label{calcdirect}
Suppose that $f$ is cofibrant.  Let $P$ denote a fibrant replacement
$$
U\cup ^VW \hookrightarrow P.
$$
Then there is a natural morphism $u: \varphi \rightarrow P$ which is a
direct limit.
Thus we can say that the morphism
$$
U\cup ^VW \rightarrow U\cup ^V_{\rm ho}W
$$
is a weak equivalence, or equivalently that the morphism of $n$-categories
$$
Cat(U\cup ^VW )\rightarrow U\cup ^V_{\rm ho}W
$$
is an equivalence.
\end{lemma}

The proof is similar to the proof of \ref{calcinverse} and is left as an
exercise.

This lemma provides justification {\em a posteriori}
\footnote{And therefore running a certain risk of being circular\ldots }
for
having said in \cite{nCAT} that $Cat(U\cup ^VW )$ is the ``categorical pushout
of $U$ and $W$ over $V$''.  It also shows that this pushout, which occurs in the
generalized Seifert-Van Kampen theorem of \cite{nCAT}, is the same as the
homotopy pushout.

\bigskip

\subnumero{Representable functors and internal $Hom$}

Suppose $A$ is an $n+1$-category. Recall that $A^o$ is the first opposed
category obtained by switching the directions of the $1$-arrows but not the
rest (this comes from the inversion automorphism on the first simplicial factor
of $\Delta ^{n+1}$).
The ``arrow family'' is a family
$$
Arr(A): A^o\times A \rightarrow nCAT',
$$
associating to $X\in A^o$ and $Y\in A$ the $n$-category $A_{1/}(X,Y)$.
We will not discuss here the existence and uniqueness of this family
(there is not actually a natural way to define this family in Tamsamani's point
of view on $n$-categories, so it must be done by constructing the family by
hand making choices of various morphisms when necessary).

The arrow family gives two functors
$$
\alpha : A \rightarrow \underline{Hom}(A^o, nCAT')
$$
and
$$
\beta : A^o \rightarrow \underline{Hom}(A, nCAT').
$$

\begin{conjecture}
That $\alpha$ and $\beta$ are fully faithful (as is the case for $n=0$).
\end{conjecture}

We say that an object of $\underline{Hom}(A^o, nCAT')$ (resp.
$\underline{Hom}(A, nCAT')$ is {\em
representable} if it comes from an object of $A$ (resp. $A^o$).
Note that such objects are themselves functors $A\rightarrow nCAT'$ or
$A^o\rightarrow nCAT'$, and we call them {\em representable functors}.

\begin{conjecture}
\label{representable}
Suppose that an $n+1$-category $A$ admits arbitrary direct and inverse limits.
Then  a functor $h: A^o\rightarrow nCAT'$ is representable by an object of $A$
if and only it transforms direct limits into inverse limits.
A functor $g: A \rightarrow nCAT'$ is representable by an object of $A^o$ if
and only if it transforms inverse limits into inverse limits.
\end{conjecture}

\begin{parag}
\label{internal}
If this conjecture is true we would obtain the following corollary:
that if an $n+1$-category  $A$ admits arbitrary direct and inverse limits,
then $A$ has an internal $\underline{Hom}$.
To see this, fix objects $x,y\in A$. Denote by $\times $ the functor $A\times
A\rightarrow A$ which associates to $(u,v)$ the direct product of $u$ and $v$
(considered as an inverse limit). This functor comes from
Theorem \ref{variation} as described in (\ref{usevar}).

Now the functor $u\mapsto Arr (A)(x\times u, y)$ from $A^o$ to $nCAT$
transforms direct limits to inverse limits
(this uses one direction of Conjecture \ref{representable}, and I suppose
without
proof that the functor $u\mapsto x\times u$ is known to preserve direct
limits).

Therefore by (the other direction of) Conjecture \ref{representable},
the functor $u\mapsto Arr (A)(x\times u, y)$ is representable by
an object $\underline{Hom}_A(x,y)$.
\end{parag}

\begin{parag}
\label{topos}
We are obviously going toward some sort of theory of {\em $n$-topoi}:
an {\em $n$-topos} would be an $n$-category admitting arbitrary direct and
inverse limits (indexed by small $n$-categories). There may be some other
conditions that one would have to impose...
\end{parag}

\bigskip

\subnumero{$n$-stacks}

Suppose $\Xx$ is a site. Consider the underlying category as an $n+1$-category.
An {\em $n$-stack over $\Xx$} is a morphism $F: \Xx \rightarrow nCAT'$ such that
for every object $X\in \Xx$ and every sieve $\Bb \subset \Xx /X$ the morphism
$$
\Gamma (\Xx /X, F|_{\Xx /X})\rightarrow \Gamma (\Bb , F|_{\Bb})
$$
is an equivalence of $n$-categories, where $\Gamma (\Bb , F|_{\Bb})$ denotes
the inverse limit of $F|_{\Bb}$ and the same for
$\Gamma (\Xx /X, F|_{\Xx /X})$.  We define $nSTACK /\Xx $ to be the
full subcategory of the (already fibrant) $n+1$-category $\underline{Hom}(\Xx ,
nCAT')$ whose objects are the morphisms $F$ satisfying the above criterion.

The $n+1$-category $nSTACK /\Xx$ admits inverse limits---since the only thing
involved in its definition is an inverse limit and inverse limits commute
with each other. In particular we may speak of {\em homotopy fiber products}
of $n$-stacks.

\begin{conjecture}
\label{projeffstacks}
Homotopy projectors are effective for $n$-stacks, in other words given an
$n$-stack $U'$ with endomorphism $p$ such that $p\circ p \sim p$,
the ``telescope construction'' $T$ of \S 5 is again an $n$-stack.
\end{conjecture}

Assuming this conjecture, the same argument as in \S 5 would work to show that
$nSTACK /\Xx$ admits direct limits.

A {\em $n$-prestack over $\Xx$} is just a morphism
$F: \Xx \rightarrow nCAT'$ without any other condition (this makes sense for
any category $\Xx$ and in fact for any $n+1$-category, it is just our notion of
family of $n$-categories indexed by $\Xx$).  We can adopt the notation
$$
nPRESTACK /\Xx := \underline{Hom}(\Xx , nCAT').
$$
Suppose $F$ is an $n$-prestack. We define the {\em associated stack}
denoted $st(F)$ to be the universal $n$-stack to which $F$ maps. Assuming
Conjecture \ref{projeffstacks}, the associated stack  $st(F)$ exists again by
copying the argument of \S 5 above.

\begin{remark}
The inverse limit of a family of stacks is the same as the
inverse limit of the underlying family of prestacks. However this is not true
for direct limits.
\end{remark}

\begin{parag}
By \ref{internal} which is based on Conjecture \ref{representable},
the $n+1$-category $nSTACK /\Xx$ admits internal $\underline{Hom}$.
Using this (or alternatively using a direct construction which associates to
any $X\in\Xx$ the $n+1$-category $nSTACK/(\Xx /X )$) we should be able to
construct the $n+1$-stack $n\underline{STACK}/\Xx$.
\end{parag}

\begin{parag}
\label{geometric1}
Now that we have a notion of $n$-stack not necessarily of groupoids, one can
ask how to generalize the definition of {\em geometricity} given in
\cite{geometricN}, to the case where the values may not be groupoids.

If $A$ is an $n$-category and $X$, $Y$ are
sets with maps $a:X\rightarrow A$ and $b:Y\rightarrow A$ then
the pullback
$$
(a^o, b)^{\ast} (Arr (A)) =:X \times Y\rightarrow (n-1)CAT'
$$
may be considered as an $(n-1)$-category (taking the union over all of the
points of $X\times Y$)
which we denote by $Hom _A(a,b)$. However, it is  no longer the same thing
as the
fiber product $X\times _AY$.  Both of these still satisfy the
recurrence-enabling fact that they are $n-1$-categories.  Thus we can still
employ the same type of definition as in \cite{geometricN}. However, as many
common examples quickly show, the smoothness condition should only be imposed on
the fiber product, not the arrows.  Thus we say that $A$ is {\em geometric}
(resp. {\em locally geometric}) if: \newline
(GS1) for any two morphisms from schemes $a:X\rightarrow A$ and $b:Y\rightarrow
A$, the arrow $n-1$-stack $Hom _A(a,b)\rightarrow X\times Y$ and the
product $X\times _AY$ are both
geometric (resp. locally geometric); and
\newline
(GS2) there exists a smooth morphism from a
finite type scheme (resp. locally finite type scheme) $X\rightarrow A$
surjective on the truncations to $0$-stacks;
\newline
where the morphism $X\rightarrow A$ is said to be {\em smooth} if for any
morphism from a scheme of finite type $Y\rightarrow A$, the locally geometric
$n-1$-stack $X\times _AY$ is actually geometric and is smooth, this latter
condition meaning that the smooth surjection to it from GS2 comes from a smooth
scheme of finite type.
\end{parag}

\begin{parag}
\label{geometric2}
Here is an example to show what we are thinking of
(this type of example---even
if relatively unknown on ``alg-geom''---apparently comes up very often on
``q-alg''). The stack of vector bundles
on a given variety, for example, is locally geometric. It has an additional
operation, tensor product, which allows it to be considered as a monoidal (or
braided or symmetric) monoidal $1$-stack, thus allowing us to consider it as a
$2$, $3$ or $4$-stack.  In these cases there are only one object (locally
speaking) and in the $3$ and $4$ cases, only one morphism (in the $4$ case
only one $2$-morphism). The original stack comes back as an arrow stack
(possibly after iterating).  In this example, if we want a tensor product we
are forced to consider things not of finite type, so the arrow stacks should
often be allowed to be only locally geometric (also one readily sees that the
arrow stacks will not necessarily be smooth). On the other hand the finite type
and smoothness conditions in GS2 correspond in this example to the smoothness
and finite type conditions for the Picard scheme.
\end{parag}

\begin{parag}
\label{locallyP}
Suppose $\Pp$ is a property of $n$-stacks of groupoids.  Then
we say that an $n$-stack $A$ is {\em locally $\Pp$} (and we call this property
$loc \Pp$) if  $F=\tau _{\leq 0}A $ is an filtered inductive limit of open
subsheaves $F_i\subset F$ (the openness condition means that for any scheme
$X\rightarrow F$, $X\times _FF_i$ is an open subset of $X$) such that $A\times
_{F}F_i$ has property $\Pp$.

In particular we obtain notions of {\em locally presentable} and {\em locally
very presentable} $n$-stacks of groupoids.

We claim that for $\Pp = $ ``geometric'' the above definition gives the same
definition as the previous definition of locally geometric.  Suppose
that $A$ is locally $\Pp$.  Let $A_i := A\times _FF_i$. This is an open
substack of $A_i$.  Let $X_i\rightarrow A_i$ be the smooth surjections from
schemes of finite type. Then $X_i\rightarrow A$ is smooth (for example by the
formal criterion for smoothness). Thus the morphism from the disjoint union of
the $X_i$ to $A$ is a smooth surjection proving that $A$ is locally geometric
according to the old definition.

Suppose now that $A$ is locally geometric for the old definition, and let $X_i
\rightarrow A$ be the smooth morphisms from schemes of finite type which
together cover $A$.  Let $F=\tau {\leq 0} A$ and let $F_i \subset F$ be the
images of $X_i$.  Let $A_i = A\times _FF_i$.
It is clear that $X_i$ maps to
$A_i$ by a map which is, on the one hand, smooth by the formal criterion, and
on the other hand surjective on the level of $\pi _0$ by definition.
Thus the $A_i$ are geometric, i.e. have property $\Pp$.  It is clear that the
union of the $A_i$ is $A$. Finally, the $F_i$ are open subsheaves of $F$,
using smoothness of $X_i \rightarrow A$ plus Artin approximation.
\end{parag}

\begin{definition}
\label{extendingproperties}
If $\Pp$ is a property of $n$-stacks of groupoids
(say, independent of $n$...) then we can extend $\Pp$ to a property of
$n$-stacks in a minimal way such that the following conditions hold:
\newline
(A)\,\,\, If $A$ has property $\Pp$ then so does the interior groupoid $A^{\rm
int}$;
\newline
(B)\,\,\, If $A$ has property $\Pp$ and $a: X\rightarrow A$ and $b:
Y\rightarrow A$ are morphisms from schemes of finite type then
$Hom _A(a,b)$ has property $\Pp$.

That such a minimal extension exists is obvious by induction.
\end{definition}

\begin{parag}
Taking the property $\Pp$ in the above definition \ref{extendingproperties}
to be
``locally presentable'' or ``locally very presentable'' or ``locally geometric''
we obtain reasonable properties for  $n$-stacks not necessarily of groupoids.
The use of the locality properties is natural here since the composition
operation will often be something like tensor product, which does not preserve
any substack of finite type.
\end{parag}

\bigskip

\subnumero{The notion of stack, in general}

We give here a very general discussion of the notion of ``stack''.
This was called ``homotopy-sheaf'' in \cite{kobe} (cf also \cite{flexible}
which predates \cite{kobe} but which was made available much later),
however that
was not the first time that such objects were encountered---the condition of
being a homotopy sheaf is the essential part of the condition of being a fibrant
(or ``flasque'') simplicial presheaf \cite{Brown} \cite{Jardine} \cite{Joyal}.

Suppose $\Cc$ is some type of category-like object (such as an $n$-category
or $\infty$-category or other such thing).  Suppose that we have a notion of
{\em inverse limit} of a family of objects of $\Cc$ indexed by a category $\Bb$.
If we call the family $F: \Bb ^o\rightarrow \Cc$ (contravariant on $\Bb$,
for our
purposes) then we denote this limit---if it exists---by $\Gamma (\Bb , F)\in
\Cc$. This should be sufficiently functorial in that if we have a functor $\Bb
\rightarrow \Bb '$ and $F$ is the pullback of a family $F'$ on $\Bb '$
(denoted $F=F'|_{\Bb}$) then we
should obtain  a morphism of functoriality (i.e. an arrow in $\Cc$)
$$
\Gamma (\Bb ', F') \rightarrow \Gamma (\Bb , F),
$$
possibly only well-defined up to some type of homotopy in $\Cc$.
Similarly if $\Bb$ has a final object $b$ (initial for our functoriality
which is
contravariant) then the morphism (obtained from above for the inclusion
$\{ b\} \rightarrow \Bb$)
$$
\Gamma (\Bb , F) \rightarrow F(b)
$$
should be an ``equivalence'' in $\Cc$ (one has to know what that means).

With all this in hand (and note that we do not assume the existence of
arbitrary limits, only existence of limits indexed by categories with final
objects) we can define the notion of {\em stack over a site $\Xx$ with
coefficients in $\Cc$}.  This is to be a family $F$ of objects of $\Cc$ indexed
by $\Xx$ (i.e. a morphism $\Xx ^o\rightarrow \Cc$) which satisfies the
following property: for every object $X\in \Xx$ and every sieve $\Bb \subset
\Xx /X$ the morphism
$$
\Gamma (\Xx /X, F|_{\Xx /X}) \rightarrow \Gamma (\Bb , F|_{\Bb})
$$
is an equivalence in $\Cc$, meaning that the limit on the right exists.
(Note
that since $\Xx /X$ has a final object $X$, the morphism
$$
\Gamma (\Xx /X, F|_{\Xx /X}) \rightarrow F(X)
$$
is assumed to exist and to be an equivalence.)

Taking the inverse limit of a family of stacks will again be a family of stacks
because inverse limits should (when that notion is defined) commute with each
other. If $\Cc$ admits arbitrary (set-theoretically reasonable) inverse limits
then taking the inverse limit of a family of stacks gives again a stack.  Using
this we can define the {\em stack associated to a prestack}. A prestack is
just any family $F: \Xx ^o \rightarrow \Cc$ not necessarily satisfying the
stack condition. The {\em associated stack} is defined to be the inverse limit
of all stacks $G$ to which $F$ maps.  Of course this needs to be investigated
some more in any specific case, in order to get useful information.

When $\Cc$ is the $2$-category of categories we obtain the classical notion of
stack \cite{LMB} \cite{ArtinInventiones}.

When $\Cc$ is the $\infty$-category of simplicial sets we obtain the notion of
``homotopy sheaf'' which is equivalent in Jardine's terminology to a simplicial
presheaf which is flasque with respect to each object of the underlying site.
In particular, fibrant simplicial presheaves satisfy this condition, and the
condition is just that of being object-by-object weak equivalent to a fibrant
simplicial presheaf. The  process of going from a prestack to the associated
stack is basically the process of going from a simplicial  presheaf to an
equivalent fibrant simplicial presheaf.

The case where $\Cc$ is the $n+1$-category $nCAT'$ of
$n$-categories yields the notion of {\em $n$-stack} described above.

\bigskip

\subnumero{Localization}

\begin{parag}
\label{universal1}
{\em Universal morphisms with certain properties}

We often encounter the following situation.  Suppose $X\in nCAT'$ and suppose
and suppose $\Pp$ is a property of morphisms $X \rightarrow B$ in $nCAT'$. Then
we can look for a universal morphism $\nu: X \rightarrow U$ with property
$\Pp$.

The ``universal'' property  can be written out in terms of our construction
$\Upsilon$: it means that for any cofibration of $n$-precats $E'\hookrightarrow
E$  and any morphism (from an edge labeled $(02)$)
$$
f:\Upsilon (E)\rightarrow nCAT'
$$
with $f(0)= X$, $f(2)=B$
sending $\Upsilon (E_0)$ to a collection of morphisms having property $\Pp$,
together with an extension along $E'$ to a morphism
$$
g':\Upsilon ^2(\ast , E')\rightarrow nCAT'
$$
with $r_{01}(g')= \nu$ and $g'(2)= B$, there exists
$$
g:\Upsilon ^2(\ast , E)\rightarrow nCAT'
$$
extending $g'$ and
with $r_{01}(g)=\nu $  and $r_{02}(g)= f$.
\end{parag}

\begin{parag}
\label{universal1.2}
Suppose $\psi : A\rightarrow nCAT'$ is a functor and $\Pp$ a property of
morphisms $\psi \rightarrow B$ to objects $B\in nCAT'$. Then we can make a
similar definition of ``universal morphism'' $\nu : \psi \rightarrow U$
having property $\Pp$.

In this case, it also makes sense to ask for a morphism
$\nu : \psi \rightarrow U$ to an object of $nCAT'$, ``universal for morphisms
with property $\Pp$'' (the definition is the same as above but we don't require
$\nu$ to have property $\Pp$).  Note that this definition in the case of one
object $X$ is vacuous: the answer would just be the identity morphism
$1_X:X\rightarrow X$.
\end{parag}

\begin{parag}
\label{universal2}
To construct $\nu$ we can try to follow the argument of \S 5,
taking the full subcategory $M(\Pp )\subset \psi /nCAT'$ of objects under $\psi$
having property $\Pp$. As before we consider the subcategory
$M(\Pp )_{\alpha}$ of objects of $\alpha$-bounded cardinality, and let
$U$ be the inverse limit of the forgetful functor
$M(\Pp )_{\alpha}\rightarrow nCAT'$.

We now need to know four things:
\newline
\ref{universal2}(i) that the morphism $\psi \rightarrow U$ again has property
$\Pp$ (preservation of $\Pp$ by inverse limits);
\newline
\ref{universal2}(ii) that there is a factorization $\psi \rightarrow U'
\rightarrow U$ with $\psi \rightarrow U'$ again having property $\Pp$ and
$\# ^{\rm pre}U'\leq \alpha$ (for $\alpha$ chosen appropriately);
\newline
\ref{universal2}(iii) that the ``telescope'' construction of
(\ref{construction2}) preserves property $\Pp$;
and
\newline
\ref{universal2}(iv)  that if $f:\psi \rightarrow \underline{Hom}(E,B)$
is a morphism which, when restricted to every object of $E_0$ gives a morphism
$\psi \rightarrow B$ with property $\Pp$, then $f$ has property $\Pp$
(this is so that a criterion analogue to \ref{easydirect} applies).
\end{parag}

\begin{conjecture}
\label{universal3}
If we know these four things then the argument of \S 5 works to
construct a universal $\nu : \psi \rightarrow T$ with property $\Pp$.
\end{conjecture}

\begin{parag}
\label{localization}
{\em Localization:}

If $X$ is an $n$-category then we denote $Fl^i(X)$ the set of $i$-morphisms,
which is the same as $X_{1,\ldots , 1}$.  Suppose we are given a collection of
subsets $S = \{ S^i \subset Fl ^i(X)\}$. Then we can define $S^{-1}X$ to be the
universal $n$-category with map $X\rightarrow S^{-1}X$ sending the elements of
$S^i$ to $i$-morphisms in $S^{-1}X$ which are invertible up to equivalence
(i.e. morphisms which are invertible in $\tau _{\leq i}(S^{-1}X)$).
To construct $S^{-1}X$, let $\Pp$ be the property of a map $X\rightarrow B$ that
the arrows in $S_i$ become invertible in $B$.   One has to verify the properties
\ref{universal2}(i)--\ref{universal2}(iv), and then apply Conjecture
\ref{universal3}. To verify the properties (i)--(iv) use
Theorem \ref{resttoIff}.

This is the $n$-categorical analogue of \cite{GabrielZisman}.
\end{parag}

{\em Caution:} If $A$ is an $m$-category considered as an $n$-category then
$S^{-1}A$ may not be an $m$-category. In particular, note that by taking the
group completion (see below) of $1$-categories one gets all homotopy types of
$n$-groupoids. (This fact, which seems to be due to Quillen, was discussed at
length in \cite{Grothendieck}...).

\begin{parag}
\label{gc}
{\em Group completion:}

The theory of $n$-categories which are not groupoids actually has a long history
in homotopy theory, in the form of the study of topological monoids.  In Adams'
book \cite{Adams} the chapter after the one on loop-space machinery, concerns
the notion of ``group completion'', namely how to go from a topological monoid
to a homotopy-theoretic group ($H$-space).  This is a special example of going
from an $n$-category to an $n$-groupoid by ``formally inverting all arrows''.

Taking $S$ to be all of the arrows in a fibrant $n$-category $X$, the
localization $S^{-1}X$ is the {\em group completion of $X$} denoted $X^{\rm
gc}$. It is the universal $n$-groupoid to which $X$ maps. This may also be
constructed by a topological approach (which has the merit of not depending on
Conjecture \ref{universal3}), as
$$
X^{\rm gc} = \Pi _n (| X | ),
$$
using Tamsamani's realization $|X|$ and Poincar\'e $n$-groupoid
$\Pi _n$ constructions \cite{Tamsamani}.
\end{parag}

As an example, \ref{resttoIff} allows us to describe the group completion of
$I$ which is contractible, as one might expect.

\begin{corollary}
The morphism $I\rightarrow \overline{I}$ is the group completion in the context
of $n$-categories.
\end{corollary}
{\em Proof:}
This follows immediately from Theorem \ref{resttoIff}.
\eop

\begin{lemma}
\label{gcCommutesWithCoprod}
Group completion commutes with coproduct. More precisely,
suppose $B\leftarrow A \rightarrow C$ are morphisms of $n$-precats. Then the
morphism
$$
(B\cup ^A C)^{\rm gp} \rightarrow B^{\rm gp} \cup ^{A^{\rm gp}} C^{\rm gp}
$$
is an equivalence.
\end{lemma}
{\em Proof:}
This can be seen directly from the topological definition
$X^{\rm gc} = \Pi _n (| X | )$ using the results of \cite{nCAT} \S 9.
\eop

\begin{parag}
\label{interiorg}
{\em Interior groupoid}
We can do a similar type of definition as \ref{universal1} for universal maps
from $B$ to $X$ having certain properties. Applying this again to the property
that all $i$-morphisms become invertible, we get the following definition.
If $X$ is a fibrant $n$-category then its {\em interior groupoid}
$X^{\rm int}$ is the universal map $X^{\rm int}\rightarrow X$ for this
property. It is an $n$-groupoid, and may be seen as the ``largest $n$-groupoid
inside $X$''.

Without refering to conjectures, we can construct $X^{\rm int}\subset X$
explicitly as follows. Assume that $X$ is an $n$-category. First we define
$X^{k\rm -int}\subset X$ with the same objects as $X$, by setting
$$
X^{1\rm -int} _{p/}(x_0,\ldots , x_p):=
X_{p/}(x_0,\ldots , x_p)^{\rm int}
$$
(note that we use inductively the definition of $Y^{\rm int}\subset Y$
for $n-1$-categories as well as the fact that this construction takes
equivalences to equivalences).

Now let
$$
X^{\rm int}_{1/}(x,y)\subset X^{1\rm -int}_{1/}(x,y)
$$
be the full sub-$(n-1)$-category of objects corresponding to morphisms which
are invertible up to equivalence.  Let
$X^{\rm int}_{p/}(x_0,\ldots , x_p)$ be the full sub-$(n-1)$-category
of $X^{1\rm -int}_{p/}(x_0,\ldots , x_p)$ consisting of objects which
project to elements of $X^{\rm int}_{1/}(x_{i-1}, x_i)$ on the principal edges.

Another way of saying this is to note that  there is a morphism
$$
X^{1\rm -int}\rightarrow \tau _{\leq 1}(X)
$$
(cf the notation of \ref{anotherapproach}).  Then
define the ``interior $1$-groupoid'' of the $1$-category $\tau _{\leq 1}(X)$
to be the subcategory consisting only of invertible morphisms, and set
$X^{\rm int}$ to be the fiber product of
$X^{1\rm -int}\rightarrow \tau _{\leq 1}(X)$ and
interior $1$-groupoid of $\tau _{\leq 1}(X)$, over
$\tau _{\leq 1}(X)$.
\end{parag}

\bigskip

\subsubnumero{$k$-groupic completion and interior}

More generally we say that an $n$-category $B$ is {\em $k$-groupic}
for $0\leq k \leq n$ if the $n-k$-categories $B_{m_1,\ldots , m_k/}$
are groupoids.  In other words this says that the $n-k$-category whose objects
are the $k$-morphisms of $B$ should be an $n-k$-groupoid.  Note that being
$O$-groupic means that $B$ is an $n$-groupoid, and the condition of being
$n$-groupic is void of content.

We can define the {\em $k$-groupic completion} $X^{k{\rm -gp}}$ as the universal
$k$-groupic $n$-category to which $X$ maps. We can define the {\em
$k$-groupic interior} $X^{k{\rm -int}}\subset X$ to be the universal $k$-groupic
$n$-category mapping to $X$.  For $k=0$ these reduce to the group
completion and interior groupoid. For the $k$-groupic interior, we have the
following formula whenever $k\geq 1$:
$$
X^{k{\rm -int}}_{p/}(x_0,\ldots , x_p)=
X_{p/}(x_0,\ldots , x_p)^{(k-1){\rm -int}},
$$
which gives an inductive construction.

\bigskip

\subnumero{Direct images and realizations}

Suppose $F:A\rightarrow B$ is a morphism of $n+1$-categories and suppose
$\varphi
: A \rightarrow nCAT'$ is a family of $n$-categories over $A$. Then
we can look for a universal family $\psi : B\rightarrow nCAT'$ together with
morphism $\varphi \rightarrow F^{\ast}(\psi )$.  If it exists, we call
$\psi$ the
{\em direct image} and denote it by $F_{\ast}(\varphi )$.

\begin{conjecture}
\label{directimage}
The direct image $F_{\ast}(\varphi )$ always exists, and is essentially unique.
\end{conjecture}

Again, the argument of \S 5 should work to give the construction
of $F_{\ast}(\varphi )$, with several things to verify analogous to
\ref{universal2}(i-iv).

\begin{parag}
{\em Caution:} the notations ``direct image'' $F_{\ast}$ and ``inverse image''
$F^{\ast}$
are switched from the usual notations for functoriality for ``morphisms of
sites''.
\end{parag}

\begin{parag}
\label{realization1}
{\em Realization:}

Suppose $A$ is an $n+1$-category and suppose
$$
\varphi : A \rightarrow nCAT'
$$
is a family of $n$-categories, and
$$
\psi : A^o \rightarrow nCAT'
$$
is a contravariant family of $n$-categories. Then we define the {\em
realization} of this pair, denoted $\langle \varphi , \psi \rangle$, as follows.
The arrow family for $A$ corresponds to a morphism
$$
\alpha : A \rightarrow \underline{Hom}(A^o, nCAT').
$$
The direct image $\alpha _{\ast}(\varphi )$ is therefore a morphism
$$
\alpha _{\ast}(\varphi )\underline{Hom}(A^o, nCAT') \rightarrow nCAT'.
$$
Put
$$
\langle \varphi , \psi \rangle := \alpha _{\ast}(\varphi )(\psi ).
$$
\end{parag}

\begin{parag}
\label{realization2}
An example of this is when $A=\Xx$ is a site, and when $\varphi$ and $\psi$ are
families of $n$-groupoids. Then
$\langle \varphi , \psi \rangle $ is an $n$-groupoid, and we conjecture that
it corresponds to the topological space given as realization of the two
functors as defined in \cite{realization}.
\end{parag}

\begin{parag}
\label{realization3}
In the main example of
\cite{realization} one took $\Xx$ to be the site of schemes over $Spec (\cc )$
and one took $\varphi$ to be the functor associating to each scheme the
$n$-truncation of the homotopy type of the underlying topological space. Then
for any presheaf $\psi$ of $n$-truncated topological spaces one obtained the
``topological realization'' of $\psi$.
\end{parag}

\begin{parag}
\label{realization4}
One can do the operation of \ref{realization1} in the other order, using
the arrow family considered as a morphism
$$
\beta : A^o \rightarrow \underline{Hom}(A, nCAT')
$$
and looking at $\beta _{\ast}(\psi ) (\varphi )$.
\newline
{\em Conjecture---}that these two ways of defining
$\langle \varphi , \psi \rangle$ give the same answer.
\end{parag}

\begin{parag}
\label{triplecombo}
The above construction is a special case of the more general phenomenon which
we call ``triple combination''. Suppose $A$ and $B$ are $(n+1)$-categories
and suppose that we have functors
$$
F: A\rightarrow nCAT',
$$
$$
G: B\rightarrow nCAT',
$$
and
$$
H: A\times B\rightarrow nCAT'.
$$
Then we can consider $H$ as a functor
$$
H:A\rightarrow \underline{Hom}(B, nCAT')
$$
and define
$$
H(F,G):= H_{\ast}(F)(G).
$$
As above, one conjectures that $H(F,G)= H^{\sigma}(G,F)$ (applying the symmetry
$\sigma : A\times B\cong B \times A$). The previous construction
is just
$$
\langle \varphi , \psi \rangle = Arr (A)(\varphi , \psi ).
$$

The same definition of triple combination works for functors $F,G,H$ in any
fibrant $n$-category $C$ which admits limits as does $nCAT'$.
\end{parag}

\bigskip

\subnumero{Relative Malcev completion}

An example which gets more to the point of my motivation for doing all of this
type of thing is the following generalization of relative Malcev completion
\cite{Hain} to higher homotopy.

\begin{parag}
\label{malcev1}
Fix a $\qq$-algebraic group $G$. Fix an
$n$-groupoid $X$ with base-object $x$ (which is the same thing as an
$n$-truncated pointed homotopy type). Fix a representation $\rho :\pi
_1(X,x)\rightarrow G$. Let $\Cc$ be the $n+1$-category of quadruples  $(R,r, p,
f)$ where $R$ is a connected $n$-groupoid, $r$ is an object, $p:
R\rightarrow BG$
is a morphism sending $r$ to the base-object $o$, and
$f: X\rightarrow R$ is a morphism sending $x$ to $r$ such that the induced
morphism $\pi _1(X,x) \rightarrow G$ is equal to $\rho$. Let $\Cc ^{\rm uni}$
denote the subset of objects satisfying the following properties: that $\pi
_1(R)$ is a $\qq$-algebraic group and $p: \pi _1(R)\rightarrow G$ is a
surjection with unipotent kernel; and that $\pi _1(R)$ acts algebraically on the
higher homotopy groups $\pi _i(R)$ which are themselves assumed to be finite
dimensional $\qq$-vector spaces.
\end{parag}

\begin{parag}
\label{malcev2}
Inverse limits
exist in $\Cc$. To see this, note that $\Cc$ is an $n+1$-category
of morphisms $V\rightarrow nCAT'$ where $V$ is the category with objects
$v_R$, $v_r$, $v_{BG}$, $v_X$ and morphisms $v_r\rightarrow v_R$,
$v_R\rightarrow v_{BG}$, $v_X\rightarrow v_R$, $v_r\rightarrow v_X$. The
$n+1$-category $\Cc$ is the subcategory of morphisms $V\rightarrow nCAT'$ which
send $v_r$ to $\ast$, send $v_{BG}$ to $BG$ and send $v_X$ to $X$, and which
send the maps $v_r\rightarrow v_X$ to the basepoint $\ast \rightarrow X$,
similarly for the map $v_r\rightarrow v_{BG}$, and which send $v_X\rightarrow
v_{BG}$ to the map induced by $\rho$.  Our Theorem \ref{inverse} as well as
\ref{variation} and Lemma \ref{fiprod}
imply that $\Cc$ admits inverse limits.  Of course $\Cc ^{\rm uni}$ is not
closed under inverse limits. However we can still take the inverse limit in
$\Cc$ of all the objects in $\Cc ^{\rm uni}$. We call this the {\em relative
Malcev completion of the homotopy type of $X$ at $\rho$}, and denote it by
$Malc(X, \rho)$ (technically this is the notation for the underlying
$n$-groupoid which is the inverse limit of the $R$'s).
\end{parag}

\begin{parag}
\label{malcev3}
We have, for example, that $\pi _1(Malc(X, \rho ), \ast )$ is equal to the
relative Malcev completion of the fundamental group $\pi _1(X)$ at $\rho$.
For this statement we fall back into the realm of $1$-categories, where
our Malcev completion coincides with the usual notion \cite{Hain}.
\end{parag}

\begin{parag}
\label{malcev4}
We can do the same thing with stacks.  For a field $k$ (of
characteristic zero, say) an algebraic group $G$ over $k$ and a
representation $\rho : \pi _1(X,x)\rightarrow G$, let $\Cc (X, \rho )/k$ be the
$n+1$-category of quadruples $(R,r, p,f)$ where $R$ is a connected $n$-stack
of groupoids on
$Sch/k$, $r$ is a basepoint,  $p: R\rightarrow BG$, and $f:
\underline{X}\rightarrow R$ are as above. Here $\underline{X}$ is the constant
stack with values $X$. Let $\Cc ^{\rm uni}(X,\rho )/k$ be the subcategory of
objects such that $\pi _1(R)$ is an algebraic group surjecting onto $G$ and
where the  $\pi _i(R)$ are linear finite dimensional representations of $\pi
_1(R)$. Again inverse limits will exist in $\Cc (X,\rho )/k$ and we can take
the inverse limit here of the objects of $\Cc ^{\rm  uni}(X, \rho )/k$.
Call this $Malc(X, \rho )/k$.

Note that $Malc(X,\rho )/\qq )$ is an $n$-stack on $Sch /\qq $ whose
$n$-groupoid of global sections is $Malc(X,\rho )$.
\end{parag}

\begin{parag}
\label{malcev5}
Suppose $X$ is a variety and let $X_B$ be the $n$-groupoid truncation of the
homotopy type of $X^{\rm top}$. Fix a representation $\rho$. Then we obtain the
``Betti'' Malcev completion $Malc(X_B, \rho )/\cc $.  On the other hand
suppose $P$ is the principal $G$-bundle with integrable connection
(with regular singularities at infinity) corresponding
to $\rho$, then  we can define in a similar way $Malc (X_{DR}, P)/\cc $.
The GAGA results imply that these two are naturally equivalent:
$$
Malc(X_B, \rho )/\cc \cong Malc(X_{DR}, P) /\cc .
$$
SImilarly we can define, for a principal Higgs bundle $Q$ with vanishing Chern
classes, $Malc (X_{Dol}, Q)/\cc$, and (in the case $X$ smooth projective)
if $Q$ corresponds to $\rho$ then
$$
Malc(X_B, \rho )/\cc \cong Malc(X_{Dol}, Q) /\cc .
$$
Finally, suppose $\rho$ is an $\rr$-variation of Hodge structure and $Q$ the
corresponding system of Hodge bundles.  Then $\cc ^{\ast }$ acts on
$Malc(X_{Dol}, Q) /\cc $ giving rise to a ``weight filtration'' and ``Hodge
filtration''. We conjecture that these (together with the $\rr$-rational
structure $Malc(X_B,\rho )/\rr $) define a ``mixed Hodge structure''
on $Malc(X_B, \rho)/\rr$.  (One has to give this definition, specially in view
of the infinite size of $Malc(X_B, \rho)/\rr$).
\end{parag}

More generally we have the following conjecture.

\begin{conjecture}
\label{mts}
Suppose $\rho$ is a reductive representation of the fundamental group of a
projective variety $X$ (we assume it is reductive when restricted to the
fundamental group of the normalization).  Then the relative Malcev
completion of the higher homotopy type $Malc(X_B, \rho)/\cc$ defined above
carries a natural mixed twistor structure (cf \cite{twistor}).
\end{conjecture}

There should also be a statement for quasiprojective varieties, but in this
case one probably needs some additional hypotheses on the behavior of $\rho$ at
infinity.


\bigskip







\begin{thebibliography}{MM2}

\bibitem{Adams}
J. Adams. {\em Infinite Loop Spaces}, Princeton University Press {\em Annals
of Math. Studies} {\bf 90} (1978).

\bibitem{ArtinInventiones}
M. Artin.  Versal deformations and algebraic stacks, {\em Inventiones Math.}
{\bf 27} (1974), 165-189.


\bibitem{BaezDolanLetter}
J. Baez, J. Dolan. $n$-Categories, sketch of a
definition. Letter to R. Street, 29 Nov. and 3 Dec. 1995, available at
{\tt http://math.ucr.edu/home/baez/ncat.def.html}


\bibitem{BaezDolan}
J. Baez, J. Dolan. Higher-dimensional algebra and
topological quantum field theory. {\em  Jour. Math. Phys} {\bf  36} (1995),
6073-6105 (preprint dating from q-alg 95-03).

\bibitem{BaezDolanIII}
J. Baez, J. Dolan. Higher dimensional algebra III: $n$-categories and the
algebra of opetopes. Preprint available on q-alg (97-02).

\bibitem{Batanin}
M. Batanin. On the definition of weak $\omega$-category. Macquarie mathematics
report number 96/207, Macquarie University, NSW Australia.


\bibitem{Batanin2}
M. Batanin.
Monoidal globular categories as a natural environment for the theory of weak
$n$-categories. Preprint, April 1997.

\bibitem{Benabou}
J. B\'enabou. {\em Introduction to Bicategories}, Lect. Notes in Math. {\bf
47}, Springer-Verlag (1967).

\bibitem{Bousfield-Kan}
A. Bousfield, D. Kan. {\em Homotopy limits, completions and localizations. }
Springer {\em Lecture Notes in Mathematics} {\bf 304} (1972).

\bibitem{Breen}
L. Breen. On the classification of $2$-gerbs and $2$-stacks. {\em Ast\'erisque}
{\bf 225}, Soc. Math. de France (1994).

\bibitem{Brown}
K. Brown. Abstract homotopy theory and generalized sheaf cohomology. {\em
Trans. A.M.S.} {\bf 186} (1973), 419-458.

\bibitem{BrownGersten}
K. Brown, S. Gersten. {\em Algebraic $K$-theory as generalize sheaf cohomology}
Springer {\em Lecture Notes in Math.} {\bf 341} (1973), 266-292.

\bibitem{DeligneMumford}
P. Deligne, D. Mumford. On the irreducibility of the space of curves of
a given genus. {\em Publ. Math. I.H.E.S.}  {\bf 36} (1969), 75-109.

\bibitem{GabrielZisman}
P. Gabriel, M. Zisman. {\em Calculus of fractions and homotopy theory.}
Springer, New York (1967).

\bibitem{Giraud}
J. Giraud. {\em Cohomologie nonab\'elienne}, Grundelehren der Wissenschaften
in Einzeldarstellung {\bf 179} Springer-Verlag (1971).

\bibitem{Gordon-Power-Street}
R. Gordon, A.J. Power, R. Street.  Coherence for tricategories {\em  Memoirs
A.M.S.} {\bf 117} (1995), 558 ff.

\bibitem{Grothendieck}
A. Grothendieck.  {\em Pursuing Stacks}, unpublished manuscript.


\bibitem{Hain}
R. Hain. Completions of mapping class groups and the cycle $C-C^{-}$.

\bibitem{Jardine}
J.F. Jardine.  Simplicial presheaves, {\em J. Pure and Appl. Algebra} {\bf 47}
(1987), 35-87.

\bibitem{Johnson}
M. Johnson. The combinatorics of $n$-categorical pasting. {\em J. Pure and
Appl. Algebra} {\bf 62} (1989), 211-225.

\bibitem{Joyal}
A. Joyal. Letter to A. Grothendieck (refered to in Jardine's paper).

\bibitem{LMB}
G. Laumon, L. Moret-Bailly. Champs alg\'ebriques. Preprint, Orsay {\bf 42}
(1992).

\bibitem{MacLane}
S. MacLane. Categories for the working mathematician. Springer (1971).

\bibitem{May}
J. P. May. {\em Simplicial objects in algebraic topology.} Van Nostrand (1967).

\bibitem{Quillen}
D. Quillen. {\em Homotopical algebra} Springer {\em L.N.M.} {\bf 43} (1967).

\bibitem{QuillenAnnals}
D. Quillen. Rational Homotopy Theory. {\em Ann. Math.} {\bf 90} (1969), 205-295.

\bibitem{Segal}
G. Segal. Homotopy everything $H$-spaces. Preprint.

\bibitem{kobe}
C. Simpson. Homotopy over the complex numbers and generalized de Rham
cohomology. {\em Moduli of Vector Bundles}, M. Maruyama (Ed.) {\em Lecture
Notes in Pure and Applied Math.} {\bf 179}, Marcel Dekker (1996), 229-263.

\bibitem{flexible}
C. Simpson. Flexible sheaves. Preprint available on q-alg (96-08).

\bibitem{realization}
C. Simpson. The topological realization of a simplicial presheaf. Preprint,
available on q-alg 96-09.

\bibitem{geometricN}
C. Simpson. Algebraic (geometric) $n$-stacks. Preprint,
available on alg-geom 96-09.

\bibitem{nCAT}
C. Simpson. A closed model structure for $n$-categories, internal $Hom$,
$n$-stacks and generalized Seifert-Van Kampen. Preprint, available on
alg-geom 97-04.

\bibitem{twistor}
C. Simpson. Mixed twistor structures. Preprint, available on alg-geom
97-05.

\bibitem{Stasheff}
J. Stasheff. Homotopy associativity of $H$-spaces. {\em Trans. Amer. Math.
Soc.} {\bf 108} (1963), 275-312.

\bibitem{Street}
R. Street. The algebra of oriented simplexes. {\em Jour. Pure and Appl.
Algebra} {\bf 49} (1987), 283-335.

\bibitem{Tamsamani}
Z. Tamsamani.  Sur des notions de $n$-categorie et $n$-groupoide non-stricte
via des ensembles multi-simpliciaux. Thesis, Universit\'e Paul Sabatier,
Toulouse (1996) available on alg-geom (95-12 and 96-07).

\bibitem{Tanre}
D. Tanre. {\em  Homotopie Rationnelle: mod\`eles de Chen, Quillen, Sullivan.}
Springer {\em Lecture Notes in Mathematics} {\bf 1025}  (1983).


\end{thebibliography}
\end{document}